\documentclass[a4paper,12pt]{amsart}
\usepackage{amsmath,amsthm,amssymb}
\usepackage[dvipdfmx]{graphicx}
\usepackage[dvipdfmx]{hyperref}
\usepackage[all]{xy}
\usepackage{graphicx}
\usepackage{here}
\usepackage{color}
\usepackage[top=30truemm,bottom=30truemm,left=25truemm,right=25truemm]{geometry}
\usepackage{ytableau}

\newtheorem{theo}{Theorem}[subsection]
\newtheorem{defi}[theo]{Definition}
\newtheorem{lem}[theo]{Lemma} 
\newtheorem{rem}[theo]{Remark}
\newtheorem{prop}[theo]{Proposition}
\newtheorem{cor}[theo]{Corollary}

\newtheorem{ex}[theo]{Example}

\newtheorem{Thm}{Theorem}

\newcommand{\nc}{\newcommand}
\nc{\on}{\operatorname}
\nc{\C}{\mathbb{C}}
\nc{\R}{\mathbb{R}}
\nc{\Q}{\mathbb{Q}}
\nc{\Z}{\mathbb{Z}}
\nc{\N}{\mathbb{N}}

\nc{\bbH}{\mathbb{H}}
\nc{\bbK}{\mathbb{K}}

\nc{\bfa}{\mathbf{a}}
\nc{\bfA}{\mathbf{A}}
\nc{\bfb}{\mathbf{b}}
\nc{\bfi}{\mathbf{i}}
\nc{\bfK}{\mathbf{K}}
\nc{\bfk}{\mathbf{k}}
\nc{\bfone}{\boldsymbol 1}
\nc{\bfeta}{\boldsymbol \eta}
\nc{\bfkappa}{\boldsymbol \kappa}
\nc{\bfsigma}{\boldsymbol \sigma}
\nc{\bfvarsigma}{\boldsymbol \varsigma}
\nc{\bfzeta}{\boldsymbol \zeta}

\nc{\A}{\mathbf{A}}
\nc{\U}{\mathbf{U}}

\nc{\clB}{\mathcal{B}}
\nc{\clC}{\mathcal{C}}
\nc{\clF}{\mathcal{F}}
\nc{\clL}{\mathcal{L}}
\nc{\clM}{\mathcal{M}}
\nc{\clO}{\mathcal{O}}
\nc{\clU}{\mathcal{U}}
\nc{\clUi}{\clU^{\imath}}
\nc{\clX}{\mathcal{X}}
\nc{\clXs}{\clX_{\mathrm{s}}}
\nc{\clY}{\mathcal{Y}}
\nc{\clYs}{\clY_{\mathrm{s}}}
\nc{\clW}{\mathcal{W}}
\nc{\clZ}{\mathcal{Z}}

\nc{\fra}{\mathfrak{a}}
\nc{\frh}{\mathfrak{h}}
\nc{\g}{\mathfrak{g}}
\nc{\frgl}{\mathfrak{gl}}
\nc{\frk}{\mathfrak{k}}
\nc{\fram}{\mathfrak{m}}
\nc{\frn}{\mathfrak{n}}
\nc{\frp}{\mathfrak{p}}
\nc{\frs}{\mathfrak{s}}
\nc{\frt}{\mathfrak{t}}
\nc{\frsl}{\mathfrak{sl}}
\nc{\frso}{\mathfrak{so}}
\nc{\frsp}{\mathfrak{sp}}
\nc{\frsu}{\mathfrak{su}}
\nc{\fru}{\mathfrak{u}}
\nc{\frz}{\mathfrak{z}}

\nc{\fin}{\mathrm{fin}}

\nc{\inv}{^{-1}}
\nc{\qu}{\quad}

\nc{\qqu}{\qquad}
\nc{\la}{\langle}
\nc{\ra}{\rangle}

\nc{\Ker}{\on{Ker}}
\nc{\im}{\on{Im}}
\nc{\Hom}{\on{Hom}}
\nc{\End}{\on{End}}
\nc{\Span}{\on{Span}}
\nc{\id}{\on{id}}
\nc{\Aut}{\on{Aut}}
\nc{\ad}{\on{ad}}
\nc{\sgn}{\on{sgn}}
\nc{\tot}{\on{tot}}
\nc{\rk}{\on{rk}}
\nc{\rank}{\on{rank}}
\nc{\Wt}{\on{Wt}}
\nc{\diag}{\on{diag}}
\nc{\Mat}{\on{Mat}}
\nc{\tr}{\on{tr}}
\nc{\Diag}{\on{Diag}}
\nc{\GL}{\on{GL}}
\nc{\SO}{\on{SO}}
\nc{\Sp}{\on{Sp}}
\nc{\gr}{\on{gr}}
\nc{\Ind}{\on{Ind}}
\nc{\Res}{\on{Res}}
\nc{\wt}{\on{wt}}
\nc{\lt}{\on{lt}}
\nc{\lc}{\on{lc}}
\nc{\Br}{\on{Br}}
\nc{\norm}{\on{norm}}
\nc{\amp}{\on{amp}}
\nc{\Int}{\on{int}}
\nc{\Inv}{\on{inv}}
\nc{\pr}{\on{pr}}

\nc{\ev}{\mathrm{ev}}
\nc{\odd}{\mathrm{odd}}
\nc{\ch}{\mathrm{ch}}

\nc{\ol}{\overline}
\nc{\ul}{\underline}
\nc{\hf}{\frac{1}{2}}
\nc{\vphi}{\varphi}
\nc{\vrho}{\varrho}
\nc{\vpi}{\varpi}
\nc{\vep}{\varepsilon}
\nc{\eps}{\epsilon}
\nc{\lm}{\lambda}
\nc{\til}{\widetilde}

\nc{\IF}{\text{ if }}
\nc{\AND}{\text{ and }}
\nc{\OR}{\text{ or }}
\nc{\OW}{\text{ otherwise}}
\nc{\lowerterms}{\text{(lower terms)}}
\nc{\higherterms}{\text{(higher terms)}}
\nc{\lex}{\text{lex}}
\nc{\sesi}{\text{ss}}
\nc{\ST}{\text{ such that }}
\nc{\Forsome}{\text{ for some }}
\nc{\Forall}{\text{ for all }}

\nc{\Atil}{\widetilde{A}}
\nc{\Btil}{\widetilde{B}}
\nc{\Etil}{\widetilde{E}}
\nc{\Ftil}{\widetilde{F}}
\nc{\Itil}{\widetilde{I}}
\nc{\Xtil}{\widetilde{X}}
\nc{\Ytil}{\widetilde{Y}}

\nc{\Ui}{\U^{\imath}}
\nc{\Uidot}{\dot{\U}^\imath}
\nc{\UidotA}{\dot{\U}^\imath_{\bfA}}
\nc{\taui}{\tau^{\imath}}
\nc{\psii}{\psi^{\imath}}
\nc{\wti}{\wt^\imath}

\nc{\Udot}{\dot{\U}}
\nc{\UdotA}{\Udot_{\bfA}}

\nc{\plim}[1][]{\mathop{\varprojlim}\limits_{#1}}
\nc{\ilim}[1][]{\mathop{\varinjlim}\limits_{#1}}

\nc{\TBA}{{\large {\bf \textcolor{red}{To Appear}}}}
\nc{\alert}{\textcolor{red}}

\title{Based modules over the $\imath$quantum group of type AI}
\author[H. Watanabe]{Hideya Watanabe}
\address{(H. Watanabe) Research Institute for Mathematical Sciences, Kyoto University, Kyoto 606-8502, Japan}
\email{hideya@kurims.kyoto-u.ac.jp}
\date{\today}

\begin{document}
\maketitle

\begin{abstract}
This paper studies classical weight modules over the $\imath$quantum group $\Ui$ of type AI. We introduce the notion of based $\Ui$-modules by generalizing the notion of based modules over quantum groups (quantized enveloping algebras). We prove that each finite-dimensional irreducible classical weight $\Ui$-module with integer highest weight is a based $\Ui$-module. As a byproduct, a new combinatorial formula for the branching rule from $\frsl_n$ to $\frso_n$ is obtained.
\end{abstract}

\section{Introduction}
\subsection*{Based modules over quantum groups}
In representation theory of quantum groups (quantized enveloping algebras), there is a notion of based modules. Let $\U = U_q(\g)$ denote the quantum group over the field $\bfK := \C(q)$ of rational functions in one variable $q$ associated with a complex semisimple Lie algebra $\g$. A based $\U$-module is a $\U$-module equipped with two distinguished structures; a crystal basis and a canonical basis, which is also known as a global crystal basis.

Crystal basis theory, which was discovered by Kashiwara \cite{Ka90} and Lusztig \cite{L90}, is a powerful tool in combinatorial representation theory. The crystal basis $\clB_M$ of a based $\U$-module $M$ is a local basis in the sense that it forms a basis of $M$ at $q = \infty$ (in the literature, the variable $q$ in this paper corresponds to $q\inv$). To be more specific, $\clB_M$ is a $\C$-basis of $\ol{\clL}_M := \clL_M/q\inv \clL_M$, where $\clL_M$ is a free submodule of $M$ over the ring $\bfK_\infty := \C[q\inv]_{(q\inv)}$ of rational functions regular at $q=\infty$, called the crystal lattice. Although $\clB_M$ is not a genuine $\bfK$-basis of $M$, it has a lot of information about the $\U$-module structure of $M$.

Canonical basis theory was initiated by Lusztig \cite{L90}, and then Kashiwara \cite{Ka91} found a different approach. The canonical basis of a based $\U$-module $M$ is a globalization of the crystal basis in the sense that it tends to the crystal basis as $q$ goes to $\infty$.

One of the most distinguished properties of based $\U$-modules is cellularity (\cite[Chapter 27]{L10}), which gives a canonical stratification of based $\mathbf{U}$-modules.
This generalizes a well-known property of the Hecke algebra of type $A$ with the Kazhdan-Lusztig basis \cite{KL79}.

\subsection*{$\imath$Quantum groups}
Recall that the quantum group $\U$ is a $q$-deformation of the universal enveloping algebra $U(\g)$ of $\g$. Namely, $\U$ tends to $U(\g)$ as $q$ goes to $1$. There is another family of $q$-deformations of the universal enveloping algebras of complex Lie algebras, called the $\imath$quantum groups.

The $\imath$quantum groups appear in the theory of quantum symmetric pairs initiated by Letzter \cite{Le99} unifying earlier constructions by Koornwinder \cite{Ko90}, Gavrilik and Klimyk \cite{GK91}, Noumi \cite{N96}, and others. Kolb \cite{Ko14} generalized Letzter's theory to symmetrizable Kac-Moody algebras. In the theory of quantum symmetric pairs, we consider a complex semisimple Lie algebra $\g$ and an involutive Lie algebra automorphism $\theta$ on $\g$. Set $\frk := \g^\theta = \{ x \in \g \mid \theta(x) = x \}$ to be the fixed point subalgebra. A quantum symmetric pair associated with $(\g,\frk)$ is a pair $(\U,\Ui)$ consisting of the quantum group $\mathbf{U} = U_q(\mathfrak{g})$ and its coideal subalgebra $\mathbf{U}^\imath$ which tends to $U(\mathfrak{k})$ as $q$ goes to $1$.
The $\Ui$ itself is referred to as an $\imath$quantum group.

The $\imath$quantum groups $\Ui$ (with the embeddings $\Ui \hookrightarrow \U$) are thought of as generalizations of the quantum groups in the sense that a quantum group $\U$ (with the comultiplication map $\Delta : \U \hookrightarrow \U \otimes \U \simeq U_q(\g \oplus \g)$) is an $\imath$quantum group associated with $(\g \oplus \g, \g)$. Such an $\imath$quantum group is said to be of diagonal type.

Based on this viewpoint, many results in the theory of quantum groups have been generalized to the $\imath$quantum groups setting, e.g., the bar-involution \cite{BK15, BW18b}, the quasi-$K$-matrix \cite{BW18a, BW18b}, the universal $K$-matrix \cite{AV20, BK19, DCM20, RV20}, quantum Schur and Howe dualities \cite{BW18a, ES18, ST19}, and the $\imath$canonical basis \cite{BW18a, BW18b}.
However, the theory of based $\mathbf{U}$-modules has not been generalized to the $\imath$quantum groups setting.

\subsection*{Results}
The aim of this paper is to define the notion of based $\Ui$-modules, and study their properties for the $\imath$quantum group of type AI, i.e., the $\imath$quantum group associated with $(\g,\frk) \simeq (\frsl_n,\frso_n)$ (see \cite{Wa17, Wa21} for related works for the $\imath$quantum group of quasi-split type AIII).
The symmetric pair $(\mathfrak{sl}_n, \mathfrak{so}_n)$ is of split type, and hence, its structure, such as defining relations, is relatively simple.
Moreover, it is one of the generalized Onsager algebras, which have background in integrable systems (see \cite{St20} and references therein for detail).

In \cite{W19}, the notion of classical weight $\Ui$-modules was introduced, and the finite-dimensional irreducible classical weight $\Ui$-modules were classified in terms of highest weight theory for several types including AI. The isoclasses of such irreducible modules are parametrized by the set $X_{\frk}^+$ of dominant integral weights for $\frk$. For each $\nu \in X_{\frk}^+$, let $V(\nu)$ denote the corresponding irreducible $\Ui$-module.
Recall that the integral weights for $\mathfrak{k} \simeq \mathfrak{so}_n$ are divided into two groups, integer weights and half integer weights.
Let $X_{\mathfrak{k}, \text{int}}$ denote the set of integer weights, and set $X_{\mathfrak{k}, \text{int}}^+ := X^+ \cap X_{\mathfrak{k}, \text{int}}$.
The first main result in this paper is the following.

\begin{Thm}[{Theorem \ref{V(nu) is a based Ui-module; general}}]
Let $\nu \in X_{\frk}^+$ be an integer weight. Then, $V(\nu)$ has a based $\mathbf{U}^\imath$-module structure.
\end{Thm}

As a byproduct, we obtain a purely combinatorial formula for the branching rule from $\U$ to $\Ui$.
Before stating the formula, let us prepare some notation.
First of all, let us identify $X_{\frk,\Int}$ with $\Z^m$, where $m$ denotes the rank of $\frk (\simeq \frso_n)$:
$$
m = \begin{cases}
\frac{n}{2} \qu & \IF \text{$n$ is even}, \\
\frac{n-1}{2} \qu & \IF \text{$n$ is odd}.
\end{cases}
$$
Let $\Etil_i,\Ftil_i$ denote the Kashiwara operators acting on crystal bases \cite{Ka90}. Given a finite-dimensional $\U$-module and its crystal basis $\clB$, for each $b \in \clB$, set
\begin{align}
\begin{split}
&\vphi_i(b) := \max\{ k \mid \Ftil_i^k b \neq 0 \}, \\
&\vep_i(b) := \max\{ k \mid \Etil_i^k b \neq 0 \}, \\
&\deg_i(b) := \begin{cases}
\vep_i(b) \qu & \IF \vphi_i(b) \text{ is even}, \\
\vep_i(b)+1 \qu & \IF \vphi_i(b) \text{ is odd},
\end{cases} \\
&\Btil_i(b) := \begin{cases}
\Etil_i b \qu & \IF \vphi_i(b) \text{ is even}, \\
\Ftil_i b \qu & \IF \vphi_i(b) \text{ is odd}.
\end{cases}
\end{split} \nonumber
\end{align}

\begin{Thm}[{Theorem \ref{Branching rule from U to Ui}}]
Let $\lm \in X^+$ and $\nu = (\nu_1,\nu_3,\ldots,\nu_{2m-1}) \in X_{\frk,\Int}^+ = \Z^m$. Let $[\lm:\nu]$ denote the multiplicity of the irreducible $\Ui$-module $V(\nu)$ of highest weight $\nu$ in the irreducible $\U$-module $V(\lm)$ of highest weight $\lm$.
\begin{enumerate}
\item Suppose $n$ is even and $\nu_{2m-1} \neq 0$. Then, we have
\begin{align}
\begin{split}
[\lm:\nu] = \hf \sharp \{ b \in \clB(\lm) \mid &\deg_{2i-1}(b) = |\nu_{2i-1}| \Forall i \in \{ 1,\ldots,m \}, \\
&\deg_{2i}(b) = 0 \Forall i \in \{ 1,\ldots,m-1 \}, \\
&\deg_{2i+1}((\Btil_{2i}\Btil_{2i-1})^{|\nu_{2i+1}|} b) = 0 \Forall i \in \{ 1,\ldots,m-1 \} \}.
\end{split} \nonumber
\end{align}
\item Suppose $n$ is even and $\nu_{2m-1} = 0$. Then, we have
\begin{align}
\begin{split}
[\lm:\nu] = \sharp \{ b \in \clB(\lm) \mid &\deg_{2i-1}(b) = \nu_{2i-1} \Forall i \in \{ 1,\ldots,m \}, \\
&\deg_{2i}(b) = 0 \Forall i \in \{ 1,\ldots,m-1 \}, \\
&\deg_{2i+1}((\Btil_{2i}\Btil_{2i-1})^{\nu_{2i+1}} b) = 0 \Forall i \in \{ 1,\ldots,m-1 \} \}.
\end{split} \nonumber
\end{align}
\item Suppose $n$ is odd. Then, we have
\begin{align}
\begin{split}
[\lm:\nu] = \sharp \{ b \in \clB(\lm) \mid &\deg_{2i-1}(b) = \nu_{2i-1} \Forall i \in \{ 1,\ldots,m \}, \\
&\deg_{2i}(b) = 0 \Forall i \in \{ 1,\ldots,m \}, \\
&\deg_{2i+1}((\Btil_{2i}\Btil_{2i-1})^{\nu_{2i+1}} b) = 0 \Forall i \in \{ 1,\ldots,m-1 \} \}.
\end{split} \nonumber
\end{align}
\end{enumerate}
\end{Thm}

Although $\frk$ is not set-theoretically the same as $\frso_n$, the Lie algebra of $n$ by $n$ skew-symmetric matrices, this formula also tells us the branching rule from $\frsl_n$ to $\frso_n$ by a similar argument to \cite[2.4]{NS05}. This classical branching problem has been studied for a long time by many people, and several (partial) answers have been provided. Among them, Jang and Kwon \cite{JK21} recently described this branching rule in terms of Littlewood-Richardson coefficients. Their formula identifies the multiplicities with the numbers of certain combinatorial objects; Littlewood-Richardson tableaux satisfying certain conditions. On the other hand, since $\clB(\lm)$ is realized as the set of semistandard tableaux of shape $\lm$, our formula identifies the multiplicities with the numbers of certain semistandard tableaux. Hence, it would be interesting to construct an explicit bijection between such Littlewood-Richardson tableaux and semistandard tableaux.

Recall that the crystal basis $\clB_M$ of a finite-dimensional based $\U$-module $M$ tells us how $M$ decomposes into irreducible ones. In particular, the crystal basis theory provides us a branching rule from $\U \otimes \U$ to $\U$. Recall also that $(\U \otimes \U, \U)$ is an $\imath$quantum group of diagonal type. Then, our formula can be seen as a generalization of the branching rule for the $\imath$quantum groups of diagonal type to that of type AI.

\subsection*{Organization}
This paper is organized as follows. In Section \ref{section: quantum groups}, we briefly review the theory of based modules over the quantum groups. In Section \ref{section: QSP}, we formulate the notion of based modules over the $\imath$quantum groups by generalizing that of based modules over the quantum groups. From Section \ref{sectoin: QSP of type AI} on, we restrict our attention to the $\imath$quantum group of type AI. For such an $\imath$quantum group, we introduce the notion of standard $X^\imath$-weight modules, which plays a key role when proving that the finite-dimensional irreducible highest weight modules of integer highest weights are based modules. In Sections \ref{section: n=2}--\ref{section: n=4}, we analyze low-rank cases. Based on results obtained there, we finally prove our main theorems in Section \ref{section: general}.

\subsection*{Acknowledgement}
The author thanks anonymous referees for helpful comments.
This work was supported by JSPS KAKENHI Grant Number JP20K14286.

\subsection*{Notation}
Throughout this paper, we use the following notation:
\begin{itemize}
\item $\Z_{\geq 0}$: the set of nonnegative integers.
\item $\Z_{\ev}$: the set of even integers.
\item $\Z_{\odd}$: the set of odd integers.
\item $\Z_{\geq 0, \ev} := \Z_{\geq 0} \cap \Z_{\ev}$.
\item $\Z_{\geq 0, \odd} := \Z_{\geq 0} \cap \Z_{\odd}$.
\item For $a \in \Z$, $p(a) := \begin{cases}
\ev \qu & \IF a \in \Z_{\ev}, \\
\odd \qu & \IF a \in \Z_{\odd}.
\end{cases}$
\item For $a \in \Z$, $q(a) := p(a-1) = \begin{cases}
\odd \qu & \IF a \in \Z_{\ev}, \\
\ev \qu & \IF a \in \Z_{\odd}.
\end{cases}$
\item For $a,b \in \Z$, $[a,b] := \{ c \in \Z \mid a \leq c \leq b \}$.
\item For $a,b \in \Z$ and $p \in \{ \ev, \odd \}$, $[a,b]_p := [a,b] \cap \Z_p$.
\item $\Z/2\Z = \{ \ol{0}, \ol{1} \}$: the abelian group of order $2$.
\item $\ol{\cdot} : \Z \rightarrow \Z/2\Z$: the quotient map.
\item For $a \in \Z/2\Z$, $p(a) := \begin{cases}
\ev \qu & \IF a = \ol{0}, \\
\odd \qu & \IF a = \ol{1}.
\end{cases}$
\end{itemize}

\section{Quantum groups}\label{section: quantum groups}
In this section, we briefly review basic results concerning finite-dimensional based modules over a quantum group.
We refer the reader to \cite{L10}.

\subsection{Hermitian inner product}
Let $\bfK$ denote the field $\C(q)$ of complex rational functions. In this subsection, we formulate the notion of $\bfK$-valued Hermitian inner product and generalize basic results about complex metric spaces.
The proofs of statements in this 
subsection are straightforward, and hence, omitted.

Set
\begin{align}
\begin{split}
&\bfK_1 := \left\{ \frac{f}{g} \in \bfK \mid f,g \in \C[q],\ g|_{q=1} \neq 0 \right\}, \\
&\bfK_\infty := \left\{ \frac{f}{g} \in \bfK \mid f,g \in \C[q\inv],\ g|_{q\inv = 0} \neq 0 \right\}.
\end{split} \nonumber
\end{align}
Let $\ev_1 : \bfK_1 \rightarrow \bfK_1/(q-1)\bfK_1 \simeq \C$ (resp., $\ev_\infty : \bfK_\infty \rightarrow \bfK_\infty/q\inv \bfK_\infty \simeq \C$) denote the evaluation map at $q = 1$ (resp., $q\inv = 0$).

\begin{defi}\normalfont
Let $c \in \bfK^\times$, and write $c = \frac{f}{g}$ with $f = \sum_{m} a_m q^{m} \in \C[q,q\inv]$, $g \in 1 + q\inv \C[q\inv]$.
\begin{enumerate}
\item The degree $\deg(c) \in \Z$ of $c$ is defined to be $\max\{ m \mid a_m \neq 0 \}$.
\item The leading coefficient $\lc(c) \in \C$ of $c$ is defined to be $a_{\deg(c)}$.
\item The leading term $\lt(c) \in \C q^{\Z}$ of $c$ is defined to be $\lc(c)q^{\deg(c)}$.
\end{enumerate}
Furthermore, we set
$$
\deg(0) := -\infty,\ \lc(0) := 0,\ \lt(0) := 0.
$$
\end{defi}

Let $z^* \in \C$ denote the complex conjugate of $z \in \C$. Extend the notion of complex conjugate to a ring automorphism on $\bfK$ by setting $q^* = q$.

\begin{defi}\label{Definition of Hermitian inner product}\normalfont
Let $V$ be a $\bfK$-vector space. A ($\bfK$-valued) Hermitian inner product on $V$ is a map $(\cdot,\cdot) : V \times V \rightarrow \bfK$ satisfying the following:
\begin{enumerate}
\item\label{Definition of Hermitian inner product 1} $(au+bv,w) = a(u,w) + b(v,w)$ for all $a,b \in \bfK$, $u,v,w \in V$.
\item\label{Definition of Hermitian inner product 2} $(v,u) = (u,v)^*$ for all $u,v \in V$.
\item\label{Definition of Hermitian inner product 3} $\lc((v,v)) \geq 0$ for all $v \in V$.
\item\label{Definition of Hermitian inner product 4} $\lc((v,v)) = 0$ implies $v = 0$.
\item\label{Definition of Hermitian inner product 5} $\deg((v,v)) \in \Z_{\ev} \sqcup \{-\infty\}$ for all $v \in V$.
\end{enumerate}
\end{defi}

\begin{defi}\normalfont
Let $V$ be a $\bfK$-vector space equipped with a Hermitian inner product. For $v \in V$, we set
$$
\deg(v) := \hf \deg((v,v)), \qu \lc(v) := \sqrt{\lc((v,v))}, \qu \lt(v) := \lc(v) q^{\deg(v)}.
$$
Here, we understand that $\hf (-\infty) = -\infty$ and $q^{-\infty} = 0$.
\end{defi}

\begin{defi}\normalfont
Let $V$ be a $\bfK$-vector space equipped with a Hermitian inner product.
\begin{enumerate}
\item $v \in V$ is said to be almost normal if $\lt(v) = 1$.
\item $u,v \in V$ are said to be orthogonal (resp., almost orthogonal) if $(u,v) = 0$ (resp., $(u,v) \in q\inv \bfK_\infty$).
\item A basis of $V$ is said to be orthogonal if it consists of vectors which are pairwise orthogonal.
\item A basis of $V$ is said to be almost orthonormal if it consists of almost normal vectors which are pairwise almost orthogonal.
\end{enumerate}
\end{defi}

\begin{rem}\normalfont
Because of conditions \eqref{Definition of Hermitian inner product 3}--\eqref{Definition of Hermitian inner product 5} in Definition \ref{Definition of Hermitian inner product}, every nonzero vector can be normalized to an almost normal vector proportional to it.
\end{rem}

The following is an easy analogue of an elementary result about complex metric spaces.

\begin{prop}\label{existence of orthogonal basis}
Let $V$ be a finite-dimensional $\bfK$-vector space equipped with a Hermitian inner product $(\cdot,\cdot)$, and $W \subseteq V$ a subspace.
\begin{enumerate}
\item\label{existence of orthogonal basis 1} $V$ possesses an orthogonal and almost orthonormal basis.
\item\label{existence of orthogonal basis 2} The restriction $(\cdot,\cdot)_W := (\cdot,\cdot)|_{W \times W}$ is a Hermitian inner product on $W$.
\item\label{existence of orthogonal basis 3} Let $B_W$ be an almost orthonormal basis of $W$. Then, there exists an almost orthonormal basis $B_V$ of $V$ which extends $B_W$.
\end{enumerate}
\end{prop}

Let $V$ and $W$ be as in Proposition \ref{existence of orthogonal basis}. In the sequel, unless otherwise stated, whenever we consider a Hermitian inner product of $W$, we assume that it is $(\cdot,\cdot)_W$ in \ref{existence of orthogonal basis} \eqref{existence of orthogonal basis 2}.

For each $\bfK$-vector space $V$ equipped with a Hermitian inner product, we use the following notation throughout this paper:
\begin{itemize}
\item $\clL_V := \{ v \in V \mid (v,v) \in \bfK_\infty \}$.
\item $\ol{\clL}_V := \clL_V/q\inv \clL_V$.
\item $\ev_\infty : \clL_V \rightarrow \ol{\clL}_V$ the quotient map.
\end{itemize}

\begin{prop}\label{characterization of local space in terms of almost orthonormal basis}
Let $V$ be a finite-dimensional $\bfK$-vector space equipped with a Hermitian inner product $(\cdot,\cdot)$ and an almost orthonormal basis $B$.
\begin{enumerate}
\item\label{characterization of local space in terms of almost orthonormal basis 1} We have $\clL_V = \bfK_\infty B$.
\item\label{characterization of local space in terms of almost orthonormal basis 2} We have $\clL_V = \{ v \in V \mid (u,v) \in \bfK_\infty \Forall u \in \clL_V \}$.
\item\label{characterization of local space in terms of almost orthonormal basis 3} $(\cdot,\cdot)$ induces a $\C$-valued Hermitian inner product (denoted by the same symbol) on $\ol{\clL}_V$ given by
$$
(\ev_\infty(u), \ev_\infty(v)) := \ev_\infty((u,v)).
$$
\item\label{characterization of local space in terms of almost orthonormal basis 4} $\ev_\infty(B)$ is an orthonormal basis of $\ol{\clL}_V$.
\end{enumerate}
\end{prop}

\begin{prop}\label{orthogonal decomposition and local space}
Let $V$ be a finite-dimensional $\bfK$-vector space equipped with a Hermitian inner product. Suppose that $V$ admits an orthogonal decomposition $V = \bigoplus_{\omega \in \Omega} V_\omega$. Then, we have orthogonal decompositions
$$
\clL_V = \bigoplus_{\omega \in \Omega} \clL_{V_\omega}, \qu \ol{\clL}_V = \bigoplus_{\omega \in \Omega} \ol{\clL}_{V_\omega},
$$
here, we identify $\clL_{V_\omega}/q\inv \clL_V$ with $\clL_{V_\omega}/ (q\inv \clL_V \cap \clL_{V_\omega})$, which equals $\ol{\clL}_{V_\omega}$.
\end{prop}

\begin{defi}\normalfont
Let $V,W$ be $\bfK$-vector spaces equipped with Hermitian inner products $(\cdot,\cdot)_V$, $(\cdot,\cdot)_W$, respectively.
\begin{enumerate}
\item We say that a linear map $f \in \Hom_{\bfK}(V,W)$ almost preserves the metrics if for each $u,v \in \clL_V$, we have $(f(u),f(v))_W - (u,v)_V \in q\inv \bfK_\infty$.
\item We say that $V$ and $W$ are almost isometric if there exists an isomorphism $V \rightarrow W$ of vector spaces which almost preserves the metrics. Such an isomorphism is called an almost isometry.
\end{enumerate}
\end{defi}

The following is an easy analog of an elementary result about complex metric spaces.

\begin{prop}\label{characterization of almost isometry}
Let $V,W$ be $\bfK$-vector spaces equipped with Hermitian inner products $(\cdot,\cdot)_V$, $(\cdot,\cdot)_W$, respectively. Let $f \in \Hom_{\bfK}(V,W)$. Then, the following are equivalent.
\begin{enumerate}
\item\label{characterization of almost isometry 1} $f$ almost preserves the metrics.
\item\label{characterization of almost isometry 2} For each almost orthonormal basis $B$ of $V$, its image $f(B)$ forms an almost orthonormal basis of $f(V)$.
\item\label{characterization of almost isometry 3} There exists an almost orthonormal basis $B$ of $V$ such that $f(B)$ forms an almost orthonormal basis of $f(V)$.
\end{enumerate}
\end{prop}

%
%

Let $V$ be a $\bfK$-vector space equipped with a Hermitian inner product, and $W \subseteq V$ a subspace. Throughout this paper, we always equip the quotient space $V/W$ with a Hermitian inner product $(\cdot,\cdot)_{V/W}$ as follows. Let $B_W \subseteq B_V$ be almost orthonormal bases of $W \subseteq V$. Set $W' := \bfK(B_V \setminus B_W)$. For each $v \in V = W \oplus W'$, let $v_1 \in W$ and $v_2 \in W'$ be such that $v = v_1+v_2$. Also, let $[v]$ denote the image of $v$ in $V/W$. Then, define 
$$
([u],[v])_{V/W} := (u_2,v_2).
$$
This is a Hermitian inner product on $V/W$. Clearly, $[B_V] := \{ [b] \mid b \in B_V \setminus B_W \}$ forms an almost orthonormal basis of $V/W$. Then, by Proposition \ref{characterization of local space in terms of almost orthonormal basis} \eqref{characterization of local space in terms of almost orthonormal basis 1}, we have
$$
\clL_{V/W} = \bfK_\infty [B_V] = \bfK_\infty B_V/(W \cap \bfK_\infty B_V) = \clL_{V}/\clL_W.
$$
This implies that although the Hermitian inner product thus constructed depends on $B_W$ and $B_V$, the $\bfK_\infty$-submodule $\clL_{V/W}$ does not.

\subsection{Quantum groups}
For $n \in \Z$ and $a \in \Z_{> 0}$, set
$$
[n]_{q^a} := \frac{q^{an} - q^{-an}}{q^a-q^{-a}}.
$$
Also, for $m \geq n \geq 0$, set
$$
[n]_{q^a}! := [n]_{q^a} \cdots [2]_{q^a}[1]_{q^a}, \qu {m \brack n}_{q^a} := \frac{[m]_{q^a}!}{[m-n]_{q^a}! [n]_{q^a}!},
$$
where we understand that $[0]_{q^a}! = 1$. When $a = 1$, we often omit the subscript $q^a$.

Let $A$ be an associative algebra over $\bfK$. For $x,y \in A$, $z \in A^\times$, $a \in \Z_{> 0}$ and $b \in \Z$, set
$$
[z;n]_{q^a} := \frac{q^{an}z-q^{-an}z\inv}{q^a - q^{-a}}, \qu \{ z;n \}_{q^a} := q^{an}z + q^{-an}z\inv, \qu [x,y]_{q^b} := xy - q^byx.
$$
When $a = 1$ or $b = 0$, we often omit the subscripts $q^a$ or $q^b$, respectively.

Let $ (a_{i,j})_{i,j \in I}$ be a Cartan matrix of finite type. We identify the index set $I$ with the set of vertices of the Dynkin diagram of the Cartan matrix, or with the Dynkin diagram itself. Let $\{ \alpha_i \}_{i \in I}$ denote the set of simple roots, $\{ h_i \}_{i \in I}$ the set of simple coroots, $\{ \vpi_i \}_{i \in I}$ the set of fundamental weights, $X := \bigoplus_{i \in I} \Z \vpi_i$ the weight lattice, $Y := \bigoplus_{i \in I} \Z h_i$ the coroot lattice, $\la \cdot, \cdot \ra : Y \times X \rightarrow \Z$ the perfect pairing given by $\la h_i, \vpi_j \ra = \delta_{i,j}$. Let $d_i \in \Z_{> 0}$ be such that $d_i a_{i,j} = d_j a_{j,i}$ for all $i,j \in I$.

Let $\g = \g(I)$ denote the finite-dimensional semisimple Lie algebra over $\C$ associated with the Dynkin diagram $I$. Let $e_i,f_i,h_i$, $i \in I$ denote the Chevalley generators.

The quantum group $\U = U_q(\g)$ is a unital associative algebra over $\bfK$ generated by $E_i,F_i,K_h$, $i \in I$, $h \in Y$ subject to the following relations; for each $i,j \in I$ and $h,h' \in Y$,
\begin{align}
\begin{split}
&K_0 = 1,\qu K_h K_{h'} = K_{h+h'}, \\
&K_h E_i = q^{\la h,\alpha_i \ra}E_i K_h, \qu K_h F_i = q^{-\la h,\alpha_i \ra} F_i K_h, \\
&[E_i,F_j] = \delta_{i,j} [K_i;0]_{q_i}, \\
&\sum_{k=0}^{1-a_{i,j}} (-1)^k {1-a_{i,j} \brack k}_{q_i} E_i^{1-a_{i,j}-k} E_j E_i^{k} = 0 \qu \IF i \neq j, \\
&\sum_{k=0}^{1-a_{i,j}} (-1)^k {1-a_{i,j} \brack k}_{q_i} F_i^{1-a_{i,j}-k} F_j F_i^{k} = 0 \qu \IF i \neq j,
\end{split} \nonumber
\end{align}
where $K_i := K_{d_ih_i}$ and $q_i := q^{d_i}$.

The quantum group $\U$ admits a Hopf algebra structure with comultiplication $\Delta$ defined by
$$
\Delta(E_i) := E_i \otimes 1 + K_i \otimes E_i, \qu \Delta(F_i) := 1 \otimes F_i + F_i \otimes K_i\inv, \qu \Delta(K_h) = K_h \otimes K_h.
$$

Let $\wp$ denote the algebra anti-involution on $\U$ defined by
$$
\wp(E_i) := q_i\inv F_iK_i, \qu \wp(F_i) := q_i\inv E_iK_i\inv, \qu \wp(K_h) := K_h.
$$
By the definitions of $\Delta$ and $\wp$, one can verify that
$$
\Delta \circ \wp = (\wp \otimes \wp) \circ \Delta.
$$

Let $(\cdot)^* : \U \rightarrow \U$ denote the ring automorphism which extends the complex conjugate on $\bfK$ in a way such that
$$
E_i^* = E_i,\ F_i^* = F_i,\ K_h^* = K_h, \Forall i \in I,\ h \in Y.
$$
Then, we have $\wp \circ * = * \circ \wp$, and $\Delta \circ * = (* \otimes *) \circ \Delta$. Therefore, $\wp^* := \wp \circ *$ is an algebra anti-involution, and satisfies
\begin{align}\label{Delta and wpstar}
\Delta \circ \wp^* = (\wp^* \otimes \wp^*) \circ \Delta.
\end{align}

Let $\Br = \Br(I)$ denote the braid group associated with the Dynkin diagram $I$ with generators $T_i$, $i \in I$. The braid group $\Br$ acts on $\U$ via Lusztig's automorphisms $T''_{i,1}$, $i \in I$ (\cite[37.1.3]{L10}):
\begin{align}
\begin{split}
&T_i(E_j) := \begin{cases}
-F_iK_i \qu & \IF j = i, \\
\sum_{r+s=-a_{i,j}} (-1)^r q_i^{-r} E_i^{(s)} E_j E_i^{(r)} \qu & \IF j \neq i,
\end{cases} \\
&T_i(F_j) := \begin{cases}
-K_i\inv E_i \qu & \IF j = i, \\
\sum_{r+s=-a_{i,j}} (-1)^r q_i^{r} F_i^{(r)} F_j F_i^{(s)} \qu & \IF j \neq i,
\end{cases} \\
&T_i K_h := K_{h - \la h,\alpha_i \ra h_i},
\end{split} \nonumber
\end{align}
where $E_i^{(n)} := \frac{1}{[n]_{q_i}!} E_i^n$ and $F_i^{(n)} := \frac{1}{[n]_{q_i}!} F_i^n$ are divided powers.

\subsection{Weight $\U$-modules}
Given a $\U$-module $M$ and a weight $\lm \in X$, the subspace
$$
M_\lm := \{ m \in M \mid K_h m = q^{\la h,\lm \ra}m \ \Forall h \in Y \},
$$
is called the weight space of $M$ of weight $\lm$, and each element of $M_\lm$ is said to be a weight vector of weight $\lm$.
A $\U$-module is said to be a weight module if it admits a weight space decomposition $M = \bigoplus_{\lambda \in X} M_\lambda$.

For each $\lm,\mu \in X$, we write $\mu \leq \lm$ to indicate that $\lm - \mu \in \sum_{i \in I} \Z_{\geq 0} \alpha_i$. This defines a partial order on $X$, called the dominance order.

Let $X^+ := \{ \lm \in X \mid \la h_i,\lm \ra \geq 0 \ \Forall i \in I \}$ denote the set of dominant weights.

\begin{theo}[{see e.g. \cite[Part I]{L10}}]\label{facts about weight U-modules}
\ \begin{enumerate}
\item For each weight $\U$-module $M$, $\lm \in X$, and $i \in I$, we have
$$
E_i M_\lm \subseteq M_{\lm + \alpha_i}, \qu F_i M_\lm \subseteq M_{\lm-\alpha_i}.
$$
\item Each finite-dimensional $\U$-module is completely reducible.
\item The isoclasses of finite-dimensional irreducible $\U$-modules are parametrized by $X^+$. Let $V(\lambda)$ denote the finite-dimensional irreducible $\mathbf{U}$-module of highest weight $\lambda \in X^+$, and $v_\lambda \in V(\lambda)$ the highest weight vector.
\item\label{spanning property} For each $\lm \in X^+$, we have
$$
V(\lm) = \Span_{\bfK} \{ F_{i_1} \cdots F_{i_r} v_\lm \mid r \geq 0,\ i_1,\ldots,i_r \in I \}.
$$
\item\label{weight constraint} For each $\lm \in X^+$ and $\mu \in X$, we have $\dim V(\lm)_{\lm} = 1$, and $\dim V(\lm)_\mu = 0$ unless $\mu \leq \lm$.
\item For each $\lm \in X^+$ and $\mu \in X$ such that $\mu < \lm$, we have
$$
V(\lm)_\mu = \sum_{i \in I} F_i V(\lm)_{\mu + \alpha_i}.
$$
\end{enumerate}
\end{theo}

\subsection{$\U$-modules with contragredient Hermitian inner products}
\begin{defi}\label{Definition of contragredient Hermitian inner products}\normalfont
Let $M$ be a $\U$-module equipped with a Hermitian inner product $(\cdot,\cdot)$. We say that $(\cdot,\cdot)$ is contragredient if it satisfies
\begin{align}\label{Property1}
(xu,v) = (u,\wp^*(x)v) \qu \Forall x \in \U,\ u,v \in M.
\end{align}
\end{defi}

\begin{prop}\label{contragredient Hermitian inner product on V(lm)}
Let $\lm \in X^+$, and consider the irreducible $\U$-module $V(\lm)$. Then, there exists a unique contragredient Hermitian inner product $(\cdot,\cdot)_\lm$ on $V(\lm)$ such that $(v_\lm,v_\lm)_\lm = 1$.
\end{prop}

\begin{proof}
Let $\U_{\Q}$ denote the $\Q(q)$-subalgebra of $\U$ generated by $E_i,F_i,K_h$, $i \in I$, $h \in Y$, and set $V(\lm)_{\Q} := \U_{\Q} v_\lm$. Then, we have
$$
\U = \U_{\Q} \otimes_{\Q(q)} \bfK, \qu V(\lm) = V(\lm)_{\Q} \otimes_{\Q(q)} \bfK.
$$

By \cite[Proposition 19.1.2]{L10}, there exists a unique symmetric bilinear form $(\cdot,\cdot)'_\lm : V(\lm)_{\Q} \times V(\lm)_{\Q} \rightarrow \Q(q)$ such that $(v_\lm,v_\lm)'_\lm = 1$, and $(xu,v)'_\lm = (u,\wp(x)v)'_\lm$ for all $x \in \U_{\Q}$ and $u,v \in V(\lm)_{\Q}$. Define a map $(\cdot,\cdot)_\lm : V(\lm) \times V(\lm) \rightarrow \bfK$ by
$$
(au,bv)_\lm := ab^*(u,v)'_\lm, \qu a,b \in \bfK,\ u,v \in V(\lm)_{\Q}.
$$
Then, $(\cdot,\cdot)_\lm$ satisfies conditions \eqref{Definition of Hermitian inner product 1} and \eqref{Definition of Hermitian inner product 2} in Definition \ref{Definition of Hermitian inner product}, condition \eqref{Property1} in Definition \ref{Definition of contragredient Hermitian inner products}, and $(v_\lm,v_\lm)_\lm = 1$.

Moreover, by \cite[Lemma 19.1.4]{L10}, there exists a basis $B$ of $V(\lm)_{\Q}$ such that
$$
(b,b')'_\lm \in \delta_{b,b'} + q\inv \bfK_\infty.
$$
This shows that $(\cdot,\cdot)_\lm$ satisfies conditions \eqref{Definition of Hermitian inner product 3}--\eqref{Definition of Hermitian inner product 5} in Definition \ref{Definition of Hermitian inner product}. Thus the proof completes.
\end{proof}

\begin{lem}\label{Tensor product of Hermitian form}
Let $M$ and $N$ be $\U$-modules with contragredient Hermitian inner products $(\cdot,\cdot)_M$ and $(\cdot,\cdot)_N$. Then, the form $(\cdot,\cdot)_{M,N}$ on $M \otimes N$ given by
$$
(m_1 \otimes n_1, m_2 \otimes n_2)_{M,N} := (m_1,m_2)_M (n_1,n_2)_N
$$
is a contragredient Hermitian inner product.
\end{lem}

\begin{proof}
The assertion follows from equation \eqref{Delta and wpstar} on page \pageref{Delta and wpstar}.
\end{proof}

\begin{lem}\label{orthogonality of different weight spaces}
Let $M$ be a weight $\U$-module equipped with a contragredient Hermitian inner product. Then, we have $(M_\lm,M_\mu) = 0$ for all $\lm,\mu \in X$ such that $\lm \neq \mu$.
\end{lem}

\begin{proof}
The assertion follows from the fact that $\wp^*(K_h) = K_h$ for all $h \in Y$.
\end{proof}

\begin{prop}\label{orthogonal irreducible decomposition}
Let $M$ be a finite-dimensional $\U$-module equipped with a contragredient Hermitian inner product. Then, there exists an irreducible decomposition $M = \bigoplus_{k=1}^r M_k$ satisfying the following:
\begin{enumerate}
\item $(M_k,M_l) = 0$ unless $k = l$.
\item For each $k \in [1,r]$, there exists a unique $\lm_k \in X^+$ and an isomorphism $\phi_k : M_k \rightarrow V(\lm_k)$ of $\U$-modules which is an almost isometry.
\end{enumerate}
\end{prop}

\begin{proof}
For each $\lm \in X^+$, set $H_\lm := \{ m \in M_\lm \mid E_i m = 0 \Forall i \in I \}$. Then, by Proposition \ref{existence of orthogonal basis} \eqref{existence of orthogonal basis 1}, $H_\lm$ has an orthogonal basis. By Lemma \ref{orthogonality of different weight spaces}, we have $(H_\lm,H_\mu) = 0$ if $\lm \neq \mu$. Then one can take an orthogonal basis $m_1,\ldots,m_r$ of $\bigoplus_{\lm \in X^+} H_\lm$ consisting of weight vectors. Setting $M_k := \U m_k$, we obtain an irreducible decomposition $M = \bigoplus_{k=1}^r M_k$ of $M$.

For each $k \in [1,r]$, let $\lm_k \in X^+$ be such that $m_k \in H_{\lm_k}$. Then, we have $M_k \simeq V(\lm_k)$. Each element of $M_k$ is of the form $xm_k$ for some $x \in \U$. Let $k \neq l$ and $x,y \in \U$. Without loss of generality, we may assume that $\lm_k \geq \lm_l$. Then, we have
$$
(xm_k,ym_l) = (m_k, \wp^*(x)y m_l).
$$
Since $\wp^*(x)y m_l \in M_l \simeq V(\lm_l)$, by Theorem \ref{facts about weight U-modules} \eqref{weight constraint}, it is a linear combination of weight vectors of weight less than or equal to $\lm_l$, while $m_k$ is of weight $\lm_k (\geq \lm_l)$. Hence, if $\lm_l \neq \lm_k$, then we have $(m_k, \wp^*(x)y m_l) = 0$ by Lemma \ref{orthogonality of different weight spaces}. Otherwise, one can write as
$$
\wp^*(x)y m_l = c m_l + m',
$$
where $c \in \bfK$, and $m'$ is a sum of weight vectors of weights less than $\lm_l (= \lm_k)$. Then, we have
$$
(m_k, \wp^*(x)y m_l) = c(m_k,m_l) = 0
$$
since $m_k$ and $m_l$ are orthogonal. Thus, we obtain that $(M_k,M_l) = 0$, as desired.

For a proof of the second assertion, replace $m_k$ with $\frac{1}{\lt(m_k)} m_k$ in the argument above. Then, we have $\lt(m_k) = 1$. Since $m_k$ is a highest weight vector of weight $\lm_k$, there exists an isomorphism $\phi_k : M_k \rightarrow V(\lm_k)$ of $\U$-modules such that $\phi_k(m_k) = v_{\lm_k}$. Then, for each $u,v \in M_k$, we have
$$
(u,v) = (m_k,m_k) (\phi_k(u),\phi_k(v))_{\lm_k}.
$$
In particular, for each $u,v \in \clL_{M_k}$, we have
$$
(\phi_k(u),\phi_k(v))_{\lm_k} - (u,v) = ((m_k,m_k)\inv - 1)(u,v) \in q\inv \bfK_\infty.
$$
This implies that $\phi_k$ is an almost isometry. Thus, the proof completes.
\end{proof}

\subsection{Crystal bases of $\U$-modules}
\begin{defi}\normalfont
Let $V$ be a $\bfK$-vector space and $A$ a subring of $\bfK$. An $A$-lattice of $V$, or a lattice of $V$ over $A$, is a free $A$-submodule $L$ of $V$ such that $L \otimes_A \bfK = V$.
\end{defi}

\begin{ex}\normalfont
Let $V$ be a finite-dimensional $\bfK$-vector space equipped with a Hermitian inner product. Then, $\clL_V$ is a $\bfK_\infty$-lattice of $V$.
\end{ex}

Let $M$ be a finite-dimensional $\U$-module equipped with a contragredient Hermitian inner product. Let $\Etil_i,\Ftil_i$, $i \in I$ denote the Kashiwara operators. They preserve $\clL_M$, and hence, induce $\C$-linear operators on $\ol{\clL}_M$. A $\C$-basis of $\ol{\clL}_M$ satisfying certain conditions is called a crystal basis of $M$ (see e.g. \cite{Ka90} for detail).
For each $b \in \clB_M$ and $i \in I$, set
$$
\vphi_i(b) := \max\{ k \mid \Ftil_i^k b \neq 0 \}, \qu \vep_i(b) := \max\{ k \mid \Etil_i^k b \neq 0 \}.
$$

Let $M$ and $N$ be $\U$-modules with crystal bases $\clB_M$ and $\clB_N$. Then, $\clB_M \otimes \clB_N$ becomes a crystal basis of $M \otimes N$ on which the Kashiwara operators and $\vphi_i,\vep_i$ act as
\begin{align}\label{tensor product rule for usual crystal}
\begin{split}
&\Ftil_i(b_1 \otimes b_2) = \begin{cases}
b_1 \otimes \Ftil_i(b_2) \qu & \IF \vep_i(b_1) < \vphi_i(b_2), \\
\Ftil_i(b_1) \otimes b_2 \qu & \IF \vep_i(b_1) \geq \vphi_i(b_2),
\end{cases} \\
&\Etil_i(b_1 \otimes b_2) = \begin{cases}
b_1 \otimes \Etil_i(b_2) \qu & \IF \vep_i(b_1) \leq \vphi_i(b_2), \\
\Etil_i(b_1) \otimes b_2 \qu & \IF \vep_i(b_1) > \vphi_i(b_2),
\end{cases} \\
&\vphi_i(b_1 \otimes b_2) = \vphi_i(b_1) + \max(0, \vphi_i(b_2)-\vep_i(b_1)), \\
&\vep_i(b_1 \otimes b_2) = \vep_i(b_2) + \max(0, \vep_i(b_1)-\vphi_i(b_2)).
\end{split}
\end{align}
These formulas can be found, for example, in \cite[Theorem 4.1]{HK02}; one should note that $\Ftil_i,\Etil_i,\vphi_i,\vep_i$ in this paper are $\Etil_i,\Ftil_i,\vep_i,\vphi_i$ there.

When we consider the irreducible $\U$-module $V(\lm)$, $\lm \in X^+$ equipped with the contragredient Hermitian inner product $(\cdot,\cdot)_\lm$ constructed in Proposition \ref{contragredient Hermitian inner product on V(lm)}, we set
$$
\clL(\lm) := \clL_{V(\lm)}, \qu \ol{\clL}(\lm) := \ol{\clL}_{V(\lm)}, \qu b_\lm := \ev_\infty(v_\lm).
$$
By \cite[Theorem 2]{Ka91}, $V(\lm)$ possesses a unique crystal basis $\clB(\lm) \subseteq \ol{\clL}(\lm)$ containing $b_\lm$.

Now, let us turn back to a general $\U$-module $M$. By Proposition \ref{orthogonal irreducible decomposition}, we have an orthogonal irreducible decomposition $M = \bigoplus_{k=1}^r M_k$. Then, by Proposition \ref{orthogonal decomposition and local space}, we have
$$
\clL_M = \bigoplus_{k=1}^r \clL_k, \qu \ol{\clL}_M = \bigoplus_{k=1}^r \ol{\clL}_k,
$$
where $\clL_k := \clL_M \cap M_k$ and $\ol{\clL}_k := \clL_k/q\inv \clL_k$. Moreover, the almost isometries $\phi_k$ in Proposition \ref{orthogonal irreducible decomposition} induce isomorphisms
$$
\clL_k \simeq \clL(\lm_k), \qu \ol{\clL}_k \simeq \ol{\clL}(\lm_k).
$$
Let $\clB_k$ denote the basis of $\ol{\clL}_k$ corresponding to $\clB(\lm_k)$ under the isomorphism $\ol{\clL}_k \simeq \ol{\clL}(\lm_k)$ above. Then, we obtain a crystal basis
$$
\clB_M := \bigsqcup_{k=1}^r \clB_k
$$
of $M$. Furthermore, every crystal basis of $M$ is of this form.

For each $\lm \in X^+$, set $M[\lm] := \bigoplus_{\substack{k \in [1,r] \\ \lm_k = \lm}} M_k$ to be the sum of all submodules of $M$ isomorphic to $V(\lm)$. Set
\begin{align}
\begin{split}
&\clL_M[\lm] := \clL_M \cap M[\lm] = \bigoplus_{\substack{k \in [1,r] \\ \lm_k = \lm}} \clL_k, \\
&\ol{\clL}_M[\lm] := \clL_M[\lm]/q\inv \clL_M[\lm] = \bigoplus_{\substack{k \in [1,r] \\ \lm_k = \lm}} \ol{\clL}_k, \\
&\clB_M[\lm] := \clB_M \cap \ol{\clL}_M[\lm] = \bigsqcup_{\substack{k \in [1,r] \\ \lm_k = \lm}} \clB_k.
\end{split} \nonumber
\end{align}
Also, set
$$
M[\geq \lm] := \bigoplus_{\mu \geq \lm} M[\mu], \qu M[> \lm] := \bigoplus_{\mu > \lm} M[\mu],
$$
and define $\clL_M[\geq \lm], \clL_M[> \lm],\ol{\clL}_M[\geq \lm], \ol{\clL}_M[> \lm],\clB_M[\geq \lm], \clB_M[> \lm]$ in the obvious way. Then, $M[\geq \lm]$, $M[> \lm]$, and $M[\geq \lm]/M[> \lm]$ are $\U$-modules possessing crystal bases $\clB_M[\geq \lm]$, $\clB_M[> \lm]$, and $\{ [b] \mid b \in \clB_M[\lm] \}$, respectively, where $[b]$ denotes the image of $b \in \ol{\clL}[\geq \lm]$ in $\ol{\clL}[\geq \lm]/\ol{\clL}[> \lm]$.

\subsection{Based $\U$-modules}
Let $\Udot = \bigoplus_{\lm,\mu \in X} {}_{\lm}\U_\mu = \bigoplus_{\lm \in X} \Udot \mathbf{1}_\lm$ denote the modified quantum group, where the $\mathbf{1}_\lambda$ are the associated orthogonal idempotents (see \cite[Chapter 23]{L10} for a precise definition).

\begin{defi}[{\cite[23.1.4]{L10}}]\normalfont
A unital $\Udot$-module is a $\Udot$-module $M$ such that for each $m \in M$, we have the following:
\begin{enumerate}
\item $\mathbf{1}_\lm m = 0$ for all but finitely many $\lm \in X$.
\item $\sum_{\lm \in X} \mathbf{1}_\lm m = m$.
\end{enumerate}
\end{defi}
As explained in \cite[23.1.4]{L10}, a weight $\mathbf{U}$-module has a unital $\Udot$-module structure, and vice versa.

Set
$$
\bfA := \C[q,q\inv] \subseteq \bfK.
$$
Let $\UdotA$ denote the $\bfA$-form of $\Udot$. It is an $\bfA$-subalgebra of $\Udot$ generated by the divided powers $E_i^{(n)} \mathbf{1}_\lm$, $F_i^{(n)} \mathbf{1}_\lm$, $i \in I$, $\lm \in X$, $n \in \Z_{\geq 0}$.

\begin{defi}\normalfont
Let $M$ be a weight $\U$-module. An $\bfA$-form of $M$ is an $\bfA$-lattice $M_{\bfA}$ of $M$ which is simultaneously a $\UdotA$-submodule.
\end{defi}

\begin{defi}\normalfont
Let $M$ be a $\U$-module. A bar-involution $\psi_M$ on $M$ is an involutive $\C$-linear automorphism such that
$$
\psi_M(xm) = \psi(x) \psi_M(m) \qu \Forall x \in \U,\ m \in M.
$$
\end{defi}

\begin{defi}\normalfont
Let $M$ be a weight $\U$-module equipped with a contragredient Hermitian inner product $(\cdot,\cdot)_M$, an $\bfA$-form $M_{\bfA}$, a bar-involution $\psi_M$, and a crystal basis $\clB_M$. We say that $(M,(\cdot,\cdot)_M, M_{\bfA}, \psi_M, \clB_M)$, or simply $M$, is a based $\U$-module if the following two conditions are satisfied:
\begin{enumerate}
\item The quotient map $\ev_\infty : \clL_M \rightarrow \ol{\clL}_M$ restricts to an isomorphism
$$
\clL_M \cap M_{\bfA} \cap \psi_M(\clL_M) \rightarrow \ol{\clL}_M
$$
of $\C$-vector spaces; let $G : \ol{\clL}_M \rightarrow \clL_M \cap M_{\bfA} \cap \psi_M(\clL_M)$ denote its inverse.
\item For each $b \in \clB_M$, we have $\psi_M(G(b)) = G(b)$.
\end{enumerate}
\end{defi}

\begin{rem}\normalfont
  This definition is equivalent to \cite[27.1.2]{L10}.
  We choose this definition so that we can straightforwardly generalize it to the $\imath$quantum groups setting.
\end{rem}

\begin{ex}\normalfont
Let $\lm \in X^+$. Then, the irreducible $\U$-module $V(\lm)$ possesses a unique bar-involution $\psi_\lm$ fixing the highest weight vector $v_\lm$. Set $V(\lm)_{\bfA} := \UdotA v_\lm$. Then, $(V(\lm), (\cdot,\cdot)_\lm, V(\lm)_{\bfA}, \psi_\lm, \clB(\lm))$ is a based $\U$-module. The set $G(\lm) := G(\clB(\lm))$ is called the canonical basis of $V(\lm)$.
\end{ex}

A finite-dimensional based $\U$-module is cellular in the following sense.

\begin{prop}[{\cite[Propositions 27.1.7 and 27.1.8]{L10}}]\label{cellularity of canonical basis of based U-module}
Let $(M, (\cdot,\cdot)_M, M_{\bfA}, \psi_M)$ be a finite-dimensional based $\U$-module, and $\clB_M$ a crystal basis. Then, for each $\lm \in X^+$, the $\U$-modules $M[\geq \lm]$, $M[> \lm]$, and $M[\geq \lm]/M[> \lm]$ are based $\U$-modules with the same (or induced) contragredient Hermitian inner products and bar-involutions as $M$, and the $\bfA$-forms spanned by $G(\clB_M[\geq \lm])$, $G(\clB_M[> \lm])$, $[G(\clB_M[\lm])] := \{ G(b) + M[> \lm] \mid b \in \clB_M[\lm] \}$, respectively. Furthermore, there exists an isomorphism
$$
M[\geq \lm]/M[> \lm] \rightarrow V(\lm)^{\oplus m_\lm}
$$
of $\U$-modules which restricts to a bijection
$$
[G(\clB_M[\lm])] \rightarrow G(\lm)^{m_\lm}
$$
between the canonical bases, where $m_\lm := \dim_{\bfK} \Hom_{\U}(V(\lm),M)$ denotes the multiplicity of $V(\lm)$ in $M$.
\end{prop}

\section{Quantum symmetric pairs}\label{section: QSP}
In this section, we formulate the $\imath$quantum groups following Kolb \cite{Ko14}. Then, we introduce the notion of based modules over the $\imath$quantum groups generalizing that over the quantum groups.

\subsection{$\imath$Quantum groups}
Let $(I,I_\bullet,\tau)$ be a Satake diagram, and $\frk$ the associated symmetric pair subalgebra (see e.g., \cite[Table 4]{BW18b} for the list of Satake diagrams). Namely, $I_\bullet$ is a subdiagram of $I$, and $\tau$ is a diagram involution on $I$ satisfying certain conditions (see \cite[Definition 2.3]{Ko14}), and $\frk$ is the Lie subalgebra of $\g$ generated by the following elements:
\begin{itemize}
\item $e_i$, $i \in I_\bullet$.
\item $h^\tau_i := \begin{cases}
h_i \qu & \IF i \in I_\bullet, \\
h_i - h_{\tau(i)} \qu & \IF i \in I_\circ.
\end{cases}$
\item $b_i := \begin{cases}
f_i \qu & \IF i \in I_\bullet, \\
f_i + \ol{\varsigma}_i \ol{T}_{w_\bullet}(e_{\tau(i)}) \qu & \IF i \in I_\circ.
\end{cases}$
\end{itemize}
Here, $I_\circ := I \setminus I_\bullet$, $w_\bullet$ is the longest element of the Weyl group of $I_\bullet$, $\ol{T}_{w_\bullet}$ is the braid group action on $\g$ in terms of triple exponentials, and $\ol{\varsigma}_i \in \C$ are certain scalars. There exists a Lie algebra involution $\theta$ on $\g$ such that
$$
\frk = \g^\theta := \{ x \in \g \mid \theta(x) = x \}.
$$
It is known that $\frk$ is a reductive Lie algebra (\cite[Proposition 1.13.3]{D96}).

Let $\Ui = \Ui_{\bfvarsigma,\bfkappa}$ denote the $\imath$quantum group associated with the Satake diagram $(I,I_\bullet,\tau)$ and parameters $\bfvarsigma = (\varsigma_i)_{i \in I_\circ} \in (\bfK_1^\times)^{I_\circ}$, $\bfkappa = (\kappa_i)_{i \in I_\circ} \in \bfK_1^{I_\circ}$ obeying constraints \cite[Equations (3.4)--(3.8) and the assumption in Proposition 4.6]{BW18b}. Namely, $\Ui$ is the subalgebra of $\U$ generated by $E_i$, $i \in I_\bullet$, $B_i,k_i$, $i \in I$ where
\begin{align}
\begin{split}
&k_i^{\pm 1} := \begin{cases}
K_i^{\pm 1} \qu & \IF i \in I_\bullet, \\
(K_i K_{\tau(i)}\inv)^{\pm 1} \qu & \IF i \in I_\circ,
\end{cases} \\
&B_i := \begin{cases}
F_i \qu & \IF i \in I_\bullet, \\
F_i + \varsigma_i T_{w_\bullet}(E_{\tau(i)}) K_i\inv + \kappa_i K_i\inv \qu & \IF i \in I_\circ.
\end{cases}
\end{split}\nonumber
\end{align}

One of the most distinguished properties of $\Ui$ is that it is a right coideal of $\U$, i.e., $\Delta(\Ui) \subseteq \Ui \otimes \U$. Therefore, the tensor product of a $\Ui$-module and a $\U$-module becomes a $\Ui$-module via $\Delta$.

\begin{prop}[{\cite[Proposition 4.6]{BW18b}}]\label{wp on Ui}
For each $i \in I_\circ$, we have $\wp(B_i) \in \Ui$, and
$$
\wp(B_i) = q_i\inv \varsigma_{\tau(i)}\inv T_{w_\bullet}\inv(B_{\tau(i)}) K_{w_\bullet(h_{\tau(i)})}K_i\inv.
$$
Consequently, the involution $\wp^*$ on $\U$ preserves $\Ui$.
\end{prop}

\begin{prop}[{\cite[Theorem 3.11]{BK15}}]
There exists a unique $\C$-algebra involution $\psii$, called the $\imath$bar-involution on $\Ui$ such that
$$
\psii(E_i) = E_i, \qu \psii(B_i) = B_i, \qu \psii(k_i) = k_i\inv, \qu \psii(q) = q\inv.
$$
\end{prop}

\subsection{$\imath$Canonical bases of based $\U$-modules}
\begin{defi}\normalfont
Let $M$ be a $\Ui$-module. An $\imath$bar-involution on $M$ is a $\C$-linear involution $\psii_M$ on $M$ such that
$$
\psii_M(xm) = \psii(x) \psii_M(m) \qu \Forall x \in \Ui \AND m \in M.
$$
\end{defi}

Let $M$ be a finite-dimensional $\U$-module equipped with a bar-involution $\psi_M$. Let $\Upsilon$ denote the quasi-$K$-matrix (see \cite{BW18b} for a precise definition). Then, $\psii_M := \Upsilon \circ \psi_M$ becomes an $\imath$bar-involution of the $\Ui$-module $M$. We call it the $\imath$bar-involution associated with $\psi_M$.

\begin{theo}[{\cite[Theorem 5.7]{BW18b}}]\label{icanonical basis of based Umodules}
Let $M$ be a finite-dimensional based $\U$-module. Then, for each $b \in \clB_M$, there exists a unique $G^\imath(b) \in \clL_M \cap M_{\bfA}$ satisfying the following:
\begin{enumerate}
\item $\psii_M(G^\imath(b)) = G^\imath(b)$.
\item $\ev_\infty(G^\imath(b)) = b$.
\end{enumerate}
\end{theo}

We call $G^\imath(\clB_M)$ the $\imath$canonical basis associated with the canonical basis $G(\clB_M)$. When $\clB_M = \clB(\lm)$ for some $\lm \in X^+$, we denote by $G^\imath(\lm)$ the $\imath$canonical basis associated with $G(\lm)$.

\begin{lem}\label{weight spaces iCB elements live}
Let $M$ be a finite-dimensional based $\U$-module. Let $M = \bigoplus_{k=1}^r M_k$ be the irreducible decomposition such that $\clB_M = \bigsqcup_{k=1}^r \clB_k$, where $\clB_k := \clB_M \cap \ol{\clL}_{M_k}$. Let $\lm_k \in X^+$ be such that $M_k \simeq V(\lm_k)$. Then, for each $b \in \clB_k$, we have
$$
G^\imath(b) \in M_k \oplus \bigoplus_{\substack{l \in [1,r] \\ \lm_l > \lm_k}} M_l.
$$
\end{lem}

\begin{proof}
From Proposition \ref{cellularity of canonical basis of based U-module}, we see that
$$
G(b') \in M_{l} \oplus \bigoplus_{\substack{l' \in [1,r] \\ \lm_{l'} > \lm_l}} M_{l'}.
$$
for all $l \in [1,r]$, $b' \in \clB_l$. Then, by the construction of $G^\imath(b)$ in \cite[Theorem 5.7]{BW18b}, we have
$$
G^\imath(b) \in M_k + \bigoplus_{\substack{l \in [1,r] \\ \lm_l > \lm_k}} \bigoplus_{b' \in \clB_l} \bfK G(b').
$$
Combining these two facts, we obtain the assertion.
\end{proof}

\subsection{Modified $\imath$quantum groups}
Following \cite[3.1]{BW18b}, we set
\begin{align}\label{Xi and Yi}
X^\imath := X/\{ \lm + w_\bullet \circ \tau(\lm) \mid \lm \in X \}, \qu Y^\imath := \{ h \in Y \mid h + w_\bullet \circ \tau(h) = 0 \}.
\end{align}
Then, the perfect pairing on $Y \times X$ induces a bilinear map $Y^\imath \times X^\imath \rightarrow \Z$ given by
$$
\la h,\ol{\lm} \ra := \la h,\lm \ra, \qu h \in Y^\imath, \lm \in X,
$$
where $\ol{\lm}$ denotes the image of $\lm$ in $X^\imath$.

Recall that $Q = \sum_{i \in I} \Z \alpha_i$ is the root lattice. From the definition of quantum groups, we see that $\U$ is a $Q$-graded algebra such that $E_i,F_i,K_i^{\pm 1}$ are homogeneous of degree $\alpha_i,-\alpha_i,0$, respectively. Set
$$
Q^\imath := \{ \ol{\lm} \mid \lm \in Q \} \subseteq X^\imath.
$$
Then, $\U$ is $Q^\imath$-graded, and $E_i,B_i,k_i^{\pm 1}$ are homogeneous of degree $\ol{\alpha_i}, -\ol{\alpha_i}, \ol{0}$, respectively. Hence, $\Ui$ is a $Q^\imath$-graded algebra. Let $\Ui(\ol{\lm})$ denote the homogeneous part of degree $\ol{\lm} \in Q^\imath$. For $\zeta \in X^\imath \setminus Q^\imath$, set $\Ui(\zeta) := 0$. Then, we have
$$
\Ui = \bigoplus_{\ol{\lm} \in Q^\imath} \Ui(\ol{\lm}) = \bigoplus_{\zeta \in X^\imath} \Ui(\zeta).
$$

For each $\zeta,\eta \in X^\imath$, let $\pi_{\zeta,\eta}$ denote the composition
$$
\pi_{\zeta,\eta} : \Ui(\zeta-\eta) \hookrightarrow \Ui \twoheadrightarrow \Ui/(\sum_{h \in Y^\imath}(K_h - q^{\la h,\zeta \ra})\Ui + \sum_{h \in Y^\imath} \Ui(K_h - q^{\la h,\eta \ra})),
$$
and set ${}_\zeta \Ui_\eta$ to be the image of $\pi_{\zeta,\eta}$. Clearly, we have ${}_\zeta \Ui_\eta = 0$ if $\zeta - \eta \notin Q^\imath$.

The modified $\imath$quantum group $\Uidot$ is defined by
$$
\Uidot := \bigoplus_{\substack{\zeta,\eta \in X^\imath \\ \zeta - \eta \in Q^\imath}} {}_\zeta \Ui_\eta = \bigoplus_{\zeta,\eta \in X^\imath} {}_\zeta \Ui_\eta.
$$

\begin{rem}\normalfont
  This definition slightly differs from the one in \cite{BW18b}. However, it is what the authors of \cite{BW18b} actually considered. The author would like to thank Weiqiang Wang for a fruitful discussion on this matter.
\end{rem}

The modified $\imath$quantum group $\Uidot$ admits an associative algebra structure just like $\dot{\mathbf{U}}$ does (see \cite[Chapter 23]{L10}). For each $\zeta_1,\zeta_2,\eta_1,\eta_2 \in X^\imath$ and $x_i \in \Ui(\zeta_i-\eta_i)$, set
$$
\pi_{\zeta_1,\eta_1}(x_1) \pi_{\zeta_2,\eta_2}(x_2) := \delta_{\eta_1,\zeta_2} \pi_{\zeta_1,\eta_2}(x_1x_2).
$$

For each $\zeta \in X^\imath$, set
$$
\mathbf{1}_\zeta := \pi_{\zeta,\zeta}(1).
$$
These $\mathbf{1}_\zeta$'s form a family of orthogonal idempotents:
$$
\mathbf{1}_\zeta \mathbf{1}_\eta = \delta_{\zeta,\eta} \mathbf{1}_\zeta.
$$
Also, we have
$$
{}_\zeta \Ui_\eta = \mathbf{1}_\zeta \Uidot \mathbf{1}_\eta.
$$

The modified $\imath$quantum group $\Uidot$ is equipped with a $\Ui$-bimodule structure defined as follows: For each $x,z \in \Ui$ and $y \in \Uidot$, set
$$
xyz := \sum_{\zeta,\eta \in X^\imath} (x \mathbf{1}_\zeta) y (\mathbf{1}_\eta z).
$$
Especially, we have for each $i \in I$ and $j \in I_\bullet$,
\begin{align}\label{commutation relations for B,E,k with 1zeta}
\begin{split}
&E_j \mathbf{1}_\zeta = \mathbf{1}_{\zeta + \ol{\alpha_j}} E_j, \\
&B_i \mathbf{1}_\zeta = \mathbf{1}_{\zeta - \ol{\alpha_i}} B_i, \\
&k_i^{\pm 1} \mathbf{1}_\zeta = \mathbf{1}_\zeta k_i^{\pm 1} = q^{\pm \la h^\tau_i,\zeta \ra} \mathbf{1}_\zeta.
\end{split}
\end{align}

\begin{defi}\normalfont
A $\Ui$-module $M$ is said to be an $X^\imath$-weight module if it admits an $X^\imath$-gradation
$$
M = \bigoplus_{\zeta \in X^\imath} M_{\zeta}
$$
satisfying the following:
\begin{enumerate}
\item $E_j M_\zeta \subseteq M_{\zeta + \ol{\alpha_j}}$ for all $j \in I_\bullet$, $\zeta \in X^\imath$.
\item $B_i M_\zeta \subseteq M_{\zeta - \ol{\alpha_i}}$ for all $i \in I$, $\zeta \in X^\imath$.
\item $k_i^{\pm 1} m = q^{\pm \la h^\tau_i,\zeta \ra} m$ for all $i \in I$, $\zeta \in X^\imath$, $m \in M_\zeta$.
\end{enumerate}
The subspace $M_\zeta$ (resp., a vector in $M_\zeta$) is referred to as the $X^\imath$-weight space (resp., an $X^\imath$-weight vector) of weight $\zeta$.
\end{defi}

\begin{ex}\label{canonical Xi weight module structure}\normalfont
Let $M = \bigoplus_{\lm \in X} M_\lm$ be a weight $\U$-module. Then, it admits an $X^\imath$-weight module structure such that
$$
M_\zeta = \bigoplus_{\substack{\lm \in X \\ \ol{\lm} = \zeta}} M_\lm
$$
for all $\zeta \in X^\imath$. We call it the canonical $X^\imath$-weight module structure of $M$.
\end{ex}

\begin{rem}\normalfont
Since the pairing $Y^\imath \times X^\imath \rightarrow \Z$ is not perfect, there may exist $\zeta,\eta \in X^\imath$ such that $\la h,\zeta \ra = \la h,\eta \ra$ for all $h \in Y^\imath$, but $\zeta \neq \eta$. Hence, the condition ``$k_i^{\pm 1} m = q^{\pm \la h^\tau_i,\zeta \ra} m$ for all $i \in I$'' cannot ensure that $m \in M_\zeta$.
\end{rem}

\begin{defi}\normalfont
A unital $\Uidot$-module is a $\Uidot$-module $M$ such that for each $m \in M$, we have
\begin{enumerate}
\item $\mathbf{1}_\zeta m = 0$ for all but finitely many $\zeta \in X^\imath$.
\item $\sum_{\zeta \in X^\imath} \mathbf{1}_\zeta m = m$.
\end{enumerate}
\end{defi}

As in the quantum groups case, each unital $\Uidot$-module is equipped with an $X^\imath$-weight module structure, and vice versa.

Let $M$ be a unital $\Uidot$-module. Then, it admits a $\Ui$-module structure. This action is compatible with the $\Ui$-bimodule structure of $\Uidot$.

\begin{defi}\normalfont
Let $M$ be an $X^\imath$-weight module. A $\Ui$-submodule $N \subseteq M$ is said to be an $X^\imath$-weight submodule if we have $N = \bigoplus_{\zeta \in X^\imath} (N \cap M_\zeta)$.
\end{defi}

The following is immediate from the definition.

\begin{prop}\label{sub and quot of Xi-weight module}
Let $M$ be an $X^\imath$-weight module, and $N \subseteq M$ an $X^\imath$-weight submodule.
\begin{enumerate}
\item\label{sub and quot of Xi-weight module 1} $N$ itself is an $X^\imath$-weight module such that $N_\zeta = N \cap M_\zeta$ for all $\zeta \in X^\imath$.
\item\label{sub and quot of Xi-weight module 2} The quotient $\Ui$-module $M/N$ admits an $X^\imath$-weight module structure such that $(M/N)_\zeta = M_\zeta/N = M_\zeta/N_\zeta$ for all $\zeta \in X^\imath$.
\end{enumerate}
\end{prop}

\begin{prop}\label{Ui-submodule generated by Xi-weight vectors is an Xi-weight submodule}
Let $M$ be an $X^\imath$-weight module, and $N \subseteq M$ a $\Ui$-submodule. Suppose that $N$ is generated by $X^\imath$-weight vectors of $M$. Then, $N$ is an $X^\imath$-weight submodule.
\end{prop}

\begin{proof}
As one can easily see, a sum of $X^\imath$-weight submodules of $M$ is an $X^\imath$-weight submodule of $M$. Hence, it suffices to show the assertion for the case when $N$ is generated by a single $X^\imath$-weight vector $m \in M_\eta$ for some $\eta \in X^\imath$. Let  $\zeta' \in Q^\imath$ and $x \in \Ui(\zeta')$. Then, by equation \eqref{commutation relations for B,E,k with 1zeta}, we have
$$
\mathbf{1}_\zeta x = x \mathbf{1}_{\zeta - \zeta'}.
$$
Hence, it follows that
$$
\mathbf{1}_\zeta (xm) = x \mathbf{1}_{\zeta-\zeta'} m = \delta_{\zeta-\zeta',\eta} xm \in N.
$$
Therefore, $N$ is a $\Uidot$-submodule of $M$. Hence, we have
$$
\mathbf{1}_\zeta N = N \cap M_\zeta.
$$
On the other hand, since $\mathbf{1}_\zeta$'s are orthogonal idempotents and $M$ is a unital $\Uidot$-module, we obtain
$$
N = \bigoplus_{\zeta \in X^\imath} \mathbf{1}_\zeta N.
$$
This proves the assertion.
\end{proof}

The modified quantum group $\Udot$ admits a $\Uidot$-bimodule structure given as follows: For each $x_1,x_2 \in \Ui$, $\zeta_i \in X^\imath$, and $y \in \Udot$, set
$$
(x_1 \mathbf{1}_{\zeta_1}) y (\mathbf{1}_{\zeta_2} x_2) := \sum_{\substack{\lm_i \in X \\ \ol{\lm_i} = \zeta_i}} (x_1 \mathbf{1}_{\lm_1}) y (\mathbf{1}_{\lm_2} x_2).
$$
In the right-hand side, we regard $x_i$ as an element of $\U$ via the inclusion $\Ui \hookrightarrow \U$.

Let $M$ be a unital $\Udot$-module. Then, it admits a weight module structure $M = \bigoplus_{\lm \in X} M_\lm$, and hence, the canonical $X^\imath$-weight module structure $M = \bigoplus_{\zeta \in X^\imath} M_\zeta$ (see Example \ref{canonical Xi weight module structure}). Therefore, it admits a unital $\Uidot$-module structure. This structure is compatible with the $\Uidot$-bimodule structure of $\Udot$.

More generally, let us consider a Satake subdiagram.

\begin{defi}\normalfont
  A Satake subdiagram of a Satake diagram $(I, I_\bullet, \tau)$ is a triple $(J, J_\bullet, \sigma)$ satisfying the following:
  \begin{itemize}
    \item $(J,J_\bullet,\sigma)$ itself is a Satake diagram.
    \item $J \subseteq I$, $J_\bullet \subseteq I_\bullet$, $\sigma(j) \in \{ j, \tau(j) \}$ for all $j \in J$.
    \item $w_\bullet(\alpha_{\tau(j)}) = w_\bullet(J)(\alpha_{\sigma(j)})$ for all $j \in J \cap I_\circ$.
  \end{itemize}
  Here and after, $A(J)$ denotes an object $A$ constructed from $(J,J_\bullet,\sigma)$ instead of $(I,I_\bullet,\tau)$.
\end{defi}

Let $(J, J_\bullet, \sigma)$ be a Satake subdiagram of our Satake diagram $(I, I_\bullet, \tau)$, and $\bfvarsigma(J), \bfkappa(J)$ parameters for $\mathbf{U}^\imath(J)$ such that $\varsigma_j = \varsigma(J)_j$, $\bfkappa_j = \bfkappa(J)_j$ for all $j \in J \cap I_\circ$.
We claim that $\mathbf{U}^\imath(J) \subseteq \mathbf{U}$.
In fact, for each $j \in J \cap I_\circ$, we have
$$
T_{w_\bullet}(E_{\tau(j)}) = T_{w_\bullet(J)} T_{w_\bullet(J)^{-1} w_\bullet}(E_{\tau(j)}) = T_{w_\bullet(J)}(E_{\sigma(j)}).
$$
The last equality follows from the identity $w_\bullet(J)^{-1} w_\bullet(\alpha_{\tau(j)}) = \alpha_{\sigma(j)}$ and \cite[Proposition 8.20]{J96}.
Hence, we obtain
$$
B_j(J) = F_j - \varsigma(J)_j T_{w_\bullet(J)}(E_{\sigma(j)})K_j^{-1} = F_j - \varsigma_j T_{w_\bullet}(E_{\tau(j)})K_j^{-1} = B_j \in \mathbf{U}^\imath.
$$
Similar statements for the other generators are verified more easily.

\begin{ex}\normalfont
Let $I$ be a Dynkin diagram, and consider two Satake diagrams $(I,\emptyset,\id)$ and $(I,I,-w_0)$, where $-w_0$ denotes the diagram automorphism induced by the longest element of the Weyl group of $I$. Namely, the first diagram consists of only black nodes, while the second one only white nodes with the same underlying Dynkin diagram $I$. Then, $(I,\emptyset,\id)$ is a Satake subdiagram of $(I,I,-w_0)$.
\end{ex}

For each $\lm \in X$, define $\lm' \in X(J)$ by
$$
\la h_j, \lm' \ra = \begin{cases}
-\hf \la w_\bullet(J)(h_j) - h_j, \lm \ra \qu & \IF j \in J_\circ \cap I_\bullet, \\
\la h_j,\lm \ra \qu & \OW.
\end{cases}
$$
By \cite[(5.6)]{Ko14}, for each $j \in J_\circ$, we have
$$
w_\bullet(J)(h_j) - h_j = w_\bullet(J)(h_{\sigma(j)}) - h_{\sigma(j)} \in \sum_{k \in J_\bullet} \Z h_k.
$$
Then, it follows that
$$
\la h_j, \lm + w_\bullet \circ \tau(\lm) \ra = \la h_j, \lm' + w_\bullet(J) \circ \sigma(\lm') \ra \qu \Forall j \in J.
$$
This shows that there exists a well-defined $\Z$-linear map $\cdot|_{\frh(J)} : X^\imath \rightarrow X^\imath(J)$ making the following diagram commute
$$
\xymatrix{
X \ar[r]^{\ol{\ \cdot\ }} \ar[d]_{\cdot|_{\frh(J)}} & X^\imath \ar[d]^{\cdot |_{\frh(J)}} \\
X(J) \ar[r]^{\ol{\ \cdot\ }} & X^\imath(J),
}
$$
where the left vertical arrow $\cdot|_{\frh(J)} : X \rightarrow X(J)$ denotes the ordinary restriction map from the Cartan subalgebra $\frh$ of $\g$ to the one $\frh(J)$ of $\g(J)$. Then, as before, $\Uidot$ admits a $\Uidot(J)$-bimodule structure: For each $x_1,x_2 \in \Ui(J)$, $\eta_i \in X^\imath(J)$, and $y \in \Uidot$, set
$$
(x_1 \mathbf{1}_{\eta_1}) y (\mathbf{1}_{\eta_2} x_2) = \sum_{\substack{\zeta_i \in X^\imath \\ \zeta_i|_{\frh(J)} = \eta_i}} (x_1 \mathbf{1}_{\zeta_1}) y (\mathbf{1}_{\zeta_2} x_2).
$$

Let $M$ be a unital $\Uidot$-module. Then, it admits an $X^\imath$-weight module structure $M = \bigoplus_{\zeta \in X^\imath} M_\zeta$, where $M_\zeta = \mathbf{1}_\zeta M$. By restriction, we obtain a $\Ui(J)$-module structure on $M$. Furthermore, it becomes an $X^\imath(J)$-weight module such that $M = \bigoplus_{\eta \in X^\imath(J)} M_\eta$, where
$$
M_\eta := \bigoplus_{\substack{\zeta \in X^\imath \\ \zeta|_{\frh(J)} = \eta}} M_\zeta.
$$
Therefore, $M$ admits a unital $\Uidot(J)$-module structure. Again, this structure is compatible with the $\Uidot(J)$-bimodule structure of $\Uidot$.

\subsection{Based $\Ui$-modules}
Recall that $\Uidot$ acts on each weight $\U$-module via the canonical $X^\imath$-weight module structure. Set
$$
\UidotA := \{ x \in \Uidot \mid x V(\lm)_{\bfA} \subseteq V(\lm)_{\bfA} \ \Forall \lm \in X^+ \}.
$$

\begin{defi}\label{A-form of Ui-modules}\normalfont
Let $M$ be an $X^\imath$-weight module. An $\bfA$-form of $M$ is an $\bfA$-lattice $M_{\bfA}$ of $M$ which is simultaneously a $\UidotA$-submodule.
\end{defi}

\begin{defi}\normalfont
Let $M$ be a $\Ui$-module equipped with a Hermitian inner product $(\cdot,\cdot)_M$. We say that $(\cdot,\cdot)_M$ is contragredient if it satisfies
\begin{align}\label{Property2}
(xu,v)_M = (u,\wp^*(x)v)_M \qu \Forall x \in \Ui,\ u,v \in M.
\end{align}
\end{defi}

\begin{rem}\normalfont
Each contragredient Hermitian inner product on a $\U$-module $M$ is a contragredient Hermitian inner product of the $\Ui$-module $M$.
\end{rem}

\begin{defi}\normalfont
Let $M$ be an $X^\imath$-weight module equipped with a contragredient Hermitian inner product $(\cdot,\cdot)_M$, an $\bfA$-form $M_{\bfA}$, an $\imath$bar-involution $\psii_M$, and a basis $\clB_M$ of $\ol{\clL}_M$. We say that $(M,(\cdot,\cdot)_M, M_{\bfA}, \psii_M, \clB_M)$, or simply $M$, is a based $\Ui$-module if the following two conditions are satisfied:
\begin{enumerate}
\item The quotient map $\ev_\infty : \clL_M \rightarrow \ol{\clL}_M$ restricts to an isomorphism
$$
\clL_M \cap M_{\bfA} \cap \psii_M(\clL_M) \rightarrow \ol{\clL}_M
$$
of $\C$-vector spaces; let $G^\imath : \ol{\clL}_M \rightarrow \clL_M \cap M_{\bfA} \cap \psii_M(\clL_M)$ denote its inverse.
\item For each $b \in \clB_M$, we have $\psii_M(G^\imath(b)) = G^\imath(b)$.
\end{enumerate}
\end{defi}

\begin{ex}\label{V(lm) with icanonical basis is a based Ui-module}\normalfont
Let $\lm \in X^+$, and consider the irreducible $\U$-module $V(\lm)$. Then, by Theorem \ref{icanonical basis of based Umodules}, the quintuple $(V(\lm), (\cdot,\cdot)_\lm, V(\lm)_{\bfA}, \psii_\lm, \clB(\lm))$ is a based $\Ui$-module.
\end{ex}

\begin{lem}\label{characterization of icanonical basis element}
Let $M$ be a based $\Ui$-module and $v \in \clL_M \cap M_{\bfA}$. Suppose that $\psii_M(v) = v$. Then, we have
$$
v = G^\imath(\ev_\infty(v)).
$$
In particular, if $v \in q\inv \clL_M$, then we have $v = 0$.
\end{lem}

\begin{proof}
By our assumption that $v = \psii_M(v) \in \psii_M(\clL_M)$, we have
$$
v \in \clL_M \cap M_{\bfA} \cap \psii_M(\clL_M).
$$
Since $G^\imath$ is the inverse of $\ev_\infty$, the assertion follows.
\end{proof}

\begin{prop}\label{almost orthonormality of global bases}
Let $M$ be a finite-dimensional based $\Ui$-module with $\clB_M$ being an orthonormal basis of $\ol{\clL}_M$. Then, $G^\imath(\clB_M)$ forms an almost orthonormal $\bfK$-basis of $M$, a free $\bfK_\infty$-basis of $\clL_M$, and a free $\bfA$-basis of $M_{\bfA}$.
\end{prop}

\begin{proof}
Since $\dim_{\C} \ol{\clL}_M = \dim_{\bfK} M$, $G^\imath(\clB_M)$ forms a $\bfK$-basis of $M$. For each $b,b' \in \clB_M$, we have
$$
\ev_\infty((G^\imath(b), G^\imath(b'))_M) = (\ev_\infty(G^\imath(b)), \ev_\infty(G^\imath(b')))_M = (b,b')_M = \delta_{b,b'}.
$$
This implies that $G^\imath(\clB_M)$ is almost orthonormal.
\end{proof}

\begin{defi}\normalfont
Let $M,N$ be based $\Ui$-modules and $f \in \Hom_{\Uidot}(M,N)$.
\begin{enumerate}
\item We say that $f$ is a homomorphism of based ($\Ui$-)modules if
$$
f(G^\imath(b)) \in G^\imath(\clB_N) \sqcup \{0\}
$$
for all $b \in \clB_M$.
\item We say that $M$ is a based submodule of $N$ if $f$ is the inclusion map and a homomorphism of based modules.
\item We say that $N$ is a based quotient module of $M$ if $f$ is the quotient map and a homomorphism of based modules.
\end{enumerate}
\end{defi}

\begin{prop}\label{based submodule and based quotient module}
Let $M$ be a finite-dimensional based $\Ui$-module with $\clB_M$ being an orthonormal basis, and $N \subseteq M$ a $\Uidot$-submodule. Suppose that $N = \bfK G^\imath(\clB_N)$, where $\clB_N := \clB_M \cap \ol{\clL}_N$.
\begin{enumerate}
\item\label{based submodule and based quotient module 1} Set $N_{\bfA} := N \cap M_{\bfA}$, and $\psii_N := \psii_M|_N$. Then, $(N, (\cdot,\cdot)_N, N_{\bfA}, \psii_N, \clB_N)$ is a based submodule of $M$.
\item\label{based submodule and based quotient module 2} Let $\psii_{M/N}$ denote the $\imath$bar-involution on $M/N$ induced from $\psii_M$, and $\clB_{M/N} := \{ b + \ol{\clL}_N \mid b \in \clB_M {\setminus} \clB_N \}$. Then, $(M/N, (\cdot,\cdot)_{M/N}, M_{\bfA}/N_{\bfA}, \psii_{M/N}, \clB_{M/N})$ is a based quotient module of $M$.
\end{enumerate}
\end{prop}

\begin{proof}
By Proposition \ref{almost orthonormality of global bases}, $G^\imath(\clB_N)$ and $G^\imath(\clB_M)$ form $\psii_M$-invariant almost orthonormal free $\bfA$-bases of $N_{\bfA}$ and $M_{\bfA}$, respectively. Furthermore, $\{ [G^\imath(b)] \mid b \in \clB_M \setminus \clB_N \}$ forms a $\psii_{M/N}$-invariant almost orthonormal free $\bfA$-basis of $M_{\bfA}/N_{\bfA}$. Then, the assertion follows if we prove that $N_{\bfA}$ is a $\UidotA$-submodule.

Let $b \in \clB_N$ and $x \in \UidotA$. Then, we have $x G^\imath(b) \in N \cap M_{\bfA} = N_{\bfA}$. Thus, the proof completes.
\end{proof}

\section{$\imath$Quantum group of type AI}\label{sectoin: QSP of type AI}
In the remainder, we restrict our attention to the $\imath$quantum group of type AI with parameters $\varsigma_i = q\inv$ and $\kappa_i = 0$ for all $i \in I_\circ$. It is an $\imath$quantum group associated with Satake diagram $(I,\emptyset,\id)$, where $I = [1,n-1]$ is the Dynkin diagram of type $A_{n-1}$. We assume that $n \geq 3$ unless otherwise specified. The case when $n = 2$ will be studied separately because $\frk$ is not a simple Lie algebra in this case.

Let $\g = \g(I)$ denote the complex simple Lie algebra associated with $I$. Recall that $Y = \bigoplus_{i=1}^{n-1} \Z h_i$ and $X = \bigoplus_{i=1}^{n-1} \Z \vpi_i$ are the coroot lattice and the weight lattice for $\g$. We often identify $\lm = \sum_{i=1}^{n-1} \lm_i \vpi_i \in X$ with $(\lm_1,\ldots,\lm_{n-1}) \in \Z^{n-1}$. For each $\lm = (\lm_1,\ldots,\lm_{n-1}) \in X$, we have $\lm \in X^+$ if and only if $\lm_i \geq 0$ for all $i \in I$.

\subsection{Symmetric pair of type AI}\label{subsection: symmetric pair of type AI}
Let $\frk$ be the subalgebra of $\g$ generated by $b_i := f_i + e_i$, $i \in I$. The defining relations of $\frk$ for the generators $b_1,\ldots,b_{n-1}$ are as follows (see \cite[Definition 2.5 and Theorem 2.7]{St20}):
\begin{align}\label{Defining Relations for k}
\begin{split}
&[b_i,b_j] = 0 \qu \IF |i-j| > 1, \\
&[b_i,[b_i,b_j]] = b_j \qu \IF |i-j| = 1.
\end{split}
\end{align}

Let us identify $\g$ with the special linear Lie algebra $\frsl_n := \{ X \in \Mat_n(\C) \mid \tr X = 0 \}$ in the usual way. Consider the special orthogonal Lie algebra $\frso_n := \{ X \in \g \mid X + {}^t X = 0 \} \subseteq \g$. It is generated by $b'_i := f_i - e_i$, $i \in I$. Let $m$ denote the rank of $\frso_n$:
$$
m = \begin{cases}
\frac{n}{2} \qu & \IF n \in \Z_{\ev}, \\
\frac{n-1}{2} \qu & \IF n \in \Z_{\odd}.
\end{cases}
$$
Then, $\bigoplus_{i=1}^m \C b'_{2i-1}$ forms a Cartan subalgebra of $\frso_n$ since it forms an abelian subalgebra of dimension $m$.

Set $I_{\frk} := [1,m]$, and identify it with the Dynkin diagram of type $D_m$ if $n \in \Z_{\ev}$ or $B_m$ if $n \in \Z_{\odd}$ whose vertices are labeled by $I_\frk$ in the same manner as \cite[Section 11.4]{H72}. Let $e'_i,f'_i,h'_i$, $i \in I_{\frk}$ be the Chevalley generators of $\g(I_{\frk})$.

Let $I_\otimes := I \cap \Z_{\odd} = \{ 2i-1 \mid i \in I_{\frk} \}$, and set $\frh_{\frk}$ to be the subalgebra of $\frk$ generated by $b_i$, $i \in I_\otimes$. As before, we see that $\frh_{\frk}$ is a Cartan subalgebra of $\frk$. Let $\{ b^i \mid i \in I_\otimes \} \subseteq \frh_{\frk}^*$ denote the dual basis of $\{ b_i \mid i \in I_\otimes \}$. For each $i \in I_{\frk}$, set
$$
\gamma_i := \begin{cases}
b^{2i-1} -  b^{2i+1} \qu & \IF i \neq m, \\
b^{2m-3} + b^{2m-1} \qu & \IF i = m \AND n \in \Z_{\ev}, \\
b^{2m-1} \qu & \IF i = m \AND n \in \Z_{\odd}.
\end{cases}
$$

For each $i \in I_{\frk}$ such that $2i < n$, set
$$
b_{2i,\pm} := \frac{1}{2}(b_{2i} \pm [b_{2i-1},b_{2i}]).
$$
Also, for each $i \in I_{\frk}$ such that $2i < n-1$ and for each $e \in \{ +,- \}$, set
$$
b_{2i,e\pm} := \frac{1}{2}(b_{2i,e} \pm [b_{2i+1},b_{2i,e}]).
$$
Now, for each $i \in I_{\frk}$, set
\begin{align}
\begin{split}
&x_i := \begin{cases}
2b_{2i,+-} \qu & \IF i \neq m, \\
2b_{2m-2,++} \qu & \IF i = m \AND n \in \Z_{\ev}, \\
2b_{2m,+} \qu & \IF i = m \AND n \in \Z_{\odd},
\end{cases} \\
&y_i := \begin{cases}
2b_{2i,-+} \qu & \IF i \neq m, \\
2b_{2m-2,--} \qu & \IF i = m \AND n \in \Z_{\ev}, \\
2b_{2m,-} \qu & \IF i = m \AND n \in \Z_{\odd},
\end{cases} \\
&w_i := \begin{cases}
b_{2i-1}-b_{2i+1} \qu & \IF i \neq m, \\
b_{2m-3}+b_{2m-1} \qu & \IF i = m \AND n \in \Z_{\ev}, \\
2b_{2m-1} \qu & \IF i = m \AND n \in \Z_{\odd}.
\end{cases}
\end{split} \nonumber
\end{align}
Using relations \eqref{Defining Relations for k}, one can verify that there exists an isomorphism $\g(I_{\frk}) \rightarrow \frk$ of Lie algebras sending $e'_i,f'_i,h'_i$ to $x_i,y_i,w_i$, respectively (see also \cite[Sections 4.1--4.2]{W19}).

Define $x'_i,y'_i,w'_i \in \frso_n$, $i \in I_{\frk}$ in the same way as $x_i,y_i,w_i$ by replacing $b_i$ with $\sqrt{-1}b'_i$. Then, there exists an isomorphism $\g(I_{\frk}) \rightarrow \frso_n$ of Lie algebras sending $e'_i,f'_i,h'_i$ to $x'_i,y'_i,w'_i$, respectively. However, one should note that $\frso_n \neq \frk$ set-theoretically.

Set $Y_{\frk} := \bigoplus_{i \in I_{\frk}} \Z w_i$, $X_{\frk} := \Hom_{\Z}(Y_{\frk},\Z)$, $X_{\frk}^+ := \{ \nu \in X_{\frk} \mid \la w_i,\nu \ra \geq 0 \Forall i \in I_{\frk} \}$. Then, the isoclasses of finite-dimensional irreducible $\frk$-modules are parametrized by $X_{\frk}^+$. Let $V_{\frk}(\nu)$ denote the irreducible $\frk$-module of highest weight $\nu \in X_{\frk}^+$.

We often identify $\nu \in X_{\frk}$ with $(\nu_1,\nu_3,\ldots,\nu_{2m-1}) \in (\hf\Z)^{I_\otimes}$, where
$$
\nu_i := \la b_i,\nu \ra.
$$
Here, we extend $\nu$ to a linear form on $\frh_{\frk} = Y_{\frk} \otimes_{\Z} \C$. With this notation, we have
$$
\la w_i,\nu \ra = \begin{cases}
\nu_{2i-1} - \nu_{2i+1} \qu & \IF i \neq m, \\
\nu_{2m-3} + \nu_{2m-1} \qu & \IF i = m \AND n \in \Z_{\ev}, \\
2\nu_{2m-1} \qu & \IF i = m \AND n \in \Z_{\odd}.
\end{cases}
$$
Note that we have either $\nu \in \Z^{I_\otimes}$ or $\nu \in (\hf + \Z)^{I_\otimes}$. We call $\nu \in \Z^{I_\otimes}$ an integer weight. Set $X_{\frk,\Int}$ to be the set of integer weights. Let $\nu = (\nu_1,\nu_3,\ldots,\nu_{2m-1}) \in X_{\frk}$. When $n \in \Z_{\ev}$, we have $\nu \in X_{\frk}^+$ if and only if $\nu_1 \geq \cdots \geq \nu_{2m-3} \geq |\nu_{2m-1}|$, while when $n \in \Z_{\odd}$, we have $\nu \in X_{\frk}^+$ if and only if $\nu_1 \geq \cdots \geq \nu_{2m-1} \geq 0$.

\subsection{$\imath$Quantum group of type AI}
Let $\Ui$ denote the $\imath$quantum group of type AI with parameters $\varsigma_i = q\inv$, $\kappa_i = 0$ for all $i \in I$. It is a subalgebra of $\U$ generated by $B_i = F_i + q\inv E_iK_i\inv$, $i \in I$. The defining relations of $\Ui$ for $B_1,\ldots,B_{n-1}$ are as follows (see e.g. \cite[Section 7]{Ko14}):
\begin{align}\label{Defining Relations for Ui}
\begin{split}
&[B_i,B_j] = 0 \qu \IF |i-j| > 1, \\
&[B_i,[B_i,B_j]_q]_{q\inv} = B_j \qu \IF |i-j| = 1.
\end{split}
\end{align}

Let us recall the notion of classical weight $\Ui$-modules, which was introduced in \cite{W19}.

\begin{defi}\normalfont
Let $M$ be a $\Ui$-module. $M$ is said to be a classical weight module if for each $i \in I_\otimes$, $B_i$ acts on $M$ diagonally with eigenvalues belonging to $\bfK_1$.
\end{defi}

Let $M$ be a classical weight module. For each $\nu \in \frh_{\frk}^*$, set $M_\nu$ to be the subspace of $M$ consisting of $v \in M$ satisfying the following: For each $i \in I_\otimes$, $B_i v = av$ for some $a \in \bfK_1$ such that $\ev_{1}(a) = \la b_i, \nu \ra$. Then, we have by definition
$$
M = \bigoplus_{\nu \in \frh_{\frk}^*} M_\nu.
$$
The subspace $M_\nu$ (resp., an element of $M_\nu$) is referred to as the $X_{\frk}$-weight space of $M$ (resp., an $X_{\frk}$-weight vector) of weight $\nu$.

\begin{prop}[{\cite[Proposition 3.1.4 and Remark 3.1.5]{W19}}]
Let $M$ be a $\U$-module of finite dimension. Then, as a $\Ui$-module, it is a classical weight module. Moreover, for each $i \in I_\otimes$, the eigenvalues of $B_i$ are of the form $[a]$, $a \in \Z$.
\end{prop}

\begin{prop}[{\cite[Proposition 3.1.6]{W19}}]\label{classical weight tensor fd is classical weight}
Let $M$ be a classical weight module, and $N$ a finite-dimensional $\U$-module. Then, $M \otimes N$ is a classical weight module.
\end{prop}

\begin{rem}\normalfont
A weight vector of a finite-dimensional $\U$-module $M$ is not necessarily an $X_{\frk}$-weight vector of $M$. Similarly, the tensor product of an $X_{\frk}-$weight vector of a classical weight module $M$ and a weight vector of a finite-dimensional $\U$-module $N$ is not necessarily an $X_{\frk}$-weight vector of $M \otimes N$.
\end{rem}

In \cite{W19}, we defined linear operators $l_j^{\pm 1}$ for $j \in I_\otimes$, and $B_{i,e}$, $B_{i,e_1e_2}$ for $i \in I \setminus I_\otimes$, $e,e_1,e_2 \in \{ +,- \}$ acting on each classical weight module $M$. By definition, for each $v \in M$ and $e \in \{ +,- \}$, we have
\begin{align}
  \begin{split}
    &l_j^{\pm 1} v = q^{\pm a} v \qu \IF B_j v = [a]v, \\
    &B_{2i,\pm} = (B_{2i}l_{2i-1}^{\pm 1} \pm [B_{2i-1},B_{2i}]_q) \frac{1}{\{ l_{2i-1};0 \}}, \\
    &B_{2i,e\pm} = (B_{2i,e}l_{2i+1}^{\pm 1} \pm [B_{2i+1},B_{2i,e}]_q) \frac{1}{\{ l_{2i+1};0 \}}.
  \end{split} \nonumber
\end{align}
Set
$$
\Itil_{\frk} := \begin{cases}
\{ (2i,\pm) \mid i \in [1,m-1] \} \qu & \IF n \in \Z_{\ev}, \\
\{ (2i,\pm) \mid i \in [1,m-1] \} \sqcup \{ 2m \} \qu & \IF n \in \Z_{\odd}.
\end{cases}
$$
For each $j \in \Itil_{\frk}$, set
\begin{align}
\begin{split}
&X_j := \begin{cases}
B_{2i,+\pm}\{ l_{2i-1};0 \} \qu & \IF j = (2i,\pm), \\
B_{2m,+}\{ l_{2m-1} \} \qu & \IF j = 2m,
\end{cases} \\
&Y_j := \begin{cases}
B_{2i,-\mp}\{ l_{2i+1};0 \} \qu & \IF j = (2i,\pm), \\
B_{2m,-}\{ l_{2m-1} \} \qu & \IF j = 2m,
\end{cases} \\
&\gamma_j := \begin{cases}
b^{2i-1} \pm b^{2i+1} \qu & \IF j = (2i,\pm), \\
b^{2m-1} \qu & \IF j = 2m.
\end{cases}
\end{split} \nonumber
\end{align}
For each $\nu,\xi \in X_{\frk}$, we write $\xi \leq \nu$ to indicate that $\nu - \xi \in \sum_{i \in I_{\frk}} \Z_{\geq 0} \gamma_i = \sum_{j \in \Itil_{\frk}} \Z_{\geq 0} \gamma_j$. This defines a partial order on $X_{\frk}$, called the dominance order.

\begin{theo}[{\cite[Corollaries 4.1.4, 4.1.7, 4.2.2, and 4.2.4]{W19}}]\label{facts about classical weight Ui-modules}
\ \begin{enumerate}
\item For each $i \in I_\otimes$, $l_i$ (resp., $B_i$) acts on each finite-dimensional classical weight module diagonally with eigenvalues of the form $q^a$ (resp., $[a]$), $a \in \Z$.
\item For each classical weight module $M$, $\nu \in X_{\frk}$, and $j \in \Itil_{\frk}$, we have
$$
X_j M_\nu \subseteq M_{\nu + \gamma_j}, \qu Y_j M_\nu \subseteq M_{\nu-\gamma_j}.
$$
\item Each finite-dimensional classical weight module is completely reducible.
\item The isoclasses of finite-dimensional irreducible classical weight modules are para\-metrized by the set $X_{\frk,\Int}^+$: The finite-dimensional irreducible classical weight module $V(\nu)$ corresponding to $\nu \in X_{\frk,\Int}^+$ is characterized by the condition that there exists a unique (up to nonzero scalar multiple) $v_\nu \in V(\nu)_\nu \setminus \{0\}$ such that $X_j v_\nu = 0$ for all $j \in \Itil_{\frk}$; the vector $v_\nu$ is referred to as a highest weight vector.
\item\label{spanning property for AI} For each $\nu \in X_{\frk,\Int}^+$, we have
$$
V(\nu) = \Span_{\bfK} \{ Y_{j_1} \cdots Y_{j_r} v_\nu \mid r \geq 0,\ j_1,\ldots,j_r \in \Itil_{\frk} \}.
$$
\item For each $\nu \in X_{\frk,\Int}^+$ and $\xi \in X_{\frk}$, we have $\dim V(\nu)_{\nu} = 1$, and $\dim V(\nu)_\xi = 0$ unless $\xi \leq \nu$.
\item\label{facts about classical weight Ui-modules 7} For each $\nu \in X_{\frk,\Int}^+$ and $\xi \in X_{\frk}$ such that $\xi < \nu$, we have
$$
V(\nu)_\xi = \sum_{j \in \Itil_{\frk}} Y_j V(\nu)_{\xi + \gamma_j}.
$$
\end{enumerate}
\end{theo}

Given a finite-dimensional classical weight module $M$ and $\nu = (\nu_1,\nu_3,\ldots,\nu_{2m-1}) \in X_{\frk,\Int}$, the weight space $M_\nu$ is described as
$$
M_{\nu} = \{ v \in M \mid l_{i} v = q^{\nu_i} v \Forall i \in I_\otimes \} = \{ v \in M \mid B_i v = [\nu_i]v \Forall i \in I_\otimes \}.
$$
Define the character $\ch_{\Ui} M \in \Z[X_{\frk,\Int}]$ of $M$ by
$$
\ch_{\Ui} M := \sum_{\nu \in X_{\frk,\Int}} (\dim_{\bfK} M_\nu) e^\nu.
$$

\begin{prop}
Let $\nu \in X_{\frk,\Int}^+$. Then, we have
$$
\ch_{\Ui} V(\nu) = \ch_{\frk} V_{\frk}(\nu),
$$
where the right-hand side denotes the character of the irreducible $\frk$-module $V_{\frk}(\nu)$ of highest weight $\nu$.
\end{prop}

\begin{proof}
The assertion essentially follows from \cite[arguments before Proposition 3.3.9]{W19}.
\end{proof}

Let us recall from \cite{W19} some formulas which will be frequently used below.
\begin{align}
&B_{2i-1} = [l_{2i-1};0], \label{formula for AI 1}\\
&B_{2i,\pm} = (B_{2i}l_{2i-1}^{\pm 1} \pm [B_{2i-1},B_{2i}]_q) \frac{1}{\{ l_{2i-1};0 \}}, \label{formula for AI 11}\\
&B_{2i} = B_{2i,+} + B_{2i,-}, \label{formula for AI 2}\\
&l_{2j-1} B_{2i,\pm} = q^{\pm \delta_{i,j}} B_{2i,\pm} l_{2j-1} \qu \IF j \neq i+1, \label{formula for AI 5}\\
&[B_{2i,+}\{ l_{2i-1};0 \}, B_{2i,-}\{ l_{2i-1};0 \}] = [l_i^2;0], \label{formula for AI 7}\\
&B_{2i,e\pm} = (B_{2i,e}l_{2i+1}^{\pm 1} \pm [B_{2i+1},B_{2i,e}]_q) \frac{1}{\{ l_{2i+1};0 \}}, \label{formula for AI 12}\\
&B_{2i,\pm} = B_{2i,\pm +} + B_{2i,\pm -}, \label{formula for AI 3}\\
&l_{2j-1} B_{2i,e_1e_2} = q^{e_1 \delta_{i,j} + e_2 \delta_{i,j-1}} B_{2i,e_1e_2} l_{2j-1}, \label{formula for AI 6}\\
&[B_{2i,+\pm}\{ l_{2i-1};0 \}, B_{2i,-\mp}\{ l_{2i+1};0 \}] = (1 + (q-q\inv)^2 B_{2i,-\pm} B_{2i,+\mp})[l_{2i-1}l_{2i+1}^{\pm 1};0], \label{formula for AI 8}\\
&[B_{2i,+\pm} \{ l_{2i-1};0 \}, B_{2i,-\pm} \{ l_{2i-1};0 \}] = 0, \label{formula for AI 9}\\
&[B_{2i,\pm+} \{ l_{2i+1};0 \}, B_{2i,\pm-} \{ l_{2i+1};0 \}] = 0. \label{formula for AI 10}
\end{align}

\begin{lem}\label{wp for AI}
Let $M$ be a finite-dimensional classical weight module equipped with a contragredient Hermitian inner product. Then, for each $u,v \in M$, $i \in I_{\frk}$, we have (if $B_{2i,\pm}$, etc. is defined)
\begin{align}
\begin{split}
&(l_{2i-1} u,v) = (u, l_{2i-1} v), \\
&(B_{2i,\pm} u,v) = (u, B_{2i,\mp} v), \\
&(B_{2i,+\pm} u,v) = (u, B_{2i,-\mp}v), \\
&(B_{2i,-\pm} u,v) = (u, B_{2i,+\mp}v).
\end{split} \nonumber
\end{align}
\end{lem}

\begin{proof}
The first assertion is reduced to the case when both $u$ and $v$ are $X_{\frk}$-weight vectors. The latter is clear from the definitions.

Applying Proposition \ref{wp on Ui} to our case, we have
$$
\wp^*(B_i) = B_i
$$
for all $i \in I$. Then, the remaining assertions are directly verified by using the first assertion and equations \eqref{formula for AI 1}--\eqref{formula for AI 6}.
\end{proof}

\begin{lem}\label{ibar-involution on l,Bpm}
Let $M$ be a finite-dimensional classical weight module equipped with an $\imath$bar-involution $\psii_M$. Then, for each $v \in M$, $i \in I_{\frk}$, we have (if $B_{2i,\pm}$, etc. is defined)
\begin{align}
\begin{split}
&\psii_M(l_{2i-1} v) = l_{2i-1}\inv \psii_M(v), \\
&\psii_M(B_{2i,\pm} v) = B_{2i,\pm} \psii_M(v), \\
&\psii_M(B_{2i,+\pm} v) = B_{2i,+\pm} \psii_M(v), \\
&\psii_M(B_{2i,-\pm} v) = B_{2i,-\pm} \psii_M(v).
\end{split} \nonumber
\end{align}
\end{lem}

\begin{proof}
Similar to the proof of Lemma \ref{wp for AI}.
\end{proof}

\subsection{Contragredient Hermitian inner products and $\imath$bar-involutions on irreducible modules}
In this subsection, we construct a contragredient Hermitian inner product and an $\imath$bar-involution on $V(\nu)$, $\nu \in X_{\frk,\Int}$.

Let $V_\natural$ denote the natural representation of $\U$. Namely, it possesses a basis $\{u_1,\ldots,u_n\}$ such that
$$
E_i u_j = \delta_{i,j-1} u_i, \qu F_i u_j := \delta_{i,j} u_{i+1}, \qu K_i u_j := q^{\delta_{i,j} - \delta_{i,j-1}} u_j.
$$
The natural representation $V_\natural$ is equipped with a contragredient Hermitian inner product $(\cdot,\cdot)_\natural$ given by
$$
(u_i,u_j)_{\natural} = \delta_{i,j}.
$$
Set
$$
\clL_\natural := \clL_{V_{\natural}}, \qu \ol{\clL}_\natural := \ol{\clL}_{V_\natural}, \qu \clB_\natural := \{ \ol{u}_1,\ldots,\ol{u}_n \},
$$
where $\ol{u}_i := \ev_\infty(u_i)$. Note that $\clB_\natural$ forms a crystal basis of $V_\natural$ on which the Kashiwara operators act as
$$
\Etil_i \ol{u}_j = \delta_{i,j-1} \ol{u}_i, \qu \Ftil_i \ol{u}_j = \delta_{i,j} \ol{u}_{i+1}.
$$
The following follows from a direct calculation.

\begin{lem}\label{X_k-weight vectors in M tensor Vnatural}
Let $M$ be a $\Ui$-module, $v \in M$. Suppose that $B_i v = [a_i] v$ for some $a_i \in \Z$ for all $i \in I_\otimes$. Let $i \in I_\otimes$, and set $v'_{i,\pm} := v \otimes (u_i \pm q^{\pm a_i} u_{i+1})$. Also, set $v_{0} := v \otimes u_n$ if $n \in \Z_{\odd}$. Then, for each $j \in I_\otimes$, we have
$$
B_j v'_{i,\pm} = [a_j \pm \delta_{i,j}] v'_{i,\pm}, \qu B_j v_{0} = [a_j] v_0.
$$
\end{lem}

\begin{prop}\label{spectra of Vnatural tensor}
Let $d \geq 0$ and $\nu \in X_{\frk,\Int}$. Set $M := V_\natural^{\otimes d}$.
\begin{enumerate}
\item If $n \in \Z_{\ev}$, then, $\dim M_\nu$ equals the number of sequences
$$
(e_1,\ldots,e_d) \in \{ (i,+),(i,-) \mid i \in I_\otimes \}^d
$$
such that $\sharp\{ k \mid e_k = (i,+) \} - \sharp\{ k \mid e_k = (i,-) \} = \nu_i$ for all $i \in I_\otimes$.
\item If $n \in \Z_{\odd}$, then, $\dim M_\nu$ equals the number of sequences
$$
(e_1,\ldots,e_d) \in \{ (i,+),(i,-),0 \mid i \in I_\otimes \}^d
$$
such that $\sharp\{ k \mid e_k = (i,+) \} - \sharp\{ k \mid e_k = (i,-) \} = \nu_i$ for all $i \in I_\otimes$.
\end{enumerate}
\end{prop}

\begin{proof}
Using Lemma \ref{X_k-weight vectors in M tensor Vnatural}, the assertions follow by induction on $d$.
\end{proof}

Let $d \geq 0$ and set $M := V_\natural^{\otimes d}$. Let $(\cdot,\cdot)_M$ denote the contragredient Hermitian inner product on $M$ obtained from $(\cdot,\cdot)_\natural$ by means of Lemma \ref{Tensor product of Hermitian form}. Set
$$
M_{\bfA_1} := \Span_{\bfA_1} \{ u_{i_1} \otimes \cdots \otimes u_{i_d} \mid i_1,\ldots,i_d \in [1,n] \}.
$$
Then, $M_1 := M_{\bfA_1}/(q-1)M_{\bfA_1}$ is a $\g (= \frsl_n)$-module isomorphic to $(\C^n)^{\otimes d}$, where $\C^n$ is equipped with the natural $\frsl_n$-module structure. As a $\frk$-module, we have an irreducible decomposition
$$
M_1 \simeq \bigoplus_{\nu \in X_{\frk,\Int}^+} V_{\frk}(\nu)^{\oplus m_\nu}
$$
for some $m_\nu \geq 0$.

On the other hand, by induction on $d$ and Lemma \ref{X_k-weight vectors in M tensor Vnatural}, we see that each $X_{\frk}$-weight space $M_\nu$, $\nu \in X_{\frk}$ is spanned by vectors in $M_{\bfA_1}$. Therefore, we obtain
$$
\ch_{\Ui} M = \ch_{\frk} M_1 = \sum_\nu m_\nu \ch_{\frk} V_{\frk}(\nu).
$$

Furthermore, we have an irreducible decomposition
$$
M \simeq \bigoplus_{\nu \in X_{\frk,\Int}^+} V(\nu)^{\oplus m'_\nu}
$$
as a $\Ui$-module for some $m'_\nu \geq 0$. Then, we have
$$
\sum_\nu m_\nu \ch_{\frk} V_{\frk}(\nu) = \ch_{\Ui} M = \sum_\nu m'_\nu \ch_{\Ui} V(\nu) = \sum_\nu m'_\nu \ch_{\frk} V_{\frk}(\nu).
$$
Since the irreducible characters of $\frk$ are linearly independent, we conclude that $m_\nu = m'_\nu$ for all $\nu \in X_{\frk,\Int}^+$.

\begin{prop}\label{Tensor power contains all highest weight vectors}
Let $\nu \in X_{\frk,\Int}^+$. Then, there exists $d \geq 0$ such that $V_\natural^{\otimes d}$ possesses a highest weight vector of weight $\nu$.
\end{prop}

\begin{proof}
It is a classical result that the irreducible $\frk$-module $V_{\frk}(\nu)$ of highest weight $\nu \in X_{\frk,\Int}^+$ appears as a submodule of $(\C^n)^{\otimes d}$ for some $d \geq 0$. This implies, by the argument above, that the irreducible $\Ui$-module $V(\nu)$ appears as a submodule of $V_\natural^{\otimes d}$.
\end{proof}

\begin{cor}
Let $\nu \in X_{\frk,\Int}^+$.
\begin{enumerate}
\item There exists a unique contragredient Hermitian inner product $(\cdot,\cdot)_\nu$ on $V(\nu)$ such that $(v_\nu,v_\nu)_\nu = 1$.
\item There exists a unique $\imath$bar-involution $\psii_\nu$ on $V(\nu)$ such that $\psii_\nu(v_\nu) = v_\nu$.
\end{enumerate}
\end{cor}

\begin{proof}
Let $d \geq 0$ be such that $M := V_\natural^{\otimes d}$ possesses a highest weight vector $v$ of weight $\nu$ (see Proposition \ref{Tensor power contains all highest weight vectors}). By Lemma \ref{ibar-involution on l,Bpm}, we see that $\psii_M(v)$ is also a highest weight vector of weight $\nu$. Replacing $v$ with $v + \psii_M(v)$, we may assume that $\psii_M(v) = v$. Furthermore, replacing $v$ with $\frac{1}{\lt(v)} v$, we may assume that $v$ is almost normal. Let $\phi : V(\nu) \rightarrow M$ be an $\Ui$-isomorphism such that $\phi(v_\nu) = v$. Define maps $(\cdot,\cdot)_\nu : V(\nu) \times V(\nu) \rightarrow \bfK$ and $\psii_\nu : V(\nu) \rightarrow V(\nu)$ by
\begin{align}
\begin{split}
&(u_1,u_2)_\nu := (v,v)_M\inv \cdot (\phi(u_1),\phi(u_2))_M, \\
&\psii_\nu(u) := \phi\inv(\psii_M(\phi(u))).
\end{split} \nonumber
\end{align}
Then, it is straightforwardly verified that these are a desired Hermitian inner product and an $\imath$bar-involution.

The uniqueness is easily verified. Hence, the proof completes.
\end{proof}

Let $\nu \in X_{\frk,\Int}^+$. As in the quantum groups setting, with respect to the contragredient Hermitian inner product $(\cdot,\cdot)_\nu$, we write
$$
\clL(\nu) := \clL_{V(\nu)}, \qu \ol{\clL}(\nu) := \ol{\clL}_{V(\nu)}, \qu b_\nu := \ev_\infty(v_\nu).
$$

\subsection{$\imath$Divided powers}
Following \cite{BeW18}, for each $i \in I$ and $k \in \Z_{\geq 0}$, set
\begin{align}
\begin{split}
&B_{i,\ev}^{(k)} := \begin{cases}
\frac{1}{[2a]!} B_i \prod_{b=-a+1}^{a-1} (B_i - [2b]) \qu & \IF k = 2a \in \Z_{\ev}, \\
\frac{1}{[2a+1]!} \prod_{b=-a}^{a} (B_i - [2b]) \qu & \IF k = 2a+1 \in \Z_{\odd}, \\
\end{cases} \\
&B_{i,\odd}^{(k)} := \begin{cases}
\frac{1}{[2a]!} \prod_{b=-a+1}^a (B_i - [2b-1]) \qu & \IF k = 2a \in \Z_{\ev}, \\
\frac{1}{[2a+1]!} B_i \prod_{b=-a+1}^a (B_i - [2b-1]) \qu & \IF k = 2a+1 \in \Z_{\odd}.
\end{cases}
\end{split} \nonumber
\end{align}

\begin{prop}[{\cite[equations (2.5) and (3.2)]{BeW18}}]\label{recursive formula for idivided powers}
For each $p \in \{ \ev,\odd \}$ and $k \geq 0$, we have
$$
B_i B_{i,p}^{(k)} = [k+1] B_{i,p}^{(k+1)} + \delta_{p,p(k)} [k] B_{i,p}^{(k-1)},
$$
where we understand $B_{i,p}^{(-1)} = 0$.
\end{prop}

Recall $X^\imath$ from equation \eqref{Xi and Yi} in page \pageref{Xi and Yi}. In our setting, we have
$$
X^\imath = X/2X = (\Z/2\Z)^{n-1}.
$$
For each $\zeta = (\zeta_1,\ldots,\zeta_{n-1}) \in X^\imath$, set
$$
B_{i,\zeta}^{(k)} := B_{i,p(\zeta_i)}^{(k)} \mathbf{1}_{\zeta}.
$$
These elements are called $\imath$divided powers.

\begin{theo}[{\cite[Corollary 7.5]{BW18c}}]
$\Uidot_{\bfA}$ is generated by $B_{i,\zeta}^{(k)}$, $i \in I$, $\zeta \in X^{\imath}$, $k \in \Z_{\geq 0}$.
\end{theo}

\begin{defi}\normalfont
Let $M = \bigoplus_{\zeta \in X^\imath} M_\zeta$ be an $X^\imath$-weight module. We say that $M$ is standard if for each $\zeta = (\zeta_1,\ldots,\zeta_{n-1}) \in X^\imath$, we have
$$
M_\zeta = \{ v \in M \mid B_{i,p(\zeta_i)}^{(k_i)}v = 0 \Forsome k_i > 0 \text{ such that } \ol{k_i} \neq \zeta_i, \Forall i \in I \}.
$$
\end{defi}

The following two propositions are immediate from the definition.

\begin{prop}\label{Xi-weight submodule of a standard Xi-weight module is standard}
Let $M$ be a standard $X^\imath$-weight module, and $N \subseteq M$ an $X^\imath$-weight submodule. Then, $N$ is standard.
\end{prop}

\begin{prop}\label{Ui-homomorphism lifts to Uidot-homomorphism}
Let $M$ and $N$ be standard $X^\imath$-weight modules. Then, each $\Ui$-module homomorphism $f : M \rightarrow N$ preserves the $X^\imath$-weight spaces. In particular, it lifts to a $\Uidot$-module homomorphism.
\end{prop}

The rest of this subsection is devoted to proving that each finite-dimensional $\U$-module with the canonical $X^\imath$-weight module structure is standard.

\begin{lem}\label{characterization of standard Xi weight}
Let $M$ be a $\Ui$-module. Let $v \in M$, $i \in I$, and $p \in \{ \ev,\odd \}$. Then, the following are equivalent:
\begin{enumerate}
\item $B_{i,p}^{(k)}v = 0 \Forsome k > 0 \text{ such that } p(k) \neq p$.
\item $v$ is a sum of $B_i$-eigenvectors of eigenvalues of the form $[a]$ with $a \in [-k+1,k-1]_{p}$.
\end{enumerate}
\end{lem}

\begin{proof}
By the definition of $\imath$divided powers, if $p(k) \neq p$, then we have
$$
B_{i,p}^{(k)} = \prod_{a \in [-k+1,k-1]_{p}} (B_i-[a]).
$$
This implies the assertion.
\end{proof}

\begin{prop}\label{canonical Xi-weight structure is standard}
Let $M$ be a weight $\U$-module. Then, its canonical $X^\imath$-weight module structure is standard.
\end{prop}

\begin{proof}
Let $i \in I$ and $l \in \Z_{\geq 0}$, and consider the $(l+1)$-dimensional irreducible $\U_i$-module $V(l)$, where $\U_i$ denotes the subalgebra of $\U$ generated by $E_i,F_i,K_i^{\pm 1}$. Let $v \in V(l)$ be a highest weight vector. By \cite[Theorems 2.10 and 3.6]{BeW18}, $V(l)$ has a basis $\{ B_{i,p(l)}^{(k)} v \mid k \in [0,l] \}$, and, we have $B_{i,p(l)}^{(l+1)} v = 0$. Since $B_{i,p(l)}^{(k)}$'s commute with each other, we obtain
$$
B_{i,p(l)}^{(l+1)} u = 0 \Forall u \in V(l).
$$
By Lemma \ref{characterization of standard Xi weight}, we see that each vector in $V(l)$ is a sum of $B_i$-eigenvectors of eigenvalues of the form $[a]$, $a \in [-l,l]_{p(l)}$.

Let $\lm \in X$, $v \in M_\lm$, and $i \in I$. We show that $B_{i,p(\lm_i)}^{(k_i)} v = 0$ for some $k_i > 0$ such that $\ol{k_i} \neq \ol{\lm_i}$. As a $\U_i$-module, $M$ decomposes as
$$
M = \bigoplus_{l \geq 0} M[l],
$$
where $M[l]$ denotes the isotypic component of $M$ of type $V(l)$. By weight consideration, we have $v \in \bigoplus_{l \in \Z_{p(\lm_i)}} M[l]$. By argument above, $v$ is a sum of $B_i$-eigenvectors of eigenvalues of the form $[a]$, $a \in \Z_{p(\lm_i)}$. By Lemma \ref{characterization of standard Xi weight}, we have $B_{i,p(\lm_i)}^{(k_i)} v = 0$ for some $k_i > 0$ such that $\ol{k_i} \neq \ol{\lm_i}$.

For each $\zeta \in X^\imath$, set
$$
M'_\zeta := \{ v \in M \mid B_{i,p(\zeta_i)}^{(k_i)} v = 0 \Forsome k_i > 0 \text{ such that } \ol{k_i} \neq \zeta_i \Forall i \in I \}.
$$
By above, we have
$$
M_\zeta = \bigoplus_{\ol{\lm} = \zeta} M_\lm \subseteq M'_\zeta \subseteq M.
$$
Summing up through $\zeta \in X^\imath$, we obtain
$$
M = \bigoplus_{\zeta \in X^\imath} M_\zeta = \sum_{\zeta \in X^\imath} M'_\zeta.
$$
Since $M'_\zeta$'s intersects trivially with each other, we conclude that
$$
M_\zeta = M'_\zeta.
$$
This proves the assertion.
\end{proof}

\subsection{Classical weight modules with contragredient Hermitian inner products}
Let $M$ be a finite-dimensional classical weight module equipped with a contragredient Hermitian inner product. The following are easy analogs of Lemma \ref{orthogonality of different weight spaces} and Proposition \ref{orthogonal irreducible decomposition}.

\begin{lem}
Let $\nu \neq \xi \in X_{\frk,\Int}$. Then, we have $(M_\nu,M_\xi) = 0$. Consequently, we have $\clL_M = \bigoplus_{\nu \in X_{\frk,\Int}} \clL_{M,\nu}$ and $\ol{\clL}_M = \bigoplus_{\nu \in X_{\frk,\Int}} \ol{\clL}_{M,\nu}$, where $\clL_{M,\nu} := \clL_M \cap M_\nu$ and $\ol{\clL}_{M,\nu} := \clL_{M,\nu}/q\inv \clL_{M,\nu}$.
\end{lem}

\begin{prop}\label{existence of orthogonal irreducible decomposition for AI}
There exists an orthogonal irreducible decomposition $M = \bigoplus_{k=1}^r M_k$ of $M$.
\end{prop}

Because of Proposition \ref{existence of orthogonal irreducible decomposition for AI}, we can define $M[\nu], M[\geq \nu], M[> \nu]$, etc. in the same way as in the quantum groups setting.

Suppose further that $M$ is a standard $X^\imath$-weight module. Then, its distinct $X^\imath$-weight spaces are orthogonal to each other:

\begin{lem}
For each $\zeta \neq \eta \in X^\imath$, we have $(M_\zeta,M_\eta) = 0$.
\end{lem}

\begin{proof}
Suppose that $\zeta \neq \eta$. Then, there exists $i \in I$ such that $\zeta_i \neq \eta_i$. Let $u \in M_\zeta$ and $v \in M_\eta$. By Lemma \ref{characterization of standard Xi weight}, we can write $u = \sum_k u_k$ (resp., $v = \sum_l v_l$) with $B_i u_k = [a_k] u_k$, $\ol{a_k} = \zeta_i$ (resp., $B_i v_l = [b_l] v_l$, $\ol{b_l} = \eta_i$). Then, for each $k,l$, we have
$$
[a_k](u_k,v_l) = (B_i u_k,v_l) = (u_k, B_iv_l) = [b_l](u_k,v_l).
$$
Since $a_k \neq b_l$, it follows that $(u_k,v_l) = 0$. Hence, we obtain $(u,v) = \sum_{k,l}(u_k,v_l) = 0$, as required.
\end{proof}

From this lemma together with Proposition \ref{orthogonal decomposition and local space}, we obtain orthogonal decompositions
$$
\clL_M = \bigoplus_{\zeta \in X^\imath} \clL_{M,\zeta}, \qu \ol{\clL}_M = \bigoplus_{\zeta \in X^\imath} \ol{\clL}_{M,\zeta},
$$
where $\clL_{M,\zeta} := \clL_M \cap M_\zeta$ and $\ol{\clL}_{M,\zeta} := \clL_{M,\zeta}/q\inv \clL_{M,\zeta}$. Hence, for each $\zeta \in X^\imath$, the idempotent $\mathbf{1}_\zeta$ defines a projection
$$
\mathbf{1}_\zeta : \ol{\clL}_M \rightarrow \ol{\clL}_{M,\zeta}.
$$

\section{$n = 2$ case}\label{section: n=2}
Sections \ref{section: n=2}--\ref{section: n=4} are devoted to investigating low rank cases; $n = 2,3,4$. The results obtained in these sections will be used to prove our main theorems.

In this section, we consider the $n=2$ case. In this case, we can identify $X = \Z$, $X^+ = \Z_{\geq 0}$, $X^\imath = \Z/2\Z$. Set $X_{\frk,\Int} := \Z$. For each $\nu \in X_{\frk,\Int}$, let $\bfK_\nu = \bfK v_\nu$ denote the $1$-dimensional $\Ui$-module such that $B_1 v_\nu = [\nu] v_\nu$. Clearly, $\bfK_\nu$ possesses a unique contragredient Hermitian inner product $(\cdot,\cdot)_\nu$ such that $(v_\nu,v_\nu)_\nu = 1$. Set $\clL(\nu) := \bfK_\infty v_\nu$, $\ol{\clL}(\nu) := \clL(\nu)/q\inv \clL(\nu)$, and $b_\nu := \ev_\infty(v_\nu)$.

In the $n = 2$ case, we define a classical weight module to be a $\Ui$-module isomorphic to $\bigoplus_{\nu \in X_{\frk,\Int}} \bfK_\nu^{\oplus m_\nu}$ for some $m_\nu \geq 0$. Then, each classical weight module $M$ admits an $X_{\frk}$-weight space decomposition $M = \bigoplus_{\nu \in X_{\frk,\Int}} M_\nu$ as in the $n \geq 3$ case.

\subsection{Modified action of $B_1$}\label{subsec: Btilde}
Let $M = \bigoplus_{\nu \in X_{\frk,\Int}} M_\nu$ be a finite-dimensional classical weight module equipped with a contragredient Hermitian inner product.

Let $\Btil_1$ be a linear operator on $M$ defined by
$$
\Btil_1 v = \begin{cases}
v \qu & \IF v \in M_\nu \Forsome \nu > 0, \\
0 \qu & \IF v \in M_0, \\
-v \qu & \IF v \in M_\nu \Forsome \nu < 0.
\end{cases}
$$
Then, $\Btil_1$ preserves $\clL_M$. Hence, it induces a $\C$-linear operator $\Btil_1$ on $\ol{\clL}_M$.

Let $\nu \in X_{\frk,\Int}$ and consider the irreducible $\Ui$-module $\bfK_\nu = \bfK v_\nu$. Let $V_\natural = \bfK u_1 \oplus \bfK u_2$ denote the natural representation of $\U = U_q(\frsl_2)$. Recall from Proposition \ref{classical weight tensor fd is classical weight} that $\bfK_\nu \otimes V_\natural$ is a classical weight module. Set
\begin{align}\label{definition of vnupm}
\begin{split}
&v'_{\nu,\pm} := v_\nu \otimes (u_1 \pm q^{\pm \nu} u_2), \\
&v_{\nu,\pm} := \frac{1}{\lt(v'_{\nu,\pm})} v'_{\nu,\pm} = \begin{cases} v_\nu \otimes (q^{\mp \nu} u_1 \pm u_2) \qu & \IF \pm \nu > 0, \\
\frac{1}{\sqrt{2}} v_\nu \otimes (u_1 \pm u_2) \qu & \IF \nu = 0, \\
v_\nu \otimes (u_1 \pm q^{\pm \nu} u_2) \qu & \IF \pm \nu < 0.
\end{cases}
\end{split}
\end{align}
By Lemma \ref{X_k-weight vectors in M tensor Vnatural}, we have
\begin{align}\label{B_1-eigenvalue of v_pm}
B_1v_{\nu,\pm} = [\nu \pm 1] v_{\nu,\pm}.
\end{align}
Since $\dim \bfK_\nu \otimes V_\natural = 2$, we obtain the irreducible decomposition
$$
\bfK_\nu \otimes V_\natural = \bfK v_{\nu,+} \oplus \bfK v_{\nu,-} \simeq \bfK_{\nu+1} \oplus \bfK_{\nu-1}.
$$
By the definition \eqref{definition of vnupm} of $v_{\nu,\pm}$, we see that
\begin{align}\label{v_pm at infty}
\ev_\infty(v_{\nu,\pm}) = \begin{cases}
\pm b_\nu \otimes \ol{u}_2 \qu & \IF \pm \nu > 0, \\
\frac{1}{\sqrt{2}} b_\nu \otimes (\ol{u}_1 \pm \ol{u}_2) \qu & \IF \nu = 0, \\
b_\nu \otimes \ol{u}_1 \qu & \IF \pm \nu < 0.
\end{cases}
\end{align}

\begin{prop}
Let $\nu \in X_{\frk,\Int}$. Then, we have
\begin{align}
\begin{split}
\Btil_1(b_\nu \otimes \ol{u}_1) &= \begin{cases}
\pm b_\nu \otimes \ol{u}_1 \qu & \IF \pm \nu > 1, \\
0 \qu & \IF \nu = \pm 1, \\
b_\nu \otimes \ol{u}_2 \qu & \IF \nu = 0,
\end{cases} \\
&= \begin{cases}
\Btil_1 b_\nu \otimes \ol{u}_1 \qu & \IF |\nu| > 1, \\
0 \qu & \IF |\nu| = 1, \\
b_\nu \otimes \ol{u}_2 \qu & \IF \nu = 0,
\end{cases} \\
\Btil_1(b_\nu \otimes \ol{u}_2) &= \begin{cases}
\pm b_\nu \otimes \ol{u}_2 \qu & \IF \pm \nu >0, \\
b_\nu \otimes \ol{u}_1 \qu & \IF \nu = 0,
\end{cases} \\
&= \begin{cases}
\Btil_1 b_\nu \otimes \ol{u}_2 \qu & \IF |\nu| > 0, \\
b_\nu \otimes \ol{u}_1 \qu & \IF \nu = 0.
\end{cases}
\end{split} \nonumber
\end{align}
\end{prop}

\begin{proof}
The assertion is a consequence of formulas \eqref{B_1-eigenvalue of v_pm} and \eqref{v_pm at infty}.
\end{proof}

Let $M$ be a finite-dimensional classical weight module equipped with a contragredient Hermitian inner product.

\begin{defi}\normalfont
  A vector $b \in \ol{\clL}_M$ is said to be $B_1$-homogeneous of degree $\nu$ if it belongs to $\ol{\clL}_{M,\nu} + \ol{\clL}_{M,-\nu}$ for some $\nu \in X_{\frk,\Int}$ such that $\nu \geq 0$. If $b \in \ol{\clL}_M$ is $B_1$-homogeneous of degree $\nu$, we write $\deg_1(b) = \nu$.
\end{defi}

The following lemmas are clear from the definitions and previous results.  

\begin{lem}\label{weight decomposition of homogeneous vectors}
Let $b \in \ol{\clL}_M \setminus \{0\}$ be $B_1$-homogeneous. Then, we have exactly one of either $\Btil_1 b = 0$ or $\Btil_1^2 b = b$. Furthermore, $(1 \pm \Btil_1) b$ is an $X_{\frk}$-weight vector of weight $\pm \deg_1(b)$.
\end{lem}

\begin{lem}\label{B1 and deg on M tensor Vnatural}
Let $b \in \ol{\clL}_M$ be $B_1$-homogeneous. Then, we have
\begin{align}
\begin{split}
&\deg_1(b \otimes \ol{u}_1) = \begin{cases}
\deg_1(b) - 1 \qu & \IF \deg_1(b) > 0, \\
1 \qu & \IF \deg_1(b) = 0,
\end{cases} \\
&\deg_1(b \otimes \ol{u}_2) = \deg_1(b) + 1, \\
&\Btil_1(b \otimes \ol{u}_1) = \begin{cases}
\Btil_1 b \otimes \ol{u}_1 \qu & \IF \deg_1(b) > 1, \\
0 \qu & \IF \deg_1(b) = 1, \\
b \otimes \ol{u}_2 \qu & \IF \deg_1(b) = 0,
\end{cases} \\
&\Btil_1(b \otimes \ol{u}_2) = \begin{cases}
\Btil_1 b \otimes \ol{u}_2 \qu & \IF \deg_1(b) > 0, \\
b \otimes \ol{u}_1 \qu & \IF \deg_1(b) = 0.
\end{cases}
\end{split} \nonumber
\end{align}
\end{lem}

We want to describe the action of $\Btil_1$ on $\ol{\clL}_M \otimes \ol{\clL}_\natural^{\otimes d}$ for an arbitrary $d \geq 0$. The space $\ol{\clL}_\natural^{\otimes d}$ possesses a crystal basis $\clB_\natural^{\otimes d} := \{ \ol{u}_{i_1} \otimes \cdots \otimes \ol{u}_{i_d} \mid i_1,\ldots,i_d \in \{ 1,2 \} \}$. For the reader's convenience, we write down the crystal structure of $\clB_\natural^{\otimes d} \otimes \clB_\natural$, which can be derived from equations \eqref{tensor product rule for usual crystal} in page \pageref{tensor product rule for usual crystal}: For each $\ol{u} = \ol{u}_{i_1} \otimes \cdots \otimes \ol{u}_{i_d}$, $i_1,\ldots,i_d \in \{ 1,2 \}$, we have
\begin{align}\label{crystal structure of Vnatural tensor d+1}
\begin{split}
&\Ftil_1(\ol{u} \otimes \ol{u}_1) = \begin{cases}
\ol{u} \otimes \ol{u}_2 \qu & \IF \vep_1(\ol{u}) < 1, \\
\Ftil_1 \ol{u} \otimes \ol{u}_1 \qu & \IF \vep_1(\ol{u}) \geq 1,
\end{cases} \\
&\Etil_1(\ol{u} \otimes \ol{u}_1) = \begin{cases}
0 \qu & \IF \vep_1(\ol{u}) \leq 1, \\
\Etil_1 \ol{u} \otimes \ol{u}_1 \qu & \IF \vep_1(\ol{u}) > 1,
\end{cases} \\
&\vphi_1(\ol{u} \otimes \ol{u}_1) = \begin{cases}
\vphi_1(\ol{u}) + 1 \qu & \IF \vep_1(\ol{u}) < 1, \\
\vphi_1(\ol{u}) \qu & \IF \vep_1(\ol{u}) \geq 1,
\end{cases} \\
&\vep_1(\ol{u} \otimes \ol{u}_1) = \begin{cases}
0 \qu & \IF \vep_1(\ol{u}) \leq 1, \\
\vep_1(\ol{u})-1 \qu & \IF \vep_1(\ol{u}) > 1,
\end{cases} \\
&\Ftil_1(\ol{u} \otimes \ol{u}_2) = \Ftil_1 \ol{u} \otimes \ol{u}_2, \\
&\Etil_1(\ol{u} \otimes \ol{u}_2) = \begin{cases}
\ol{u} \otimes \ol{u}_1 \qu & \IF \vep_1(\ol{u}) = 0, \\
\Etil_1 \ol{u} \otimes \ol{u}_2 \qu & \IF \vep_1(\ol{u}) > 0,
\end{cases} \\
&\vphi_1(\ol{u} \otimes \ol{u}_2) =\vphi_1(\ol{u}), \\
&\vep_1(\ol{u} \otimes \ol{u}_2) = \vep_1(\ol{u})+1.
\end{split}
\end{align}

\begin{prop}\label{tensor product rule for deg and Btil}
Let $M$ be a finite-dimensional classical weight module equipped with a contragredient Hermitian inner product. Let $b \in \ol{\clL}_M$ be $B_1$-homogeneous and $\ol{u} = \ol{u}_{i_1} \otimes \cdots \ol{u}_{i_d} \in \clB_\natural^{\otimes d}$. Then, $b \otimes \ol{u}$ is $B_1$-homogeneous. Moreover, we have
\begin{align}
\begin{split}
&\deg_1(b \otimes \ol{u}) = \begin{cases}
\deg_1(b) - \vphi_1(\ol{u}) + \vep_1(\ol{u}) \qu & \IF \deg_1(b) > \vphi_1(\ol{u}), \\
\vep_1(\ol{u}) \qu & \IF \vphi_1(\ol{u}) - \deg_1(b) \in \Z_{\geq 0, \ev}, \\
\vep_1(\ol{u}) + 1 \qu & \IF \vphi_1(\ol{u}) - \deg_1(b) \in \Z_{\geq 0, \odd},
\end{cases} \\
&\Btil_1(b \otimes \ol{u}) = \begin{cases}
\Btil_1 b \otimes \ol{u} \qu & \IF \deg_1(b) > \vphi_1(\ol{u}), \\
b \otimes \Etil_1 \ol{u} \qu & \IF \vphi_1(\ol{u}) - \deg_1(b) \in \Z_{\geq 0, \ev}, \\
b \otimes \Ftil_1 \ol{u} \qu & \IF \vphi_1(\ol{u}) - \deg_1(b) \in \Z_{\geq 0, \odd}.
\end{cases}
\end{split} \nonumber
\end{align}
\end{prop}

\begin{proof}
We proceed by induction on $d$; the $d = 0$ case is clear. Let $d \geq 0$ and assume that the proposition is true, and let us prove the $d+1$ case. During the proof below, we often use Lemma \ref{B1 and deg on M tensor Vnatural} and equations \eqref{crystal structure of Vnatural tensor d+1} without mentioning one by one. First, we compute $\deg_1$ and $\Btil_1$ for $b \otimes \ol{u} \otimes \ol{u}_1$.

Suppose that $\vep_1(\ol{u}) > 0$. Then, we have $\vphi_1(\ol{u} \otimes \ol{u}_1) = \vphi_1(\ol{u})$ and $\vep_1(\ol{u} \otimes \ol{u}_1) = \vep(\ol{u})-1$.
\begin{itemize}
\item When $\deg_1(b) > \vphi_1(\ol{u} \otimes \ol{u}_1)$. In this case, we have
$$
\deg_1(b \otimes \ol{u}) = \deg_1(b) - \vphi_1(\ol{u}) + \vep_1(\ol{u}) \geq \vep_1(\ol{u})+1 \geq 2.
$$
Hence, we obtain
\begin{align}
\begin{split}
\deg_1(b \otimes \ol{u} \otimes \ol{u}_1) &= \deg_1(b \otimes \ol{u}_1)-1 \\
&= \deg_1(b) - \vphi_1(\ol{u}) + \vep_1(\ol{u})-1 \\
& = \deg_1(b) - \vphi_1(\ol{u} \otimes \ol{u}_1) + \vep_1(\ol{u} \otimes \ol{u}_1),
\end{split} \nonumber
\end{align}
and
\begin{align}
\begin{split}
\Btil_1(b \otimes \ol{u} \otimes \ol{u}_1) &= \Btil_1(b \otimes \ol{u}) \otimes \ol{u}_1 \\
&= \Btil_1 b \otimes \ol{u} \otimes \ol{u}_1.
\end{split} \nonumber
\end{align}
\item When $\vphi_1(\ol{u} \otimes \ol{u}_1) - \deg_1(b) \in \Z_{\geq 0, \ev}$. In this case, we have
$$
\deg_1(b \otimes \ol{u}) = \vep_1(\ol{u}) \geq 1.
$$
Hence, we obtain
\begin{align}
\begin{split}
\deg_1(b \otimes \ol{u} \otimes \ol{u}_1) &= \deg_1(b \otimes \ol{u}_1)-1 \\
&= \vep_1(\ol{u})-1 \\
&= \vep_1(\ol{u} \otimes \ol{u}_1),
\end{split} \nonumber
\end{align}
and
\begin{align}
\begin{split}
\Btil_1(b \otimes \ol{u} \otimes \ol{u}_1) &= \begin{cases}
\Btil_1(b \otimes \ol{u}) \otimes \ol{u}_1 \qu & \IF \vep_1(\ol{u}) > 1, \\
0 \qu & \IF \vep_1(\ol{u}) = 1
\end{cases} \\
&= \begin{cases}
b \otimes \Etil_1 \ol{u} \otimes \ol{u}_1 \qu & \IF \vep_1(\ol{u}) > 1, \\
0 \qu & \IF \vep_1(\ol{u}) = 1
\end{cases} \\
&= b \otimes \Etil_1(\ol{u} \otimes \ol{u}_1).
\end{split} \nonumber
\end{align}
\item When $\vphi_1(\ol{u} \otimes \ol{u}_1) - \deg_1(b) \in \Z_{\geq 0, \odd}$. In this case, we have
$$
\deg_1(b \otimes \ol{u}) = \vep_1(\ol{u}) + 1 \geq 2.
$$
Hence, we obtain
\begin{align}
\begin{split}
\deg_1(b \otimes \ol{u} \otimes \ol{u}_1) &= \deg_1(b \otimes \ol{u}_1)-1 \\
&= \vep_1(\ol{u}) \\
&= \vep_1(\ol{u} \otimes \ol{u}_1) + 1,
\end{split} \nonumber
\end{align}
and
\begin{align}
\begin{split}
\Btil_1(b \otimes \ol{u} \otimes \ol{u}_1) &= \Btil_1(b \otimes \ol{u}) \otimes \ol{u}_1 \\
&= b \otimes \Ftil_1 \ol{u} \otimes \ol{u}_1 \\
&= b \otimes \Ftil_1(\ol{u} \otimes \ol{u}_1).
\end{split} \nonumber
\end{align}
\end{itemize}
This case-by-case analysis proves the assertion for the $\vep_1(\ol{u}) > 0$ case.

Next, suppose that $\vep_1(\ol{u}) = 0$. Then, we have $\vphi_1(\ol{u} \otimes \ol{u}_1) = \vphi_1(\ol{u}) + 1$ and $\vep_1(\ol{u} \otimes \ol{u}_1) = 0$.
\begin{itemize}
\item When $\deg_1(b) > \vphi_1(\ol{u} \otimes \ol{u}_1)$. In this case, we have
$$
\deg_1(b \otimes \ol{u}) = \deg_1(b) - \vphi_1(\ol{u}) > 1.
$$
Therefore, we obtain
\begin{align}
\begin{split}
\deg_1(b \otimes \ol{u} \otimes \ol{u}_1) &= \deg_1(b \otimes \ol{u}) - 1 \\
&=\deg_1(b) - \vphi_1(\ol{u}) - 1 \\
&=\deg_1(b) - \vphi_1(\ol{u} \otimes \ol{u}_1) + \vep_1(\ol{u} \otimes \ol{u}_1),
\end{split} \nonumber
\end{align}
and
\begin{align}
\begin{split}
\Btil_1(b \otimes \ol{u} \otimes \ol{u}_1) &= \Btil_1 (b \otimes \ol{u}) \otimes \ol{u}_1 \\
&= \Btil_1 b \otimes \ol{u} \otimes \ol{u}_1.
\end{split} \nonumber
\end{align}
\item When $\vphi_1(\ol{u} \otimes \ol{u}_1) - \deg_1(\ol{u}) \in \Z_{\geq 0, \ev}$. If $\deg_1(b) > \vphi_1(\ol{u})$, then we must have $\deg_1(b) = \vphi_1(\ol{u})+1$. Hence, we have
$$
\deg_1(b \otimes \ol{u}) = \deg_1(b) - \vphi_1(\ol{u}) = 1.
$$
Therefore, we obtain
$$
\deg_1(b \otimes \ol{u} \otimes \ol{u}_1) = \deg_1(b \otimes \ol{u}_1)-1 = 0 = \vep_1(\ol{u} \otimes \ol{u}_1),
$$
and
$$
\Btil_1(b \otimes \ol{u} \otimes \ol{u}_1) = 0 = b \otimes \Etil_1(\ol{u} \otimes \ol{u}_1).
$$
On the other hand, if $\deg_1(b) \leq \vphi_1(\ol{u})$, then we have $\vphi_1(\ol{u}) - \deg_1(b) \in \Z_{\geq 0, \odd}$, and hence
$$
\deg_1(b \otimes \ol{u}) = \vep_1(\ol{u}) + 1 = 1.
$$
Therefore, $\deg_1(b \otimes \ol{u} \otimes \ol{u}_1)$ and $\Btil_1(b \otimes \ol{u} \otimes \ol{u}_1)$ are the same as before.
\item When $\vphi_1(\ol{u} \otimes \ol{u}_1) - \deg_1(b) \in \Z_{\geq 0, \odd}$. In this case, we have $\vphi_1(\ol{u}) - \deg_1(b) \in \Z_{\geq 0, \ev}$. Hence, we have
$$
\deg_1(b \otimes \ol{u}) = \vep_1(\ol{u}) = 0.
$$
Therefore, we obtain
$$
\deg_1(b \otimes \ol{u} \otimes \ol{u}_1) = 1 = \vep_1(\ol{u} \otimes \ol{u}_1)+1,
$$
and
$$
\Btil_1(b \otimes \ol{u} \otimes \ol{u}_1) = b \otimes \ol{u} \otimes \ol{u}_2 = b \otimes \Ftil_1(\ol{u} \otimes \ol{u}_1).
$$
\end{itemize}
This case-by-case analysis proves the assertion for the $\vep_1(\ol{u}) = 0$ case.

Now, let us compute $\deg_1$ and $\Btil_1$ for $b \otimes \ol{u} \otimes \ol{u}_2$.
\begin{itemize}
\item When $\deg_1(b) > \vphi_1(\ol{u} \otimes \ol{u}_2)$. In this case, we have
$$
\deg_1(b \otimes \ol{u}) = \deg_1(b) - \vphi_1(\ol{u}) + \vep_1(\ol{u}) \geq \vep_1(\ol{u})+1 \geq 1.
$$
Hence, we obtain
\begin{align}
\begin{split}
\deg_1(b \otimes \ol{u} \otimes \ol{u}_2) &= \deg_1(b \otimes \ol{u}) + 1 \\
&= \deg_1(b) - \vphi_1(\ol{u}) + \vep_1(\ol{u}) + 1 \\
&= \deg_1(b) - \vphi_1(\ol{u} \otimes \ol{u}_2) + \vep_1(\ol{u} \otimes \ol{u}_1),
\end{split} \nonumber
\end{align}
and
\begin{align}
\begin{split}
\Btil_1(b \otimes \ol{u} \otimes \ol{u}_2) &= \Btil_1(b \otimes \ol{u}) \otimes \ol{u}_2 \\
&= \Btil_1 b \otimes \ol{u} \otimes \ol{u}_2.
\end{split} \nonumber
\end{align}
\item When $\vphi_1(\ol{u} \otimes \ol{u}_2) - \deg_1(b) \in \Z_{\geq 0, \ev}$. In this case, we have
$$
\deg_1(b \otimes \ol{u}) = \vep_1(\ol{u}) \geq 0.
$$
Hence, we obtain
\begin{align}
\begin{split}
\deg_1(b \otimes \ol{u} \otimes \ol{u}_2) &= \deg_1(b \otimes \ol{u}) + 1 \\
&= \vep_1(\ol{u}) + 1 \\
&= \vep_1(\ol{u} \otimes \ol{u}_2),
\end{split} \nonumber
\end{align}
and
\begin{align}
\begin{split}
\Btil_1(b \otimes \ol{u} \otimes \ol{u}_2) &= \begin{cases}
\Btil_1(b \otimes \ol{u}) \otimes \ol{u}_2 \qu & \IF \vep_1(\ol{u}) > 0, \\
b \otimes \ol{u} \otimes \ol{u}_1 \qu & \IF \vep_1(\ol{u}) = 0
\end{cases} \\
&= \begin{cases}
b \otimes \Etil_1 \ol{u} \otimes \ol{u}_2 \qu & \IF \vep_1(\ol{u}) > 0, \\
b \otimes \ol{u} \otimes \ol{u}_1 \qu & \IF \vep_1(\ol{u}) = 0
\end{cases} \\
&= b \otimes \Etil_1(\ol{u} \otimes \ol{u}_2).
\end{split} \nonumber
\end{align}
\item When $\vphi_1(\ol{u} \otimes \ol{u}_2) - \deg_1(b) \in \Z_{\geq 0, \odd}$. In this case, we have
$$
\deg_1(b \otimes \ol{u}) = \vep_1(\ol{u}) + 1 \geq 1.
$$
Hence, we obtain
\begin{align}
\begin{split}
\deg_1(b \otimes \ol{u} \otimes \ol{u}_2) &= \deg_1(b \otimes \ol{u}) + 1 \\
&= \vep_1(\ol{u}) + 2 \\
&= \vep_1(\ol{u} \otimes \ol{u}_2) + 1,
\end{split} \nonumber
\end{align}
and
\begin{align}
\begin{split}
\Btil_1(b \otimes \ol{u} \otimes \ol{u}_2) &= \Btil_1(b \otimes \ol{u}) \otimes \ol{u}_2 \\
&= b \otimes \Ftil_1 \ol{u} \otimes \ol{u}_2 \\
&= b \otimes \Ftil_1(\ol{u} \otimes \ol{u}_2).
\end{split} \nonumber
\end{align}
\end{itemize}
Thus, the proof completes.
\end{proof}

\begin{cor}\label{deg and Btil on usual crystals}
Let $M$ be a finite-dimensional $\U$-module with a crystal basis $\clB_M$. For each $b \in \clB_M$, the following hold:
\begin{enumerate}
\item $b$ is $B_1$-homogeneous;
$$
\deg_1(b) = \begin{cases}
\vep_1(b) \qu & \IF \vphi_1(b) \in \Z_{\ev}, \\
\vep_1(b) + 1 \qu & \IF \vphi_1(b) \in \Z_{\odd}.
\end{cases}
$$
\item $\Btil_1 b \in \clB_M \sqcup \{0\}$;
$$
\Btil_1 b = \begin{cases}
\Etil_1 b \qu & \IF \vphi_1(b) \in \Z_{\ev}, \\
\Ftil_1 b \qu & \IF \vphi_1(b) \in \Z_{\odd}.
\end{cases}
$$
\end{enumerate}
\end{cor}

\begin{proof}
Since $\clB_M$ can be embedded into $\clB_\natural^{\otimes d}$ for some $d \geq 0$, it suffices to prove the assertions for $M = V_\natural^{\otimes d}$. The latter follows from Proposition \ref{tensor product rule for deg and Btil} by considering $\bfK_0 \otimes V_\natural^{\otimes d} \simeq V_\natural^{\otimes d}$.
\end{proof}

\begin{cor}\label{orthonormal weight basis of usual crystal}
Let $M$ be a finite-dimensional $\U$-module with a crystal basis $\clB_M$. Then, the following forms an orthonormal basis of $\ol{\clL}_M$ consisting of $X_{\frk}$-weight vectors:
$$
\{ b \in \clB_M \mid \vphi_1(b) \in \Z_{\ev} \AND \vep_1(b) = 0 \} \sqcup \{ \tfrac{1}{\sqrt{2}}(b \pm \Ftil_1(b)) \mid b \in \clB_M \AND \vphi_1(b) \in \Z_{\odd} \}.
$$
\end{cor}

\begin{proof}
The assertion follows from Lemma \ref{weight decomposition of homogeneous vectors} and Corollary \ref{deg and Btil on usual crystals}.
\end{proof}

\begin{rem}\normalfont
In the general $n \geq 3$ case, we define $\Btil_i$ and $\deg_i$ for all $i \in I$ in the obvious way.
\end{rem}

\subsection{Based $\Ui$-module structures of irreducible $\U$-modules}
Let $\lm \in X^+ = \Z_{\geq 0}$, and consider the irreducible $\U$-module $V(\lm)$. Since the weights of $V(\lm)$ belong to the interval $[-\lm,\lm]_{p(\lm)}$, each vector in $V(\lm)$ is an $X^\imath$-weight vector of weight $\ol{\lm}$.

Recall that the set $\{ F_1^{(k)} v_\lm \mid k \in [0,\lm] \}$ forms the canonical basis $G(\lm)$, and the set $\{ B_{1,p(\lm)}^{(k)} v_\lm \mid k \in [0,\lm] \}$ forms the $\imath$canonical basis $G^\imath(\lm)$ of $V(\lm)$. Viewing $V(\lm)$ as a based $\Ui$-module in the sense of Example \ref{V(lm) with icanonical basis is a based Ui-module}, we have $G^\imath(\Ftil_1^k b_\lm) = B_{1,p(\lm)}^{(k)} v_\lm$.

For each $m \in \Z$, set
$$
V(\lm)\{\geq m\} := \bfK \{ B_{1,p(\lm)}^{(k)} v_\lm \mid k \geq m \}.
$$
Then, we have a filtration
$$
0 = V(\lm)\{\geq \lm+1\} \subseteq V(\lm)\{\geq \lm-1\} \subseteq \cdots \subseteq V(\lm)\{\geq \lm-2\lfloor \tfrac{\lm}{2} \rfloor-1\} = V(\lm)
$$
of $X^\imath$-weight modules, where $\lfloor r \rfloor$ denotes the greatest integer not greater than $r \in \Q$. In fact, by Proposition \ref{recursive formula for idivided powers}, we have
$$
B_1 B_{1,p(\lm)}^{(\lm-2m-1)}v_\lm \equiv [\lm-2m] B_{1,p(\lm)}^{(\lm-2m)}v_\lm, \qu B_1 B_{1,p(\lm)}^{(\lm-2m)}v_\lm \equiv [\lm-2m]B_{1,p(\lm)}^{(\lm-2m-1)}v_\lm
$$
modulo $V(\lm)\{\geq \lm-2m+1\}$.

For each $\nu \in [-\lm,\lm]_{p(\lm)}$, set
$$
b_{\lm, \nu} := \begin{cases}
b_\lm \qu & \IF \nu = 0, \\
\frac{1}{\sqrt{2}} (\Ftil_1^{|\nu|-1} b_\lm \pm \Ftil_1^{|\nu|} b_\lm) \qu & \IF \pm\nu > 0.
\end{cases}
$$
Then, by Corollary \ref{orthonormal weight basis of usual crystal}, the vector $b_{\lm,\nu} \in \ol{\clL}(\lm)$ is an $X_{\frk}$-weight vector of weight $\nu$, and the set $\{ b_{\lm,\nu} \mid \nu \in [-\lm,\lm]_{p(\lm)} \}$ forms an orthonormal basis of $\ol{\clL}(\lm)$.

\begin{prop}\label{weights for icanonical basis vector for V(l)}
Let $\lm \in X^+$, and $\nu \in [-\lm,\lm]_{p(\lm)}$. Then, $G^\imath(b_{\lm,\nu})$ is a sum of $X_{\frk}$-weight vectors of weights in $\{ \nu, \pm (|\nu| + 2), \pm(|\nu| + 4), \ldots, \pm \lm \}$.
\end{prop}

\begin{proof}
For each $\nu \in [-\lm,\lm]_{p(\lm)}$, define $B_{1,\nu} \in \Ui$ by
$$
B_{1,\nu} := \begin{cases}
1 \qu & \IF \nu = 0, \\
\frac{1}{\sqrt{2}}(B_{1,p(\lm)}^{(|\nu|-1)} \pm B_{1,p(\lm)}^{(|\nu|)}) \qu & \IF \pm \nu > 0.
\end{cases}
$$
Then, by above, we have
$$
B_{1,\nu} v_\lm = G^\imath(b_{\lm,\nu}),
$$
and
$$
B_1 B_{1,\nu} v_\lm \equiv [\nu] B_{1,\nu} v_\lm \qu \pmod{V(\lm)\{\geq |\nu|+1\}}.
$$
Then, the assertion follows by descending induction on $|\nu|$.
\end{proof}

\section{$n = 3$ case}\label{section: n=3}
In this section, we consider the $n=3$ case. In this case, we can identify $X = \Z^2$, $X^+ = \Z^2_{\geq 0}$, $X_{\frk,\Int} = \Z$, $X_{\frk,\Int}^+ = \Z_{\geq 0}$, and $X^\imath = (\Z/2\Z)^2$.

\subsection{Lowering and raising operators}
Set
$$
X_2 = B_{2,+}\{ l_{1};0 \}, \qu Y_2 = B_{2,-}\{ l_{1};0 \}, \qu l := l_{1}.
$$
In this section, when there is no ambiguity, we abbreviate $X_2$ and $Y_2$ as $X$ and $Y$, respectively.

By Lemma \ref{wp for AI} and equation \eqref{formula for AI 5}, we have
$$
(Xu,v) = (u, Y \frac{\{ l;-1 \}}{\{ l;0 \}}v), \qu (Yu,v) = (u, X\frac{\{ l;1 \}}{\{ l;0 \}}v).
$$
for all $u,v \in M$ with $M$ being a $\Ui$-module equipped with a contragredient Hermitian inner product.

Let $\nu \in X_{\frk}^+$, and consider the corresponding irreducible $\Ui$-module $V(\nu)$.

\begin{lem}\label{B_+ B_-^n}
Let $k \geq 0$. Then, we have
$$
X Y^{(k)} v_\nu = [2\nu-k+1] Y^{(k-1)} v_\nu,
$$
where $Y^{(k)} := \frac{1}{[k]!} Y^k$.
\end{lem}

\begin{proof}
By equation \eqref{formula for AI 7}, we have
$$
[X,Y] = [l^2;0].
$$
Then, the assertion is verified by induction on $k$.
\end{proof}

From this lemma (with a standard argument), we see that $\{ Y^{(k)}v_\nu \mid k \in [0,2\nu] \}$ forms a basis of $V(\nu)$. Note that $Y^{(k)}v_\nu \in V(\nu)_{\nu-k}$. For each $k \in [0,2\nu]$, we can choose $c_{k,\nu} \in \bfK^\times$ in a way such that $\lt(c_{k,\nu} Y^{(k)} v_\nu) = 1$. For each $k \in \Z$, set
$$
\Ytil^k v_\nu := \begin{cases}
c_{k,\nu} Y^{(k)} v_\nu \qu & \IF k \in [0,2\nu], \\
0 \qu & \OW.
\end{cases}
$$
Then, $\{ \Ytil^k v_\nu \mid k \in [0,2\nu] \}$ forms an almost orthonormal basis of $V(\nu)$.

Using Lemma \ref{B_+ B_-^n}, we see by induction on $k$ that
\begin{align}
(Y^{(k)}v_\nu, Y^{(k)}v_\nu)_\nu = {2\nu \brack k} \frac{\{ \nu \}}{\{ \nu-k \}}. \nonumber
\end{align}
This shows that for each $k \in [0,2\nu]$, we have
\begin{align}\label{deg and lc of Yn vnu}
\begin{split}
&\deg(Y^{(k)}v_\nu) = \hf((2\nu-k)k + \nu - |\nu-k|) = \begin{cases}
\hf(2\nu-k+1)k \qu & \IF k \leq \nu, \\
\hf(2\nu-k)(k+1) \qu & \IF k \geq \nu,
\end{cases} \\
&\lc(Y^{(k)}v_\nu) = \begin{cases}
\frac{1}{\sqrt{2}} \qu & \IF k = \nu > 0, \\
1 \qu & \OW.
\end{cases}
\end{split}
\end{align}
Consequently, we obtain
$$
\lt(Y^{(k)} v_\nu) = \begin{cases}
1 \qu & \IF \nu = 0, \\
q^{\hf(2\nu-k+1)n} \qu & \IF 0 \leq k < \nu, \\
\frac{1}{\sqrt{2}} q^{\hf \nu(\nu+1)} \qu & \IF k = \nu > 0, \\
q^{\hf(2\nu-k)(k+1)} \qu & \IF \nu < k \leq 2\nu.
\end{cases}
$$

\begin{rem}\normalfont
For each $k \in [0,\nu]$, we have
\begin{align}
\begin{split}
&\deg(\prod_{j=0}^{k-1} \{ \nu-j \}) = \sum_{j=0}^{k-1}(\nu-j) = \hf(2\nu-k+1)k, \\
&\lc(\prod_{j=0}^{k-1} \{ \nu-j \}) = 1.
\end{split} \nonumber
\end{align}
Hence, we may choose
\begin{align}\label{Btil2 for AI}
c_{k,\nu} := \begin{cases}
1 \qu & \IF \nu = 0, \\
\prod_{j=0}^{k-1} \frac{1}{\{ \nu-j \}} \qu & \IF 0 \leq k < \nu, \\
\sqrt{2}\prod_{j=0}^{\nu-1} \frac{1}{\{ \nu-j \}} \qu & \IF k = \nu > 0, \\
\prod_{j=0}^{2\nu-k-1} \frac{1}{\{ \nu-j \}} \qu & \IF \nu < k \leq 2\nu.
\end{cases}
\end{align}
In the sequel, $c_{k,\nu}$ always means this value. Note that we have
$$
\psii_\nu(\Ytil^k v_\nu) = \Ytil^k v_\nu
$$
for all $k \in [0,2\nu]$, and
$$
c_{0,\nu} = c_{2\nu,\nu} = 1.
$$
\end{rem}

\begin{defi}\normalfont
The lowering operator $\Ytil$ and the raising operator $\Xtil$ are $\bfK$-linear endomorphisms on $V(\nu)$ defined by
$$
\Ytil(\Ytil^k v_\nu) := \begin{cases}
\Ytil^{k+1} v_\nu \qu & \IF k \in [0,2\nu], \\
0 \qu & \OW,
\end{cases} \qu \Xtil(\Ytil^k v_\nu) := \begin{cases}
\Ytil^{k-1} v_\nu \qu & \IF k \in [0,2\nu], \\
0 \qu & \OW.
\end{cases}
$$
\end{defi}

Clearly, the operators $\Ytil$ and $\Xtil$ preserve $\clL(\nu)$. Hence, they induce linear operators $\Ytil$ and $\Xtil$ on $\ol{\clL}(\nu)$.


Since each finite-dimensional classical weight module is completely reducible, we can extend the definitions of $\Ytil,\Xtil$ on the irreducible modules to the finite-dimensional classical weight modules. The following are counterparts of well-known results about the Kashiwara operators. The proofs are similar, so we omit them.


\begin{prop}\label{Im Yj = Im Ytilj for rank 1 AI}
Let $M$ be a finite-dimensional classical weight module. Let $\xi \in X_{\frk,\Int}$. Then, we have $Y M_\xi = \Ytil M_\xi$.
\end{prop}

\begin{prop}\label{Ftil and Etil are almost adjoint AI}
Let $M$ be a finite-dimensional classical weight module equipped with a contragredient Hermitian inner product. Then, for each $u,v \in \clL_M$, we have
\begin{align}
\begin{split}
&(\Ytil u,v) - (u, \Xtil v) \in q\inv \bfK_\infty.
\end{split} \nonumber
\end{align}
\end{prop}

\subsection{Tensor product rule}
In this subsection, we investigate the behavior of $\Ytil,\Xtil$ on tensor product modules. Let $V_\natural = \bfK u_1 \oplus \bfK u_2 \oplus \bfK u_3$ denote the natural representation of $\U = U_q(\frsl_3)$.

Let $M$ be a finite-dimensional classical weight module equipped with a contragredient Hermitian inner product, and $\nu \in X_{\frk,\Int}$. For $v \in M_\nu$, set
$$
v'_\pm := v \otimes (u_1 \pm q^{\pm \nu} u_2), \qu v_0 := v \otimes u_3,
$$
and
$$
v_\pm := \frac{\lt(v)}{\lt(v'_\pm)} v'_{\pm} = \begin{cases}
q^{\mp \nu} v'_{\pm} \qu & \IF \pm \nu > 0, \\
\frac{1}{\sqrt{2}}v'_\pm \qu & \IF \nu = 0, \\
v'_\pm \qu & \IF \pm \nu < 0.
\end{cases}
$$
Then, we have
\begin{align}\label{deg and lc of v'pm and v'0}
\begin{split}
&\deg(v'_\pm) = \deg(v)+ \max\{\pm \nu,0\},\qu \deg(v_0) = \deg(v), \\
&\lc(v'_\pm) = \begin{cases}
\lc(v) \qu & \IF \nu \neq 0, \\
\sqrt{2}\lc(v) \qu & \IF \nu = 0.
\end{cases}
\end{split}
\end{align}

\begin{lem}
We have $v'_\pm \in (M \otimes V_\natural)_{\nu \pm 1}$ and $v_0 \in (M \otimes V_\natural)_\nu$.
\end{lem}

\begin{proof}
The assertion follows from Lemma \ref{X_k-weight vectors in M tensor Vnatural}.
\end{proof}

For convenience, we write down easy formulas:
\begin{align}
&v \otimes u_1 = \frac{1}{\{\nu\}}(q^{-\nu}v'_+ + q^\nu v'_-), \label{weight decomposition of v otimes u_1}\\
&v \otimes u_2 = \frac{1}{\{\nu\}}(v'_+ - v'_-). \label{weight decomposition of v otimes u_2}
\end{align}

\begin{prop}\label{X,Y on V otimes Vnatural}
We have
\begin{align}
\begin{split}
&X v'_+ = q\inv (X v)'_+, \\
&Y v'_+ = \frac{\{\nu+1\}}{\{\nu\}}(Yv)'_+ + q^\nu\{\nu+1\}v_0 + \frac{q^\nu(q-q\inv)}{\{\nu\}}(Xv)'_-, \\
&X v'_- = \frac{q^{-\nu}(q-q\inv)}{\{\nu\}}(Yv)'_+ - q^{-\nu}\{\nu-1\}v_0 + \frac{\{\nu-1\}}{\{\nu\}}(Xv)'_-, \\
&Y v'_- = q\inv (Y v)'_-, \\
&X v_0 = v'_+ + q (X v)_0, \\
&Y v_0 = q(Y v)_0 - v'_-. \\
\end{split} \nonumber
\end{align}
\end{prop}

\begin{proof}
Since $\Delta(B_2) = B_2 \otimes K_2\inv + 1 \otimes B_2$, we have
\begin{align}
\begin{split}
&B_2 v'_{\pm} \\
&= B_2v \otimes (u_1 \pm q^{\pm \nu - 1} u_2) \pm q^{\pm \nu} v \otimes u_3 \\
&= (B_{2,+} + B_{2,-})v \otimes u_1 \pm q^{\pm \nu - 1} (B_{2,+} + B_{2,-})v \otimes u_2 \pm q^{\pm \nu} v_0 \\
&= \frac{1}{\{ \nu+1 \}}(q^{-\nu-1}(B_{2,+}v)'_+ + q^{\nu+1}(B_{2,+}v)'_-) + \frac{1}{\{ \nu-1 \}}(q^{-\nu+1}(B_{2,-}v)'_+ + q^{\nu-1}(B_{2,-}v)'_-) \\
&\qu \pm q^{\pm \nu - 1}(\frac{1}{\{ \nu+1 \}}((B_{2,+}v)'_+ - (B_{2,+}v)'_-) + \frac{1}{\{ \nu-1 \}}((B_{2,-}v)'_+ - (B_{2,-}v)'_-)) \pm q^{\pm \nu} v_0.
\end{split} \nonumber
\end{align}
For the last equality, we used equations \eqref{weight decomposition of v otimes u_1} and \eqref{weight decomposition of v otimes u_2}. Rearranging this, we obtain
\begin{align}
\begin{split}
&B_2 v'_+ = \frac{q\inv \{\nu\}}{\{ \nu+1 \}} (B_{2,+}v)'_+ + \frac{q^\nu(q-q\inv)}{\{ \nu+1 \}}(B_{2,+}v)'_- + (B_{2,-}v)'_+ + q^\nu v_0, \\
&B_2 v'_- = (B_{2,+}v)'_- + \frac{q^{-\nu}(q-q\inv)}{\{ \nu-1 \}}(B_{2,-}v)'_+ + \frac{q\inv \{\nu\}}{\{ \nu-1 \}} (B_{2,-}v)'_- - q^{-\nu}v_0.
\end{split} \nonumber
\end{align}
By weight consideration, we conclude
\begin{align}
\begin{split}
&B_{2,+}v'_+ = \frac{q\inv \{\nu\}}{\{ \nu+1 \}} (B_{2,+}v)'_+, \\
&B_{2,-}v'_+ = \frac{q^\nu(q-q\inv)}{\{ \nu+1 \}}(B_{2,+}v)'_- + (B_{2,-}v)'_+ + q^\nu v_0, \\
&B_{2,+}v'_- = (B_{2,+}v)'_- + \frac{q^{-\nu}(q-q\inv)}{\{ \nu-1 \}}(B_{2,-}v)'_+  - q^{-\nu}v_0, \\
&B_{2,-}v'_- = \frac{q\inv \{\nu\}}{\{ \nu-1 \}} (B_{2,-}v)'_-.
\end{split} \nonumber
\end{align}
These imply the first four assertions.

Also, we have
\begin{align}
\begin{split}
B_2 v_0 &= v \otimes u_2 + q B_2v \otimes u_3\\
&= \frac{1}{\{ \nu \}}(v'_+ - v'_-) + q(B_{2,+}v)_0 + q(B_{2,-}v)_0.
\end{split} \nonumber
\end{align}
Again, we used equation \eqref{weight decomposition of v otimes u_2} for the last equality. By weight consideration, we obtain
\begin{align}
\begin{split}
&B_{2,+}v_0 = \frac{1}{\{ \nu \}} v'_+ + q(B_{2,+}v)_0, \\
&B_{2,-}v_0 = -\frac{1}{\{ \nu \}} v'_- + q(B_{2,-}v)_0.
\end{split} \nonumber
\end{align}
These imply the remaining assertions. Thus, we complete the proof.
\end{proof}

\begin{lem}\label{Yn on v'pm and v'0}
Let $M$ be a finite-dimensional classical weight module, $\nu \in X_{\frk,\Int}$, and $v \in M_\nu$. Then, for each $k \geq 0$, we have
\begin{align}
\begin{split}
Y^{(k)} v'_+ &= \frac{\{\nu+1\}}{\{\nu-k+1\}}(Y^{(k)}v)'_+ + q^\nu\{\nu+1\} (Y^{(k-1)}v)_0 \\
&\qu - \frac{\{\nu+1\}}{\{\nu-k+1\}}(Y^{(k-2)}v)'_- + \frac{q^{\nu-k+1}(q-q\inv)}{\{\nu-k+1\}} (Y^{(k-1)}Xv)'_-, \\
Y^{(k)} v_0 &= q^k(Y^{(k)}v)_0 - (Y^{(k-1)}v)'_-, \\
Y^{(k)} v'_- &= q^{-k}(Y^{(k)}v)'_-.
\end{split} \nonumber
\end{align}
\end{lem}

\begin{proof}
Using Proposition \ref{X,Y on V otimes Vnatural}, the assertion is verified by induction on $k$.
\end{proof}

\begin{prop}\label{Yn on vnu+}
Let $\nu \in X_{\frk,\Int}^+$. Then, for each $k \in [0,2\nu+2]$, we have
\begin{align}
\begin{split}
&X(v_\nu)'_+ = 0, \\
&Y^{(k)} (v_\nu)'_+ = \frac{\{ \nu+1 \}}{\{ \nu-k+1 \}} (Y^{(k)}v_\nu)'_+ + q^{\nu} \{ \nu+1 \} (Y^{(k-1)} v_\nu)_0 - \frac{\{ \nu+1 \}}{\{ \nu-k+1 \}} (Y^{(k-2)} v_\nu)'_-, \\
&\ev_\infty(\Ytil^k(v_\nu)_+) = \begin{cases}
\ev_\infty((v_\nu)_+) \qu & \IF k = 0, \\
\ev_\infty((\Ytil^{k-1}v_\nu)_0) \qu & \IF 1 \leq k \leq 2\nu+1, \\
-\ev_\infty((\Ytil^{2\nu}v_\nu)_-) \qu & \IF k = 2\nu+2.
\end{cases}
\end{split} \nonumber
\end{align}
\end{prop}

\begin{proof}
The first and second assertions are verified by using Proposition \ref{X,Y on V otimes Vnatural} and Lemma \ref{Yn on v'pm and v'0}. Let us prove the third assertion. By the second equality and the definition of $\Ytil$, we can write
\begin{align}\label{eq 1}
d_0 \Ytil^k(v_\nu)_+ = d_1 (\Ytil^k v_\nu)_+ + d_2(\Ytil^{k-1}v_\nu)_0 - d_3(\Ytil^{k-2}v_\nu)_-,
\end{align}
where $d_0,d_1,d_2,d_3 \in \bfK$ such that $\lc(d_j) \in \R_{\geq 0}$ and
\begin{align}
\begin{split}
&\deg(d_0) = \deg(Y^{(k)} (v_\nu)'_+),\qu \deg(d_1) = \deg(\tfrac{\{ \nu+1 \}}{\{ \nu-k+1 \}} (Y^{(k)}v_\nu)'_+), \\
&\deg(d_2) = \deg(q^{\nu} \{ \nu+1 \} (Y^{(k-1)} v_\nu)_0),\qu \deg(d_3) = \deg(\tfrac{\{ \nu+1 \}}{\{ \nu-k+1 \}} (Y^{(k-2)} v_\nu)'_-).
\end{split} \nonumber
\end{align}
More explicitly, by equations \eqref{deg and lc of Yn vnu} and \eqref{deg and lc of v'pm and v'0}, we have
\begin{align}
\begin{split}
\deg(d_1) &= \nu+1 - |\nu-k+1| + \hf(n(2\nu-k)+\nu-|\nu-k|) + \max\{\nu-k,0\} \\
&= \nu+1 + \hf(k+1)(2\nu-k) - |\nu-k+1|, \\
\deg(d_2) &= 2\nu+1 + \hf((k-1)(2\nu-k+1)+\nu-|\nu-k+1|), \\
\deg(d_3) &= \nu+1 - |\nu-k+1| + \hf((k-2)(2\nu-k+2)+\nu-|\nu-k+2|) \\
&\qu + \max\{-\nu+k-2,0\} \\
&= \nu+1 + \hf(k-2)(2\nu-k+3) - |\nu-k+1|.
\end{split} \nonumber
\end{align}
Here, we used an easy formula
$$
\max\{x,0\} = \hf(x + |x|).
$$
Then, we compute as
\begin{align}
\begin{split}
\deg(d_2)-\deg(d_1) &= \nu+\hf(-3\nu+3k-1 + |\nu-k+1|) \\
&= \begin{cases}
k \qu & \IF 0 \leq k \leq \nu+1, \\
-\nu+2k-1 \qu & \IF \nu+1 \leq k \leq 2\nu+2,
\end{cases} \\
\deg(d_2)-\deg(d_3) &= \nu+\hf(3\nu-3k+5 + |\nu-k+1|) \\
&= \begin{cases}
3\nu-2k+3 \qu & \IF k \leq \nu+1, \\
2\nu-k+2 \qu & \IF \nu+1 \leq k \leq 2\nu+2.
\end{cases}
\end{split} \nonumber
\end{align}
From these, we obtain
\begin{align}\label{Comparison of degrees}
\begin{cases}
\deg(d_3) \underset{3(\nu-k+1)}{\leq} \deg(d_1) \underset{k}{\leq} \deg(d_2) \qu & \IF 0 \leq k \leq \nu+1, \\
\deg(d_1) \underset{-3(\nu-k+1)}{\leq} \deg(d_3) \underset{2\nu-k+2}{\leq} \deg(d_2) \qu & \IF \nu+1 \leq k \leq 2\nu+2.
\end{cases}
\end{align}
Here, $a \underset{c}{\leq} b$ means $c = b-a \geq 0$.

On the other hand, since the terms in the right-hand side of \eqref{eq 1} are orthogonal to each other, we have
$$
\deg(d_0) = \begin{cases}
\deg(d_1) \qu & \IF k = 0, \\
\max\{\deg(d_1),\deg(d_2)\} \qu & \IF k = 1, \\
\max\{\deg(d_1),\deg(d_2),\deg(d_3)\} \qu & \IF 2 \leq k \leq 2\nu, \\
\max\{\deg(d_2),\deg(d_3)\} \qu & \IF k = 2\nu+1, \\
\deg(d_3) \qu & \IF k = 2\nu+2.
\end{cases}
$$
Combining this and inequalities \eqref{Comparison of degrees}, we obtain
$$
\ev_\infty(\Ytil^n(v_\nu)_+) = \begin{cases}
\frac{\lc(d_1)}{\lc(d_0)} \ev_\infty((v_\nu)_+) \qu & \IF k = 0, \\
\frac{\lc(d_2)}{\lc(d_0)} \ev_\infty((\Ytil^{k-1}v_\nu)_0) \qu & \IF 1 \leq k \leq 2\nu+1, \\
-\frac{\lc(d_3)}{\lc(d_0)} \ev_\infty((\Ytil^{2\nu}v_\nu)_-) \qu & \IF k = 2\nu+2.
\end{cases}
$$
Noting that $\lc(d_j) \in \R_{\geq 0}$ for all $j = 0,1,2,3$, and that
$$
\lc(\Ytil^k (v_\nu)_+) = \lc((v_\nu)_+) = \lc((\Ytil^{k-1} v_\nu)_0) = \lc((\Ytil^{2\nu} v_\nu)_-) = 1,
$$
we finally obtain the required equation.
\end{proof}

The following two results can be proved in a similar way to Proposition \ref{Yn on vnu+}.

\begin{prop}
Let $\nu \in X_{\frk,\Int}^+$ be such that $\nu > 0$. Set
\begin{align}
\begin{split}
&v' := (Y v_\nu)'_+ - q\inv[2\nu] (v_\nu)_0, \\
&v := \frac{1}{\lt(v')}v'.
\end{split} \nonumber
\end{align}
Then, for each $k \in [0,2\nu]$, we have
\begin{align}
\begin{split}
&Xv' = 0, \\
&Y^{(k)}v' = \frac{\{\nu\}}{\{ \nu-k \}}([k+1](Y^{(k+1)} v_\nu)'_+ - q\inv[2\nu-2k] (Y^{(k)} v_\nu)_0 + [2\nu-k+1](Y^{(k-1)} v_\nu)'_-), \\
&\ev_\infty(\Ytil^k v) = \begin{cases}
\ev_\infty((\Ytil^{k+1} v_\nu)_+) \qu & \IF 0 \leq k < \nu, \\
\frac{1}{\sqrt{2}}\ev_\infty((\Ytil^{\nu+1} v_\nu)_+ + (\Ytil^{\nu-1} v_\nu)_-) \qu & \IF k = \nu, \\
\ev_\infty((\Ytil^{k-1} v_\nu)_-) \qu & \IF \nu < k \leq 2\nu.
\end{cases}
\end{split} \nonumber
\end{align}
\end{prop}


\begin{prop}
Let $\nu \in X_{\frk,\Int}^+$ be such that $\nu > 0$. Set
\begin{align}
\begin{split}
&v' := -[2](Y^{(2)}v_\nu)'_+ + q^{-\nu-1}[2\nu-1]\{\nu-1\}(Yv_\nu)_0 + [2\nu][2\nu-1](v_\nu)'_-, \\
&v := \frac{1}{\lt(v')} v'.
\end{split} \nonumber
\end{align}
Then, for each $k \in [0,2\nu-2]$, we have{\small
\begin{align}
\begin{split}
&Xv' = 0, \\
&Y^{(k)} v' = -\frac{[k+2][k+1]\{\nu-1\}}{\{\nu-k-1\}} (Y^{(k+2)} v_\nu)'_+ \\
&\qu + q^{-\nu-1}[k+1][2\nu-k-1]\{ \nu-1 \}(Y^{(k+1)}v_\nu)_0 +\frac{[2\nu-k][2\nu-k-1]\{\nu-1\}}{\{\nu-k-1\}}(Y^{(k)} v_\nu)'_-), \\
&\ev_\infty(\Ytil^k v) = \begin{cases}
\ev_\infty((\Ytil^{k} v_\nu)_-) \qu & \IF 0 \leq k < \nu-1, \\
-\frac{1}{\sqrt{2}}\ev_\infty((\Ytil^{\nu+1} v_\nu)_+ - (\Ytil^{\nu-1} v_\nu)_-) \qu & \IF k = \nu-1, \\
-\ev_\infty((\Ytil^{k+2} v_\nu)_+) \qu & \IF \nu-1 < k \leq 2\nu-2.
\end{cases}
\end{split} \nonumber
\end{align}}
\end{prop}


\subsection{Based module structures of irreducible $\Ui$-modules}
In this subsection, we study $V(\nu)$, $\nu \in X_{\frk,\Int}^+$, and show that it admits a based $\Ui$-module structure.

\begin{prop}\label{1zeta on Ui2-highest weight vector}
Let $\nu \in X_{\frk,\Int}^+$. Then, $V(\nu)$ is a standard $X^\imath$-weight module. Furthermore, the highest weight vector $v_\nu$ has the following $X^\imath$-weight vector decomposition
\begin{align}
\begin{split}
v_\nu &= \hf(v_\nu + (-1)^\nu Y^{(2\nu)}v_\nu) + \hf(v_\nu - (-1)^\nu Y^{(2\nu)}v_\nu) \in V(\nu)_{(\ol{\nu},\ol{0})} \oplus  V(\nu)_{(\ol{\nu},\ol{1})} \\
&= \hf(v_\nu + Y^{(2\nu)}v_\nu) + \hf(v_\nu - Y^{(2\nu)}v_\nu) \in V(\nu)_{(\ol{\nu},\ol{\nu})} \oplus  V(\nu)_{(\ol{\nu},\ol{\nu-1})}
\end{split} \nonumber
\end{align}
\end{prop}

\begin{proof}
Let $M := V_\natural^{\otimes \nu}$. Then, its canonical $X^\imath$-weight module structure is standard by Proposition \ref{canonical Xi-weight structure is standard}. Define $u^{\pm \nu} \in M$ inductively by
$$
u^{\pm 1} := u_1 \pm u_2, \qu u^{\pm k} := (u^{\pm (k-1)})_{\pm} = u^{\pm (k-1)} \otimes (q^{-k+1}u_1 \pm u_2).
$$
Explicitly, we have
$$
u^{\pm \nu} = \sum_{i_1,\ldots,i_\nu \in \{ 1,2 \}} (\pm 1)^{\sharp \{ k \mid i_k = 2 \}} \prod_{\substack{k \\ i_k = 1}} q^{-k+1} u_{i_1} \otimes \cdots \otimes u_{i_\nu}.
$$
Hence, we obtain
$$
u^\nu + u^{-\nu} = 2\sum_{\substack{i_1,\ldots,i_\nu \in \{ 1,2 \} \\ \sharp \{ k \mid i_k = 2 \} \in \Z_{\ev}}} \prod_{\substack{k \\ i_k = 1}} q^{-k+1} u_{i_1} \otimes \cdots \otimes u_{i_\nu} \in M_{(\ol{\nu},\ol{0})},
$$
and
$$
u^\nu - u^{-\nu} = 2\sum_{\substack{i_1,\ldots,i_\nu \in \{ 1,2 \} \\ \sharp \{ k \mid i_k = 2 \} \in \Z_{\odd}}} \prod_{\substack{k \\ i_k = 1}} q^{-k+1} u_{i_1} \otimes \cdots \otimes u_{i_\nu} \in M_{(\ol{\nu},\ol{1})}.
$$
These imply that
$$
u^\nu = \hf(u^\nu + u^{-\nu}) + \hf(u^\nu - u^{-\nu}) \in M_{(\ol{\nu}, \ol{0})} \oplus M_{(\ol{\nu}, \ol{1})}
$$
is the $X^\imath$-weight vector decomposition of $u^\nu$.

On the other hand, by the definition of $u^{\pm \nu}$ and Proposition \ref{Yn on vnu+}, we have
$$
X u^\nu = 0, \qu l u^\nu = q^\nu u^\nu, \qu Y^{(2\nu)}u^\nu = (-1)^\nu u^{-\nu}.
$$
Therefore, the $\Ui$-submodule $\Ui u^\nu$ generated by $u^\nu$ is isomorphic to $V(\nu)$, and it is also generated by two $X^\imath$-weight vectors $u^\nu + u^{-\nu}$ and $u^\nu - u^{-\nu}$. Hence, $\Ui u^\nu$ is a standard $X^\imath$-weight module by Propositions \ref{Ui-submodule generated by Xi-weight vectors is an Xi-weight submodule}, and \ref{Xi-weight submodule of a standard Xi-weight module is standard}. This completes the proof.
\end{proof}

Now, we set
$$
V(\nu)_{\bfA} := \Uidot_{\bfA} v_\nu,
$$
and call it the $\bfA$-form of $V(\nu)$. That it is actually an $\bfA$-form of $V(\nu)$ in the sense of Definition \ref{A-form of Ui-modules} will be proved later.

\begin{lem}\label{ibar-invariance of unu}
Let $\nu \in X_{\frk,\Int}^+$, $M := V_\natural^{\otimes \nu}$. Let $u^\nu \in M$ be as before. Then, we have
$$
\psii_M(u^\nu) = u^\nu.
$$
\end{lem}

\begin{proof}
By Proposition \ref{spectra of Vnatural tensor}, we see that $M$ has no $X_{\frk}$-weight vectors of weight greater than $\nu$ or less than $-\nu$.

The vector $u^\nu$, expanded by the basis $\{ u_{i_1} \otimes \cdots \otimes u_{i_\nu} \mid i_1,\ldots,i_\nu \in \{ 1,2,3 \} \}$, contains $u_2^{\otimes \nu}$ with coefficient $1$. Note that we have
\begin{align}
\begin{split}
&\Etil_1(\ol{u}_1^{\otimes \nu}) = 0, \\
&\Ftil_1^{\nu-1}(\ol{u}_1^{\otimes \nu}) = \ol{u}_1 \otimes \ol{u}_2^{\otimes \nu-1}, \\
&\Ftil_1^{\nu}(\ol{u}_1^{\otimes \nu}) = \ol{u}_2^{\otimes \nu}.
\end{split} \nonumber
\end{align}
This together with Proposition \ref{weights for icanonical basis vector for V(l)}, implies that $G^\imath(\ol{u}_1 \otimes \ol{u}_2^{\otimes \nu-1} + \ol{u}_2^{\otimes \nu})$ is an $X_{\frk}$-weight vector of weight $\nu$. Furthermore, the vector $G^\imath(\ol{u}_1 \otimes \ol{u}_2^{\otimes \nu-1} + \ol{u}_2^{\otimes \nu})$, expanded by the basis $\{ u_{i_1} \otimes \cdots \otimes u_{i_\nu} \mid i_1,\ldots,i_\nu \in \{ 1,2,3 \} \}$, contains $u_2^{\otimes \nu}$ with coefficient $1$. Since $\dim M_\nu = 1$ by Proposition \ref{spectra of Vnatural tensor}, we conclude that
$$
u^\nu = G^\imath(\ol{u}_1 \otimes \ol{u}_2^{\otimes \nu-1} + \ol{u}_2^{\otimes \nu}),
$$
and hence, it is $\psii_{M}$-invariant. This completes the proof
\end{proof}

\begin{prop}\label{information about Vnu}
Let $\nu \in X_{\frk,\Int}^+$.
\begin{enumerate}
\item\label{information about Vnu 1} For each $k \geq 0$, we have
$$
B_{2,p(\nu)}^{(k)} (v_\nu + Y^{(2\nu)}v_\nu), \ B_{2,q(\nu)}^{(k)} (v_\nu - Y^{(2\nu)}v_\nu) \in V(\nu)_{\bfA}.
$$
\item\label{information about Vnu 2} We have $B_{2,p(\nu)}^{(k)} (v_\nu + Y^{(2\nu)}v_\nu) = 0$ for all $k > \nu$.
\item\label{information about Vnu 3} We have $B_{2,q(\nu)}^{(k)} (v_\nu - Y^{(2\nu)}v_\nu) = 0$ for all $k > \nu-1$.
\item\label{information about Vnu 4} $\{ B_{2,p(\nu)}^{(k)} (v_\nu + Y^{(2\nu)}v_\nu) \mid k \in [0,\nu] \} \sqcup \{ B_{2,q(\nu)}^{(k-1)} (v_\nu - Y^{(2\nu)}v_\nu) \mid k \in [0,\nu-1] \}$ forms a free $\bfK_\infty$-basis of $\clL(\nu)$.
\item\label{information about Vnu 5} $\{ \ev_\infty(B_{2,p(\nu)}^{(k)} (v_\nu + Y^{(2\nu)}v_\nu)) \mid k \in [0,\nu] \} \sqcup \{ \ev_\infty(B_{2,q(\nu)}^{(k-1)} (v_\nu - Y^{(2\nu)}v_\nu)) \mid k \in [0,\nu-1] \}$ forms a $\C$-basis of $\ol{\clL}(\nu)$.
\item\label{information about Vnu 6} We have $\ev_\infty(B_{2,p(\nu)}^{(k)} (v_\nu + Y^{(2\nu)}v_\nu)) = \Ytil^k b_\nu + \Ytil^{2\nu-k}b_\nu$ for all $k \in [0,\nu]$.
\item\label{information about Vnu 7} We have $\ev_\infty(B_{2,q(\nu)}^{(k)} (v_\nu - Y^{(2\nu)}v_\nu)) = \Ytil^k b_\nu - \Ytil^{2\nu-k}b_\nu$ for all $k \in [0,\nu-1]$.
\end{enumerate}
\end{prop}

\begin{proof}
Let $M = V_\natural^{\otimes \nu}$, and $u^\nu \in M$ be as before. Then, there exists a $\Ui$-module homomorphism $\phi : V(\nu) \rightarrow M$ such that $\phi(v_\nu) = \frac{1}{\sqrt{2}} u^\nu$. We identify $V(\nu)$ with $\phi(V(\nu))$. By Proposition \ref{1zeta on Ui2-highest weight vector}, we have
\begin{align}
\begin{split}
&B_{2,p(\nu)}^{(k)} (v_\nu + Y^{(2\nu)}v_\nu) = B_{2,(\ol{\nu},\ol{\nu})}^{(k)} (v_\nu + Y^{(2\nu)}v_\nu) = 2B_{2,(\ol{\nu},\ol{\nu})}^{(k)} v_\nu \in V(\nu)_{\bfA}, \\
&B_{2,q(\nu)}^{(k)} (v_\nu - Y^{(2\nu)}v_\nu) = B_{2,(\ol{\nu},\ol{\nu-1})}^{(k)} (v_\nu - Y^{(2\nu)}v_\nu) = 2B_{2,(\ol{\nu},\ol{\nu-1})}^{(k)} v_\nu \in V(\nu)_{\bfA}.
\end{split} \nonumber
\end{align}
This proves assertion \eqref{information about Vnu 1}.

Exchanging the roles of $B_1$ with $B_2$ in Proposition \ref{spectra of Vnatural tensor}, we see that the $B_2$-eigenvalues of $M$ are $[a]$ with $a \in [-\nu,\nu]$. Since $v_\nu + Y^{(2\nu)}v_\nu \in M_{(\ol{\nu},\ol{\nu})}$, by Lemma \ref{characterization of standard Xi weight}, it is a sum of $B_2$-eigenvectors of eigenvalues $[a]$ with $a \in \Z_{p(\nu)}$. Therefore, we see that
$$
B_{2,p(\nu)}^{(k)} (v_\nu + Y^{(2\nu)}v_\nu) = 0 \qu \IF k > \nu.
$$
Similarly, since $v_\nu - Y^{(2\nu)}v_\nu \in M_{(\ol{\nu},\ol{\nu-1})}$, it is a sum of $B_2$-eigenvectors of eigenvalues $[a]$ with $a \in \Z_{q(\nu)}$. Therefore, we see that
$$
B_{2,q(\nu)}^{(k)} (v_\nu - Y^{(2\nu)}v_\nu) = 0 \qu \IF k > \nu-1.
$$
These prove assertions \eqref{information about Vnu 2} and \eqref{information about Vnu 3}.

For each $p \in \{ \ev,\odd \}$ and $k \in [0,\nu]$, we show by induction on $k$ that there exist $d_{m,k,p} \in \bfK_\infty$ such that $d_{0,k,p} = 1$ if $k < \nu$, $d_{0,\nu,p} = \frac{1}{\sqrt{2}}$, and $d_{m,k,p} \in q\inv \bfK_\infty$ if $m > 0$, and
$$
B_{2,p}^{(k)} v_\nu = \sum_{m=0}^{\lfloor \frac{k}{2} \rfloor} d_{m,k,p} \Ytil^{k-2m}v_\nu.
$$
Let $k \in [0,\nu-1]$. Then, inductively, we have
\begin{align}
\begin{split}
B_2 B_{2,p}^{(k)}v_\nu &= (Y+X)\frac{1}{\{l;0\}} \sum_{m=0}^{\lfloor \frac{k}{2} \rfloor} d_{m,k,p} c_{k-2m,\nu} Y^{(k-2m)} v_\nu \\
&= \sum_{m=0}^{\lfloor \frac{k}{2} \rfloor} \frac{d_{m,k,p} c_{k-2m,\nu}}{\{ \nu-k+2m \}} ([k-2m+1] Y^{(k-2m+1)} + [2\nu-k+2m+1]Y^{(k-2m-1)})v_\nu \\
&= \sum_{m=0}^{\lfloor \frac{k}{2} \rfloor} \frac{d_{m,k,p} c_{k-2m,\nu}}{\{ \nu-k+2m \}} (\frac{[k-2m+1]}{c_{k-2m+1,\nu}} \Ytil^{k-2m+1} + \frac{[2\nu-k+2m+1]}{c_{k-2m-1,\nu}} \Ytil^{k-2m-1})v_\nu.
\end{split} \nonumber
\end{align}
By degree consideration, the right-hand side is a sum of $\Ytil^{k-2m+1}$ with $m \in [0,\lfloor \frac{k+1}{2} \rfloor]$ whose coefficient equals $\frac{[k+1]}{\sqrt{2}^{\delta_{k,\nu-1}}}$ if $m = 0$, and belongs to $q^{k-1} \bfK_\infty$ if $m \neq 0$. Then, our claim follows from the recursive formula for the $\imath$divided powers in Proposition \ref{recursive formula for idivided powers}. Similarly, for each $k \in [0,\nu]$ we have
$$
B_{2,p}^{(k)} Y^{(2\nu)}v_\nu = \sum_{m=0}^{\lfloor \frac{k}{2} \rfloor} d_{m,k,p} \Ytil^{2\nu-k+2m}v_\nu
$$
for the same $d_{m,k,p} \in \bfK_\infty$ as before. Therefore, we obtain
$$
\ev_\infty(B_{2,p(\nu)}^{(k)}(v_\nu + Y^{(2\nu)}v_\nu)) = \Ytil^k b_\nu + \Ytil^{2\nu-k} b_\nu
$$
for each $k \in [0,\nu]$, and
$$
\ev_\infty(B_{2,q(\nu)}^{(k)}(v_\nu - Y^{(2\nu)}v_\nu)) = \Ytil^k b_\nu - \Ytil^{2\nu-k} b_\nu
$$
for each $k \in [0,\nu-1]$. These imply the remaining assertions.
\end{proof}

\begin{theo}\label{Vnu is a based Ui-module for n=3}
Let $\nu \in X_{\frk,\Int}^+$. Then, $V(\nu)$ is a based $\Ui$-module. Furthermore, for each $k \in [0,\nu]$, we have
\begin{align}
\begin{split}
&G^\imath(\Ytil^k b_\nu + \Ytil^{2\nu-k} b_\nu) = B_{2,p(\nu)}^{(k)} (v_\nu + Y^{(2\nu)}v_\nu), \\
&G^\imath(\Ytil^k b_\nu - \Ytil^{2\nu-k} b_\nu) = B_{2,q(\nu)}^{(k)} (v_\nu - Y^{(2\nu)}v_\nu).
\end{split} \nonumber
\end{align}
\end{theo}

\begin{proof}
Let $u^\nu \in M := V_\natural^{\otimes \nu}$ be as before. Then, there exists an almost isometry
$$
\phi : V(\nu) \rightarrow \Uidot u^\nu (= \Ui u^\nu)
$$
of $\Ui$-modules such that $\phi(v_\nu) = \frac{1}{\sqrt{2}}u^\nu$. Note that $u^\nu$ belongs to the $\bfA$-form $M_{\bfA} := \bigoplus_{i=1}^n \bfA u_i$. Also, by Lemma \ref{ibar-invariance of unu}, it is $\psii_{M}$-invariant. Then, by Propositions \ref{Ui-homomorphism lifts to Uidot-homomorphism} and \ref{information about Vnu} \eqref{information about Vnu 1}, we have
$$
\phi(B_{2,p(\nu)}^{(k)}(v_\nu + Y^{(2\nu)}v_\nu)), \ \phi(B_{2,q(\nu)}^{(k)}(v_\nu - Y^{(2\nu)}v_\nu)) \in \Uidot_{\bfA} u^\nu \subseteq M_{\bfA},
$$
for all $k \geq 0$, and they are $\psii_M$-invariant. Furthermore, by Proposition \ref{information about Vnu} \eqref{information about Vnu 6} and \eqref{information about Vnu 7}, we obtain
\begin{align}
\begin{split}
&\ev_\infty(\phi(B_{2,p(\nu)}^{(k)}(v_\nu + Y^{(2\nu)}v_\nu))) = \Ytil^k b + \Ytil^{2\nu-k} b \qu \IF k \in [0,\nu], \\
&\ev_\infty(\phi(B_{2,q(\nu)}^{(k)}(v_\nu - Y^{(2\nu)}v_\nu))) = \Ytil^k b - \Ytil^{2\nu-k} b \qu \IF k \in [0,\nu-1],
\end{split} \nonumber
\end{align}
where $b := \ev_\infty(\phi(v_\nu))$. Therefore, by Lemma \ref{characterization of icanonical basis element}, we conclude that
\begin{align}
\begin{split}
&\phi(B_{2,p(\nu)}^{(k)}(v_\nu + Y^{(2\nu)}v_\nu)) = G^\imath(\Ytil^k b + \Ytil^{2\nu-k} b) \qu \IF k \in [0,\nu], \\
&\phi(B_{2,q(\nu)}^{(k)}(v_\nu - Y^{(2\nu)}v_\nu)) = G^\imath(\Ytil^k b - \Ytil^{2\nu-k} b) \qu \IF k \in [0,\nu-1].
\end{split} \nonumber
\end{align}
These show that $\Ui u^\nu$ is spanned by $G^\imath(\ev_\infty(\phi(\clL(\nu))))$. Hence, by Proposition \ref{based submodule and based quotient module} \eqref{based submodule and based quotient module 1}, $\Ui u^\nu$ is a based submodule of $M$. This shows that $V(\nu) = \phi\inv(\Ui u^\nu)$ is a based $\Ui$-module. Thus, the proof completes.
\end{proof}

\begin{cor}\label{Btil on Vnu for rank 2}
Let $\nu \in X_{\frk,\Int}^+$.
\begin{enumerate}
\item\label{Btil on Vnu for rank 2 1} For each $k \in [0,\nu]$, we have
$$
\Btil_1(\Ytil^k + \Ytil^{2\nu-k})b_\nu = (\Ytil^k - \Ytil^{2\nu-k})b_\nu,
$$
and for each $k \in [0,\nu-1]$,
$$
\Btil_1(\Ytil^k - \Ytil^{2\nu-k})b_\nu = (\Ytil^k + \Ytil^{2\nu-k})b_\nu.
$$
\item\label{Btil on Vnu for rank 2 2} For each $k \in [0,\nu]$, we have
$$
\Btil_2(\Ytil^k + \Ytil^{2\nu-k})b_\nu = \begin{cases}
(\Ytil^{k-1} + \Ytil^{2\nu-k+1})b_\nu \qu & \IF \nu-k \in \Z_{\ev} \setminus \{0\}, \\
\sqrt{2} (\Ytil^{\nu-1} + \Ytil^{\nu+1}) b_\nu \qu & \IF \nu-k = 0, \\
(\Ytil^{k+1} + \Ytil^{2\nu-k-1})b_\nu \qu & \IF \nu-k \in \Z_{\odd} \setminus \{1\}, \\
\sqrt{2} \Ytil^\nu b_\nu \qu & \IF \nu-k = 1,
\end{cases}
$$
and for each $k \in [0,\nu-1]$,
$$
\Btil_2(\Ytil^k - \Ytil^{2\nu-k})b_\nu = \begin{cases}
(\Ytil^{k+1} - \Ytil^{2\nu-k-1})b_\nu \qu & \IF \nu-k \in \Z_{\ev}, \\
(\Ytil^{k-1} - \Ytil^{2\nu-k+1})b_\nu \qu & \IF \nu-k \in \Z_{\odd},
\end{cases}
$$
\end{enumerate}
\end{cor}

\begin{ex}\normalfont
We give graphical descriptions of $\Btil_1$ and $\Btil_2$ on $\ol{\clL}(\nu)$ for $\nu = 2$ (left) and $\nu = 3$ (right).
\begin{align}
\xymatrix{
(\Ytil^{0} + \Ytil^{4})b_2 \ar@{<->}[r]^-{\Btil_1} & (\Ytil^{0} - \Ytil^{4})b_2 \ar@{<->}[d]^-{\Btil_2} \\
(\Ytil^{1} + \Ytil^{3})b_2 \ar@{<->}[r]^-{\Btil_1} \ar@{<->}[d]^-{\Btil_2} & (\Ytil^{1} - \Ytil^{3})b_2 \\
\sqrt{2}\Ytil^{2}b_2 
} \qu \qu \qu
\xymatrix{
(\Ytil^{0} + \Ytil^{6})b_3 \ar@{<->}[r]^-{\Btil_1} \ar@{<->}[d]^-{\Btil_2} & (\Ytil^{0} - \Ytil^{6})b_3 \\
(\Ytil^{1} + \Ytil^{5})b_3 \ar@{<->}[r]^-{\Btil_1} & (\Ytil^{1} - \Ytil^{5})b_3 \ar@{<->}[d]^-{\Btil_2} \\
(\Ytil^{2} + \Ytil^{4})b_3 \ar@{<->}[r]^-{\Btil_1} \ar@{<->}[d]^-{\Btil_2} & (\Ytil^{2} - \Ytil^{4})b_3 \\
\sqrt{2}\Ytil^{3}b_3 
}
\nonumber
\end{align}
\end{ex}

\begin{prop}\label{Xi weight decomposition of highest weight vector}
Let $M$ be a finite-dimensional classical weight standard $X^\imath$-weight module. Let $v \in M$ be a highest weight vector of weight $\nu \in X_{\frk,\Int}^+$. Then, we have
$$
\Ui v = \Uidot v.
$$
Furthermore, for each $\zeta = (\zeta_1,\zeta_2) \in X^\imath$, we have
\begin{align}
\begin{split}
&\mathbf{1}_\zeta v = \begin{cases}
\hf(v + Y^{(2\nu)} v) \qu & \IF \zeta = (\ol{\nu},\ol{\nu}), \\
\hf(v - Y^{(2\nu)} v) \qu & \IF \zeta = (\ol{\nu},\ol{\nu-1}), \\
0 \qu & \OW.
\end{cases} 
\end{split} \nonumber
\end{align}
\end{prop}

\begin{proof}
The assertions follow from Propositions \ref{Ui-homomorphism lifts to Uidot-homomorphism} and \ref{1zeta on Ui2-highest weight vector}.
\end{proof}

\begin{prop}\label{preparation for branching rule n=3}
Let $M$ be a finite-dimensional classical weight module equipped with a contragredient Hermitian inner product. Set
\begin{align}
\begin{split}
&L_1 := \{ b \in \ol{\clL}_M \mid \Xtil b = 0 \}, \\
&L_2 := \{ b \in \ol{\clL}_M \mid \Btil_2 b = 0 \}.
\end{split} \nonumber
\end{align}
Then, the linear map
$$
L_2 \rightarrow L_1;\ b \mapsto (1 + \Btil_1)b
$$
is an isomorphism of $\C$-vector spaces, with inverse
$$
L_1 \rightarrow L_2;\ b \mapsto \sum_{\substack{\zeta \in X^\imath \\ \zeta_2 = \ol{0}}} \mathbf{1}_\zeta b.
$$
\end{prop}

\begin{proof}
Write an orthogonal irreducible decomposition of $M$ as
$$
M = \bigoplus_{k=1}^r M_k, \qu M_k \simeq V(\nu_k), \qu \nu_k \in X_{\frk,\Int}^+.
$$
This induces an orthogonal decomposition
$$
\ol{\clL}_M = \bigoplus_{k=1}^r \ol{\clL}_k, \qu \ol{\clL}_k := \ol{\clL}_{M_k} \simeq \ol{\clL}(\nu_k).
$$
Since both $\Xtil$ and $\Btil_2$ preserve $\ol{\clL}_k$ for all $k \in [1,r]$, it suffices to prove the assertion for the $M = V(\nu)$, $\nu \in X_{\frk,\Int}^+$ case.

When $M = V(\nu)$ for some $\nu \in X_{\frk,\Int}^+$, we have
$$
L_1 = \C b_\nu.
$$
Furthermore, by Corollary \ref{Btil on Vnu for rank 2} \eqref{Btil on Vnu for rank 2 2}, we see that
$$
L_2 = \C (1 + (-1)^\nu \Ytil^{2\nu}) b_\nu.
$$
Using Proposition \ref{Xi weight decomposition of highest weight vector}, we compute as
\begin{align}
\begin{split}
(1 + \Btil_1)(\sum_{\zeta_2 = \ol{0}} \mathbf{1}_\zeta b_\nu) &= \hf(1+\Btil_1)(1 + (-1)^{\nu} \Ytil^{2\nu}) b_\nu \\
&= \hf(1 + (-1)^{\nu} \Ytil^{2\nu})b_\nu + \hf(1 - (-1)^{\nu} \Ytil^{2\nu})b_\nu \\
&= b_\nu.
\end{split} \nonumber
\end{align}
This, together with dimension consideration, proves the assertion.
\end{proof}

\subsection{Applications to the general $n \geq 3$ case}
In this subsection, we consider the general $n \geq 3$ case; results obtained here are necessary for the $n = 4$ case. For a subset $J = \{ j_1,\ldots,j_r \} \subseteq I$, let $\U_J = \U_{j_1,\ldots,j_r}$ (resp., $\Ui_J = \Ui_{j_1,\ldots,j_r}$) denote the subalgebra generated by $E_j,F_j,K_j^{\pm 1}$, $j \in J$ (resp., $B_j$, $j \in J$). Then, $(J,\emptyset,\id)$ is a Satake subdiagram of $(I,\emptyset,\id)$ of type AI. Let $\Uidot_J = \Uidot_{j_1,\ldots,j_r} := \Uidot(J)$ denote the corresponding modified $\imath$quantum group.

For each $\eta = (\eta_j)_{j \in J} \in (\Z/2\Z)^J = X^\imath(J)$, let $\mathbf{1}_{J,\eta} = \mathbf{1}_{j_1,\ldots,j_r,\eta}$ denote the corresponding idempotent in $\Uidot_J$. It acts on each $X^\imath$-weight module by
$$
\mathbf{1}_{J,\eta} = \sum_{\substack{\zeta \in X^\imath \\ \zeta_j = \eta_j \Forall j \in J}} \mathbf{1}_\zeta.
$$
In particular, for each $\zeta = (\zeta_1,\ldots,\zeta_{n-1}) \in X^\imath$, we have
\begin{align}\label{factorization of idempotent}
\mathbf{1}_{1,\zeta_1} \cdots \mathbf{1}_{n-1,\zeta_{n-1}} = \mathbf{1}_\zeta
\end{align}
on each $X^\imath$-weight module.

\begin{lem}\label{V(nu) is standard if it is contained in a standard module}
Let $M$ be a finite-dimensional classical weight standard $X^\imath$-weight module. Let $v \in M$ be a highest weight vector of weight $\nu \in X_{\frk,\Int}^+$. Then, $\Ui v$ is a standard $X^\imath$-weight module.
\end{lem}

\begin{proof}
For each $\zeta \in X^\imath$, by equation \eqref{factorization of idempotent}, we have
$$
\mathbf{1}_\zeta v = \mathbf{1}_{1,\zeta_1} \cdots \mathbf{1}_{n-1,\zeta_{n-1}} v.
$$
For a proof, by Propositions \ref{Ui-submodule generated by Xi-weight vectors is an Xi-weight submodule} and \ref{Xi-weight submodule of a standard Xi-weight module is standard}, it suffices to show that $\mathbf{1}_\zeta v \in \Ui v$ for all $\zeta \in X^\imath$.

Let $\zeta \in X^\imath$ be such that $\mathbf{1}_\zeta v \neq 0$. We show that $\ol{\nu}_{2i-1} = \zeta_{2i-1}$ for all $i \in I_{\frk}$ and
$$
\mathbf{1}_\zeta v = \frac{1}{2^{m'}} (1-(-1)^{\delta_{\zeta_1,\zeta_2}} Y_2^{(2\nu_1)}) (1-(-1)^{\delta_{\zeta_3,\zeta_4}} Y_4^{(2\nu_3)}) \cdots (1-(-1)^{\delta_{\zeta_{2m'-1},\zeta_{2m'}}} Y_{2m'}^{(2\nu_{2m'-1})}) v,
$$
where
$$
m' := \begin{cases}
\frac{n}{2}-1 \qu & \IF n \in \Z_{\ev}, \\
\frac{n-1}{2} \qu & \IF n \in \Z_{\odd}.
\end{cases}
$$

We proceed by induction on $n$. The case when $n = 3$ has been already proved in Proposition \ref{Xi weight decomposition of highest weight vector}.

First, consider the case when $n \in \Z_{\ev}$. Since $B_{n-1} v = [\nu_{n-1}]v$, we obtain
$$
\mathbf{1}_{n-1,\zeta_{n-1}} v = \delta_{\ol{\nu_{n-1}},\zeta_{n-1}} v.
$$
This, together with our hypothesis that $\mathbf{1}_{\zeta} v \neq 0$, implies that $\ol{\nu_{n-1}} = \zeta_{n-1}$. As a $\Ui_{1,2,\ldots,n-2}$-module vector, $\mathbf{1}_{n-1,\zeta_{n-1}} v = v$ is a highest weight vector of weight $(\nu_1,\nu_3,\ldots,\nu_{n-3})$. Therefore, our claim follows from the induction hypothesis.

Next, consider the case when $n \in \Z_{\odd}$. By the $n = 3$ case, we have
$$
\mathbf{1}_{n-2,\zeta_{n-2}} \mathbf{1}_{n-1,\zeta_{n-1}} v = \delta_{\ol{\nu_{n-2}},\zeta_{n-2}} \hf(1 - (-1)^{\delta_{\zeta_{n-2},\zeta_{n-1}}} Y_{n-1}^{(2\nu_{n-2})})v.
$$
As before, this implies that $\ol{\nu_{n-2}} = \zeta_{n-2}$. By weight consideration, as a $\Ui_{1,2,\ldots,n-3}$-module vector, both $v$ and $Y_{n-1}^{(2\nu_{n-2})}v$ are highest weight vectors of weight $(\nu_1,\nu_3,\ldots,\nu_{n-4})$. Then, our claim follows from the induction hypothesis. Thus, the proof completes.
\end{proof}

\begin{prop}
Let $\nu \in X_{\frk,\Int}$. Then, $V(\nu)$ is a standard $X^\imath$-weight module.
\end{prop}

\begin{proof}
The assertion follows from Propositions \ref{Tensor power contains all highest weight vectors}, \ref{canonical Xi-weight structure is standard}, and Lemma \ref{V(nu) is standard if it is contained in a standard module}.
\end{proof}

For each $\nu \in X_{\frk,\Int}$, set
$$
V(\nu)_{\bfA} := \UidotA v_\nu,
$$
and call it the $\bfA$-form of $V(\nu)$.

For a based classical weight module $M$, consider the following condition:
\begin{align}\label{Property3}
G^\imath(b) \in M_\nu \oplus \bigoplus_{\substack{\xi \in X_{\frk,\Int} \setminus \{\nu\} \\ |\xi_i| \geq |\nu_i| \Forall i \in I_{\otimes}}} M_\xi \qu \Forall \nu \in X_{\frk,\Int}, b \in \ol{\clL}_{M,\nu}.
\end{align}
When $M$ satisfies this condition, for each $\nu \in X_{\frk,\Int}$ and $b \in \ol{\clL}_{M,\nu}$, set $G^\imath_0(b)$ to be the image of $G^\imath(b)$ under the projection onto the weight space $M_\nu$. This defines a $\C$-linear map
$$
G^\imath_0 : \ol{\clL}_M \rightarrow M
$$
which preserves the $X_{\frk}$-weight spaces. Note that we have
$$
\ev_\infty(G^\imath_0(b)) = b
$$
for all $b \in \ol{\clL}_M$.

\begin{prop}
Let $(M, (\cdot,\cdot)_M, M_{\bfA}, \psi_M, \clB_M)$ be a finite-dimensional based $\U$-module such that $\psii_M := \Upsilon \circ \psi_M$ is defined. Then, $(M, (\cdot,\cdot)_M, M_{\bfA}, \psii_M,\clB_M)$ is a based classical weight standard $X^\imath$-weight module satisfying condition \eqref{Property3}.
\end{prop}

\begin{proof}
Let $\clB$ be a crystal basis of $M$. Let $i \in I$, and set
$$
\clB_i := \{ b \in \clB \mid \Etil_i b = 0 \}.
$$
Then, we have
$$
\clB = \{ \Ftil_i^k b \mid b \in \clB_i,\ k \in [0,\vphi_i(b)] \}.
$$
Let $\preceq$ be a total ordering on $\clB_i$ such that for each $b_1,b_2 \in \clB_i$, $\wt(b_1) < \wt(b_2)$ implies $b_1 \prec b_2$. For each $b \in \clB_i$, set
\begin{align}
\begin{split}
&\clB_{\succeq b} := \{ \Ftil_i^k b' \in \clB \mid b' \in \clB_i,\ b' \succeq b,\ k \in [0,\vphi_i(b')] \}, \\
&\clB_{\succ b} := \{ \Ftil_i^k b' \in \clB \mid b' \in \clB_i,\ b' \succ b,\ k \in [0,\vphi_i(b')] \}, \\
&\clB_b := \clB_{\succeq b} \setminus \clB_{\succ b} = \{ \Ftil_i^k b \mid k \in [0,\vphi_i(b)] \}.
\end{split} \nonumber
\end{align}
Also, for each $b \in \clB_i$ and $\nu \in [-\vphi_i(b), \vphi_i(b)]_{p(\vphi_i(b))}$, set
$$
b_\nu := \begin{cases}
b \qu & \IF \nu = 0, \\
\frac{1}{\sqrt{2}}(\Ftil_i^{|\nu|-1} b \pm \Ftil_i^{|\nu|} b) \qu & \IF \pm \nu > 0.
\end{cases}
$$

By Proposition \ref{cellularity of canonical basis of based U-module}, for each $b \in \clB_i$, the $\U_i$-submodule $M_{\succeq b}$ (resp., $M_{\succ b}$) spanned by $\{ G(b') \mid b' \in \clB_{\succeq b} \}$ (resp., $\{ G(b') \mid b' \in \clB_{\succ b} \}$) is a based $\U_i$-submodule with a crystal basis $\clB_{\succeq b}$ (resp., $\clB_{\succ b}$). Furthermore, as a $\Ui_i$-submodule, it is a based $\Ui_i$-submodule by Lemma \ref{weight spaces iCB elements live} and Proposition \ref{based submodule and based quotient module} \eqref{based submodule and based quotient module 1}.

Now, we prove by descending induction on $b \in \clB_i$ that $M_{\succeq b}$ satisfies condition \eqref{Property3}. Let $b \in \clB_i$, and assume that our claim is true for all $b' \in \clB_i$ such that $b' \succ b$. Then, the quotient $\Ui_i$-module $M_b := M_{\succeq b}/M_{\succ b}$ is a based $\Ui_i$-module (by Proposition \ref{based submodule and based quotient module} \eqref{based submodule and based quotient module 2}) isomorphic to the $(l+1)$-dimensional irreducible $\U_i$-module $V(l)$, where $l := \vphi_i(b)$. Then, for each $\nu \in [-l,l]_{p(l)}$, we have
$$
B_{i,\nu} [G^\imath(b)] = [G^\imath(b_\nu)],
$$
where $[v]$ denotes the image of $v \in M_{\succeq b}$ in $M_b$, and $B_{i,\nu}$ is as in the proof of Proposition \ref{weights for icanonical basis vector for V(l)}. Let us write
$$
B_{i,\nu} G^\imath(b) = G^\imath(b_\nu) + \sum_{\substack{b' \in \clB_i \\ b' \succ b}} \sum_{\xi \in [-\vphi_i(b'), \vphi_i(b')]_{p(\vphi_i(b'))}} c_{b',\xi} G^\imath(b'_\xi),
$$
where $c_{b',\xi} \in \bfA_{\Inv} := \{ f \in \bfA \mid f(q\inv) = f(q) \}$. We show that $c_{b',\xi} = 0$ for all $\xi \notin \{ \nu, \pm(|\nu| + 2), \pm(|\nu| + 4), \ldots \}$. Assume contrary and take $(b',\xi)$ such that $\deg(c_{b',\xi})$ is maximal. Let us write
$$
G^\imath(b'_\xi) = \sum_{\omega \in X_{\frk}} m_\omega, \qu m_\omega \in M_\omega.
$$
Then, we have $m_\omega \in \clL_M$, $\ev_\infty(m_\omega) = \delta_{\omega,\xi} b'_\xi$. Consider
$$
(B_{i,\nu} G^\imath(b), m_\xi)_M.
$$
Since the first factor is a sum of $B_i$-eigenvector of eigenvalues in $\{ \nu, \pm(|\nu|+2), \pm(|\nu|+4), \ldots \}$, this value equals $0$. On the other hand, by the maximality of $\deg(c_{b',\xi})$, this value belongs to $c_{b',\xi} + q^{\deg(c_{b',\xi})-1} \bfK_\infty$. This contradicts our assumption that $c_{b',\xi} \neq 0$.

By above, we obtain that
$$
G^\imath(b_\nu) = B_{i,\nu} G^\imath(b) - \sum_{\substack{b' \in \clB_i \\ b' \succ b}} \sum_{\xi \in \{ \nu, \pm(|\nu|+2),\ldots \}} c_{b',\xi} G^\imath(b'_\xi).
$$
By the definition of $B_{i,\nu}$ and our induction hypothesis, the right-hand side is a linear combination of $B_i$-eigenvectors of eigenvalues in $\{ \nu, \pm(|\nu|+2),\ldots \}$. Therefore, we conclude that $M_{\succeq b}$ satisfies condition \eqref{Property3}, as desired. This completes the proof.
\end{proof}

\begin{lem}\label{cellularity of highest component}
Let $M$ be a finite-dimensional based classical weight standard $X^\imath$-weight module satisfying condition \eqref{Property3}. Let $\nu \in X_{\frk,\Int}$ be such that $M[> \nu] = 0$. Assume that $V(\nu)$ is a based $\Ui$-module. Then, there exists a based $\Ui$-submodule $N$ of $M$ isomorphic to $V(\nu)$.
\end{lem}

\begin{proof}
Let $b \in \ol{\clL}_M[\nu] \cap \ol{\clL}_{M,\nu}$. Then, by condition \eqref{Property3} and our hypothesis that $M[> \nu] = 0$, the vector $G^\imath(b)$ is a highest weight vector of weight $\nu$. Set $N := \Ui G^\imath(b)$. Then, there exists an isomorphism $\phi : V(\nu) \rightarrow N$ of $\Ui$-modules which sends $v_\nu$ to $G^\imath(b)$. It almost preserves the metrics. In particular, it induces an isomorphism $\phi : \ol{\clL}(\nu) \rightarrow \ol{\clL}_N \subseteq \ol{\clL}_M$ of $\C$-vector spaces.

For each $b' \in \ol{\clL}(\nu)$, there exists $x_{b'} \in \UidotA$ such that $x_{b'} v_\nu = G^\imath(b')$. Replacing $x_{b'}$ with $\hf(x_{b'} + \psii(x_{b'}))$, we may assume that $\psii(x_{b'}) = x_{b'}$. Then, we have
$$
\psii_M(\phi(G^\imath(b'))) = \psii_M(\phi(x_{b'} v_\nu)) = \psii_M(x_{b'} G^\imath(b)) = x_{b'} G^\imath(b) = \phi(x_{b'} v_\nu) = \phi(G^\imath(b')).
$$
Furthermore, since $\phi$ almost preserves the metrics, we obtain
$$
\ev_\infty(\phi(G^\imath(b'))) = \phi(\ev_\infty(G^\imath(b'))) = \phi(b').
$$
Then, by Lemma \ref{characterization of icanonical basis element}, we conclude that
$$
\phi(G^\imath(b')) = G^\imath(\phi(b')).
$$

This far, we have shown that $N$ is spanned by $G^\imath(\phi(\ol{\clL}(\nu)))$. Since $\phi(\ol{\clL}(\nu)) = \ol{\clL}_N$, by Proposition \ref{based submodule and based quotient module} \eqref{based submodule and based quotient module 1}, we finally see that $N$ is a based $\Ui$-submoule of $M$.
\end{proof}

\section{$n = 4$ case}\label{section: n=4}
In this section, we consider the $n=4$ case. In this case, we can identify $X = \Z^3$, $X^+ = \Z^3_{\geq 0}$, $X_{\frk,\Int} = \Z^2$, $X_{\frk,\Int}^+ = \{ (\nu_1,\nu_3) \in X_{\frk,\Int} \mid \nu_1 \geq |\nu_3| \}$, and $X^\imath = (\Z/2\Z)^3$.

\subsection{Lowering and raising operators}
Set
$$
X_{\pm} := X_{2,\pm} = B_{2,+\pm}\{ l_1;0 \}, \qu Y_\pm := Y_{2,\pm} = B_{2,-\mp}\{ l_3;0 \}.
$$

Let $\nu = (\nu_1,\nu_3) \in X_{\frk,\Int}^+$, and consider the irreducible $\Ui$-module $V(\nu)$.

\begin{lem}\label{B_+ B_-^n AI-2}
For each $l,m \in \Z_{\geq 0}$, we have
\begin{align}
\begin{split}
&Y_+^{(m)}Y_-^{(l)} v_\nu = Y_-^{(l)}Y_+^{(m)}v_\nu, \\
&X_- Y_-^{(l)}Y_+^{(m)} v_\nu = \frac{[\nu_1-\nu_3-l+1]\{ \nu_1-l+1 \}\{ \nu_3+l \}}{\{ \nu_1-l-m+1 \}\{ \nu_3+l-m \}} Y_-^{(l-1)} Y_+^{(m)} v_\nu, \\
&X_+ Y_-^{(l)}Y_+^{(m)} v_\nu = \frac{[\nu_1+\nu_3-m+1]\{ \nu_1-m+1 \}\{ \nu_3-m \}}{\{ \nu_1-l-m+1 \}\{ \nu_3+l-m \}} Y_-^{(l)} Y_+^{(m-1)} v_\nu,
\end{split} \nonumber
\end{align}
where $Y_\pm^{(k)} := \frac{1}{[k]!} Y_\pm^k$.
\end{lem}

\begin{proof}
By equations \eqref{formula for AI 6}--\eqref{formula for AI 10}, we have
\begin{align}
\begin{split}
&[Y_+,Y_-] = 0, \\
&[X_\pm,Y_\mp\frac{\{ l_1;0 \}}{\{ l_3;0 \}}] = 0, \\
&[X_+,Y_+] = [l_1l_3;0] + (q-q\inv)^2 Y_-X_- \frac{[l_1l_3;0]}{\{ l_1;0 \}\{ l_3;-1 \}}, \\
&[X_-,Y_-] = [l_1l_3\inv;0] + (q-q\inv)^2 Y_+X_+ \frac{[l_1l_3\inv;0]}{\{ l_1;0 \}\{ l_3;1 \}}.
\end{split} \nonumber
\end{align}
For each $m \geq 0$,
\begin{align}
\begin{split}
&X_\pm Y_\mp^{(m)} = Y_\mp^{(m)} X_\pm \frac{\{ l_1;1 \}\{ l_3;0 \}}{\{ l_1;-m+1 \}\{ l_3;\pm m \}}, \\
&X_{\pm} Y_{\pm}^{(m)} = Y_{\pm}^{(m-1)}[l_1l_3^{\pm 1}; -m+1] + Y_\pm^{(m)}X_\pm + Y_\pm^{(m-1)} Y_\mp X_\mp P_\pm(m),
\end{split} \nonumber
\end{align}
for some rational function $P_\pm(m) \in \Q(l_1,l_3)$ in variables $l_1,l_3$. Now, it is straightforward to obtain the desired equations by induction on $l,m$.
\end{proof}

From this lemma (with a standard argument), we see that the set
$$
\{ Y_-^{(l)} Y_+^{(m)} v_\nu \mid l \in [0,\nu_1-\nu_3],\ m \in [0,\nu_1+\nu_3] \}
$$
forms a basis of $V(\nu)$. Note that $Y_-^{(l)} Y_+^{(m)} v_\nu \in V(\nu)_{(\nu_1-l-m,\nu_3+l-m)}$. Then, for each $l \in [0,\nu_1-\nu_3]$, $m \in [0,\nu_1+\nu_3]$, we can choose $c_{l,m,\nu} \in \bfK^\times$ in a way such that $\lt(c_{l,m,\nu}Y_-^{(l)} Y_+^{(m)}v_{\nu}) = 1$. For each $l,m \in \Z$, set
$$
\Ytil_-^l \Ytil_+^m v_\nu := \begin{cases}
c_{l,m,\nu} Y_-^{(l)} Y_+^{(m)} v_\nu \qu & \IF l \in [0,\nu_1-\nu_3] \AND m \in [0,\nu_1+\nu_3], \\
0 \qu & \OW.
\end{cases}
$$
Then, $\{ \Ytil_-^{l} \Ytil_+^{m} v_\nu \mid l \in [0,\nu_1-\nu_3],\ m \in [0,\nu_1+\nu_3] \}$ forms an almost orthonormal basis of $V(\nu)$.

Using Lemma \ref{B_+ B_-^n AI-2}, we see by induction on $l$ and $m$ that
\begin{align}
\begin{split}
&(Y_{-}^{(l)}Y_{+}^{(m)}v_\nu,Y_{-}^{(l)}Y_{+}^{(m)}v_\nu)_\nu \\
&\qu = {\nu_1-\nu_3 \brack l}{\nu_1+\nu_3 \brack m} \frac{\{\nu_3-m\}}{\{\nu_3+l-m\}} \prod_{i=1}^l \frac{\{\nu_1-i+1\}\{\nu_3+i\}}{\{\nu_1-i-m+1\}\{\nu_1-i-m\}} \prod_{j=1}^m \frac{\{\nu_3-j+1\}}{\{\nu_1-j\}}.
\end{split} \nonumber
\end{align}

\begin{lem}\label{expansion of Yl}
For each $k \geq 0$, we have
$$
Y^{(k)} = \sum_{l+m=k} Y_-^{(l)} Y_+^{(m)} \frac{\{ l_3;l-m \}\prod_{j=1}^k\{ l_1;-j+1 \}}{\prod_{j=0}^k \{ l_3;l-j \}}.
$$
\end{lem}

\begin{proof}
By equation \eqref{formula for AI 3}, we have
$$
Y = (Y_- + Y_+) \frac{\{ l_1;0 \}}{\{ l_3;0 \}}.
$$
Then, the assertion follows by induction on $k$.
\end{proof}

Set
$$
Y^{(l,m)}v_\nu := \frac{\{\nu_3+l-m\}}{\{\nu_1-l-m\}}\prod_{j=0}^{l+m}\frac{\{\nu_1-j\}}{\{\nu_3+l-j\}} Y_-^{(l)} Y_+^{(m)} v_\nu.
$$
By Lemma \ref{expansion of Yl}, we see that
\begin{align}\label{definition of Ymn vnu}
Y^{(k)} v_\nu = \sum_{l+m = k} Y^{(l,m)} v_\nu.
\end{align}
Also, we have
\begin{align}
\begin{split}
&(Y^{(l,m)}v_\nu,Y^{(l,m)}v_\nu)_\nu \\
&\qu = {\nu_1-\nu_3 \brack l}{\nu_1+\nu_3 \brack m} \frac{\{\nu_1\}\{\nu_3+l-m\}}{\{\nu_3\}\{\nu_1-l-m\}} \prod_{i=1}^l \frac{\{\nu_1-i+1\}}{\{\nu_3+l-i+1\}} \prod_{j=1}^m \frac{\{\nu_1-j+1\}}{\{\nu_3-j\}}.
\end{split} \nonumber
\end{align}
This implies that for each $l \in [0,\nu_1-\nu_3]$ and $m \in [0,\nu_1+\nu_3]$, we have
\begin{align}\label{deg on Ymn vnu}
\begin{split}
&\deg(Y^{(l,m)} v_\nu) \\
&\qu = \hf((\nu_1-\nu_3-l)l + (\nu_1+\nu_3-m)m + \nu_1 - |\nu_3| - |\nu_1-l-m| + |\nu_3+l-m| \\
&\qu +\sum_{i=1}^l (|\nu_1-i+1|-|\nu_3+l-i+1|) + \sum_{j=1}^m (|\nu_1-j+1| - |\nu_3-j|)).
\end{split}
\end{align}

\begin{defi}\normalfont
The lowering operator $\Ytil_{\pm}$ and the raising operator $\Xtil_{\pm}$ are $\bfK$-linear endomorphism on $V(\nu)$ defined by
\begin{align}
\begin{split}
&\Ytil_{+}(\Ytil_{-}^{l} \Ytil_{+}^{m} v_\nu) := \begin{cases}
\Ytil_{-}^{l} \Ytil_{+}^{m+1} v_\nu \qu & \IF l \in [0,\nu_1-\nu_3] \AND m \in [0,\nu_1+\nu_3], \\
0 \qu & \OW,
\end{cases} \\
&\Ytil_{-}(\Ytil_{-}^{l} \Ytil_{+}^{m} v_\nu) := \begin{cases}
\Ytil_{-}^{l+1} \Ytil_{+}^{m} v_\nu \qu & \IF l \in [0,\nu_1-\nu_3] \AND m \in [0,\nu_1+\nu_3], \\
0 \qu & \OW,
\end{cases} \\
&\Xtil_{+}(\Ytil_{-}^{l} \Ytil_{+}^{m} v_\nu) := \begin{cases}
\Ytil_{-}^{l} \Ytil_{+}^{m-1} v_\nu \qu & \IF l \in [0,\nu_1-\nu_3] \AND m \in [0,\nu_1+\nu_3], \\
0 \qu & \OW,
\end{cases} \\
&\Xtil_{-}(\Ytil_{-}^{l} \Ytil_{+}^{m} v_\nu) := \begin{cases}
\Ytil_{-}^{l-1} \Ytil_{+}^{m} v_\nu \qu & \IF l \in [0,\nu_1-\nu_3] \AND m \in [0,\nu_1+\nu_3], \\
0 \qu & \OW.
\end{cases}
\end{split} \nonumber
\end{align}
\end{defi}

Clearly, the operators $\Ytil_\pm$ and $\Xtil_\pm$ preserve $\clL(\nu)$. Hence, they induce linear operators $\Ytil_\pm$ and $\Xtil_\pm$ on $\ol{\clL}(\nu)$.

Since each finite-dimensional classical weight module is completely reducible, we can extend the definitions of $\Ytil_{\pm},\Xtil_{\pm}$ on the irreducible modules to each finite-dimensional classical weight modules. The following can be proved in a similar way to the $n = 3$ case.

\begin{prop}\label{Im Yj = Im Ytilj for rank 2 AI}
Let $M$ be a finite-dimensional classical weight $\Ui$-module. Let $\xi \in X_{\frk,\Int}$. Then, we have $Y_\pm M_\xi = \Ytil_\pm M_\xi$.
\end{prop}

\begin{prop}\label{Ftil and Etil are almost adjoint AI 2}
Let $M$ be a finite-dimensional classical weight $\Ui$-module equipped with a contragredient Hermitian inner product. Then, for each $u,v \in \clL_M$, we have
\begin{align}
\begin{split}
&(\Ytil_\pm u,v) - (u, \Xtil_\pm v) \in q\inv \bfK_\infty.
\end{split} \nonumber
\end{align}
\end{prop}

\subsection{Based module structures of irreducible modules}
For each $l,m \in \Z$, set
$$
\Ytil^{l,m} := \Ytil_-^l \Ytil_+^m.
$$

Let $\nu \in X_{\frk,\Int}^+$, and consider the irreducible module $V(\nu)$. Note that for each $l \in [0,\nu_1-\nu_3]$ and $m \in [0,\nu_1+\nu_3]$, we have
$$
\ev_\infty(\frac{1}{\lt(Y^{(l,m)} v_\nu)} Y^{(l,m)} v_\nu) = \Ytil^{l,m} b_\nu.
$$

\begin{lem}\label{symmetry of deg of Y(m,n) vnu}
Let $l \in [0,\nu_1-\nu_3]$, $m \in [0,\nu_1+\nu_3]$ be such that $l+m \leq \nu_1$. Then, we have
$$
\deg(Y^{(l,m)}v_\nu) = \deg(Y^{(\nu_1-\nu_3-l,\nu_1+\nu_3-m)}v_\nu).
$$
\end{lem}

\begin{proof}
By equation \eqref{deg on Ymn vnu}, we have
\begin{align}
\begin{split}
&\deg(Y^{(\nu_1-\nu_3-l,\nu_1+\nu_3-m)}v_\nu) \\
&= \hf((\nu_1-\nu_3-l)l + (\nu_1+\nu_3-m)m + \nu_1- |\nu_3| - |\nu_1-l-m| + |\nu_3+l-m| \\
&\qu+ \sum_{i=1}^{\nu_1-\nu_3-l}(|\nu_1-i+1|-|\nu_1-l-i+1|) + \sum_{j=1}^{\nu_1+\nu_3-m}(|\nu_1-j+1| - |\nu_3-j|)).
\end{split} \nonumber
\end{align}
Here, we have
\begin{align}
\begin{split}
&\sum_{i=1}^{\nu_1-\nu_3-l}(|\nu_1-i+1|-|\nu_1-l-i+1|) \\
&= \begin{cases}
\sum_{s=1}^{\nu_1-\nu_3-l} |\nu_1-s+1| - \sum_{s=1}^{\nu_1-\nu_3-l} |\nu_1-m-t+1| \qu & \IF \nu_1-l \geq \nu_3+l, \\
(\sum_{s=1}^{l} |\nu_1-s+1| - \sum_{t=1}^{-\nu_1+\nu_3+2l} |\nu_3+l-t+1|) - \\ (\sum_{s=1}^{l} |\nu_3+l-s+1| - \sum_{t=1}^{-\nu_1+\nu_3+2l} |\nu_3+l-t+1|) \qu & \IF \nu_1-l \leq \nu_3+l
\end{cases} \\
&= \sum_{s=1}^{l} |\nu_1-s+1| - \sum_{s=1}^{l} |\nu_3+l-s+1| \\
&= \sum_{i=1}^l (|\nu_1-i+1| - |\nu_3+l-i+1|).
\end{split} \nonumber
\end{align}
Similarly, we obtain
$$
\sum_{j=1}^{\nu_1+\nu_3-m}(|\nu_1-j+1| - |\nu_3-j|) = \sum_{j=1}^{m}(|\nu_1-j+1| - |\nu_3-j|).
$$
Thus, the assertion follows.
\end{proof}

Set $d_\nu^{l,m} := \deg(Y^{(l,m)}v_\nu)$. By direct calculation, we obtain the following.

\begin{lem}\label{relative deg on Ymn}
Let $l \in [0,\nu_1-\nu_3]$, $m \in [0,\nu_1+\nu_3]$ be such that $l+m \leq \nu_1$.
\begin{enumerate}
\item\label{relative deg on Ymn 1} If $\nu_3+l-m \geq 0$, then we have
$$
d_\nu^{l,m} - d_\nu^{l+1,m-1} = \begin{cases}
2(\nu_3+l-m)+1 > 0 \qu & \IF \nu_3-m \leq 0, \\
(\nu_3+l-m)+l+1 > 0 \qu & \IF \nu_3-m \geq 0.
\end{cases}
$$
\item\label{relative deg on Ymn 2} If $\nu_3+l-m \leq 0$, then we have
$$
d_\nu^{l,m} - d_\nu^{l-1,m+1} = \begin{cases}
-2(\nu_3+l-m)+1 > 0 \qu & \IF \nu_3+l \geq 0, \\
-(\nu_3+l-m)+m+1 > 0 \qu & \IF \nu_3+l \leq 0.
\end{cases}
$$
\item\label{relative deg on Ymn 3} If $\nu_3+l-m = -1$, then we have
$$
d_\nu^{l,m} = d_\nu^{l+1,m-1}.
$$
\end{enumerate}
\end{lem}

\begin{prop}\label{deg3 leq 1}
Let $k \in [0,2\nu_1]$.
\begin{enumerate}
\item Suppose that $\nu_3 \geq 0$ and set
$$
(l_0,m_0) := \begin{cases}
(0,k) \qu & \IF 0 \leq k \leq \nu_3, \\
(\frac{k-\nu_3}{2},\frac{k+\nu_3}{2}) \qu & \IF \nu_3 \leq k \leq 2\nu_1-\nu_3 \AND k-\nu_3 \in \Z_{\ev}, \\
(\frac{k-\nu_3+1}{2},\frac{k+\nu_3-1}{2}) \qu & \IF \nu_3 \leq k \leq 2\nu_1-\nu_3 \AND k -\nu_3 \in \Z_{\odd}, \\
(\nu_1-\nu_3, k-\nu_1+\nu_3) \qu & \IF 2\nu_1-\nu_3 \leq k \leq 2\nu_1.
\end{cases}
$$
Then, we have
$$
\Ytil^k b_\nu = \begin{cases}
\frac{1}{\sqrt{2}}(\Ytil^{l_0,m_0} + \Ytil^{l_0-1,m_0+1})b_\nu \qu & \IF \nu_3 \leq k \leq 2\nu_1-\nu_3 \AND k -\nu_3 \in \Z_{\odd}, \\
\Ytil^{l_0,m_0} b_\nu \qu & \OW.
\end{cases}
$$
\item Suppose that $\nu_3 \leq 0$ and set
$$
(l_0,m_0) := \begin{cases}
(k,0) \qu & \IF 0 \leq k \leq -\nu_3, \\
(\frac{k-\nu_3}{2},\frac{k+\nu_3}{2}) \qu & \IF -\nu_3 \leq k \leq 2\nu_1+\nu_3 \AND k+\nu_3 \in \Z_{\ev}, \\
(\frac{k-\nu_3+1}{2},\frac{k+\nu_3-1}{2}) \qu & \IF -\nu_3 \leq k \leq 2\nu_1+\nu_3 \AND k+\nu_3 \in \Z_{\odd}, \\
(k-\nu_1-\nu_3,\nu_1+\nu_3) \qu & \IF 2\nu_1+\nu_3 \leq k \leq 2\nu_1.
\end{cases}
$$
Then, we have
$$
\Ytil^k b_\nu = \begin{cases}
\frac{1}{\sqrt{2}}(\Ytil^{l_0,m_0} + \Ytil^{l_0-1,m_0+1})b_\nu \qu & \IF -\nu_3 \leq k \leq 2\nu_1+\nu_3 \AND k +\nu_3 \in \Z_{\odd}, \\
\Ytil^{l_0,m_0} b_\nu \qu & \OW.
\end{cases}
$$
\end{enumerate}
\end{prop}

\begin{proof}
The assertion follows from the definition of $\Ytil^{l,m}$ and Lemmas \ref{symmetry of deg of Y(m,n) vnu} and \ref{relative deg on Ymn}.
\end{proof}

\begin{lem}\label{existence of based structure on highest component}
Let $\nu \in X_{\frk,\Int}^+$, $M$ a based classical weight standard $X^\imath$-weight module satisfying condition \eqref{Property3}. Suppose that $M[> \nu] = 0$. Then, $M$ contains a $\Ui$-submodule isomorphic to $V(\nu)$ which is also a based submodule of $M$.
\end{lem}

\begin{proof}
We prove the statement under the assumption that $\nu_3 \geq 0$. The $\nu_3 \leq 0$ case can be proved similarly.

Let $b \in \ol{\clL}_{M,\nu}$ and set $v := G^\imath(b)$. By our hypothesis, $v$ is a highest weight vector of weight $\nu$. Hence, in order to prove the assertion, it suffices to show that $G^\imath(\Ytil^{l,m} b) \in \Ui v$ for all $l \in [0,\nu_1-\nu_3]$ and $m \in [0,\nu_1+\nu_3]$. Below, we concentrate on the case when $l+m \leq \nu_1$. The case when $l+m \geq \nu_1$ can be proved similarly by exchanging the roles of $v$ and $Y^{(2\nu_1)} v (= Y^{(\nu_1-\nu_3,\nu_1+\nu_3)} v)$.

Set
$$
B := \{ \Ytil^{l,m} v \mid l \in [0,\nu_1-\nu_3], m \in [0,\nu_1+\nu_3] \}.
$$
Then, it forms a basis of $\Ui v$ consisting of vectors in $\clL_M$. Let $\Btil \subseteq \clL_M$ be a basis of $M$ which extends $B$.

We show that $G^\imath(\Ytil^{l,m} b) \in \Ui v$ inductively. Assume that we have $G^\imath(\Ytil^{l',m'} b) \in \Ui v$ for the following cases:
\begin{itemize}
\item $|\nu_1-(l'+m')| > \nu_1-(l+m)$, or
\item $l'+m' = l+m$ and $|\nu_3+l'-m'| < |\nu_3+l-m|$.
\end{itemize}
For such $l',m'$, we can write
$$
G^\imath(\Ytil^{l',m'} b) = \sum_{\substack{l'' \in [0,\nu_1-\nu_3] \\ m'' \in [0,\nu_1+\nu_3]}} p^{l'',m''}_{l',m'} \Ytil^{l'',m''} v
$$
for some $p^{l'',m''}_{l',m'} \in \bfK_\infty$. By condition \eqref{Property3}, we see that $p^{l'',m''}_{l',m'} = 0$ if
\begin{itemize}
\item $|\nu_1-(l''+m'')| < |\nu_1-(l'+m')|$, or
\item $|\nu_3+l''-m''| < |\nu_3+l'-m'|$.
\end{itemize}
Assume further that for each $l'',m''$ such that $l''+m'' = l'+m'$ and $|\nu_3+l''-m''| \geq |\nu_3+l'-m'|$, we have
$$
\deg(p^{l'',m''}_{l',m'}) \leq d_\nu^{l'',m''} - d_\nu^{l'_0,m'_0} + \sum_{j=1}^{\tilde{k}_{l',m'}} (|\nu_3+l''-m''|-j),
$$
where $(l'_0,m'_0)$ is as in Lemma \ref{deg3 leq 1} with $l$ being replaced with $l'+m'$, and
$$
\tilde{k}_{l',m'} := \begin{cases}
|l'-l'_0| \qu & \IF l'+m' \leq \nu_3 \OR l'+m' \geq 2\nu_1-\nu_3, \\
|\nu_3+l'-m'| \qu & \IF \nu_3 \leq l'+m' \leq 2\nu_1-\nu_3 \AND |\nu_3+m'-n'| \leq l'+m'-\nu_3, \\
l' \qu & \IF \nu_3 \leq l'+m' \leq 2\nu_1-\nu_3 \AND |\nu_3+l'-m'| \geq l'+m'-\nu_3.
\end{cases}
$$

Under the induction hypothesis above, we show that $G^\imath(\Ytil^{l,m} b) \in \Ui v$. Before doing so, let us prepare some notation. Set $k := l+m$. Then, we can write as
$$
\Ytil^k v = \sum_{l'+m' = k} r^{l',m'} \Ytil^{l',m'} v
$$
for some $r^{l',m'} \in \bfK_\infty$ such that
$$
\lt(r^{l',m'}) = q^{d_\nu^{l',m'} - d_\nu^{l_0,m_0}}.
$$
Also, for each $j \geq 0$, we have
$$
B_3^{(j)} \Ytil^k v = \sum_{l'+m' = k} r^{l',m',j} \Ytil^{l',m'} v,
$$
where $r^{l',m',j} := r^{l',m'} \frac{[\nu_3+l'-m']^j}{[j]!}$. Hence, we have
$$
\lt(r^{l',m',j}) = q^{d_\nu^{l',m'} - d_\nu^{l_0,m_0} + \sum_{i=1}^j (|\nu_3+l'-m'|-i)}.
$$
Note that we have
\begin{align}\label{deg of rmnk}
\deg(r^{l',m',0}) < \deg(r^{l',m',1}) < \cdots < \deg(r^{l',m',|\nu_3+l'-m'|-1}). 
\end{align}

Let us prove that $G^\imath(\Ytil^{l,m} b) \in \Ui v$. First, suppose that $k \leq \nu_3$. In this case, we have
$$
(l_0,m_0) = (0,k), \AND \tilde{k}_{l,m} = l.
$$
Using Lemma \ref{relative deg on Ymn}, we see that for each $l',m'$ such that $l'+m' = k$ and for $j \geq 0$, we have
\begin{align}
\begin{split}
\deg(r^{l',m',j}) &= -\hf l'(2(\nu_3-m')+l'-1) + \hf j(2(\nu_3-m')+2l'-j-1) \\
&= -\hf(l'-j)(2(\nu_3-m')+l'-j-1).
\end{split} \nonumber
\end{align}
This implies that
\begin{align}
\begin{split}\label{degree estimate}
&\deg(r^{l',m',l'}) = 0, \\
&\deg(r^{l',m',j}) < 0 \qu \IF j < l'.
\end{split}
\end{align}

By Theorem \ref{Vnu is a based Ui-module for n=3}, we have
$$
G^\imath(\Ytil^{0,k}b) = G^\imath(\Ytil^k b) \in \Ui v.
$$
By the proof of Proposition \ref{information about Vnu}, we see that
$$
G^\imath(\Ytil^k b) \in \Ytil^k v + \bigoplus_{\substack{j < k \OR j > 2\nu_1-k \\ \ol{j} = \ol{k}}} q\inv \bfK_{\infty} \Ytil^j v.
$$
This shows that $G^\imath(\Ytil^k b)$ is a sum of $B_3$-eigenvectors of eigenvalues $[a]$, $a \in \Z_{p(\nu_3-k)}$. Since $M$ is a standard $X^\imath$-weight module, we see that
$$
\mathbf{1}_\zeta G^\imath(\Ytil^k b) = 0
$$
for all $\zeta \in X^\imath$ such that $\zeta_3 \neq \ol{\nu_3-k}$. Therefore, we obtain
$$
\sum_{\substack{\zeta \in X^\imath \\ \zeta_3 = \ol{\nu_3-k}}} \mathbf{1}_\zeta G^\imath(\Ytil^k b) = \sum_{\zeta \in X^\imath} \mathbf{1}_\zeta G^\imath(\Ytil^k b) = G^\imath(\Ytil^k b).
$$
Hence, we have
$$
B_{3,p(\nu_3-k)}^{(l)} G^\imath(\Ytil^k b) = \sum_{\zeta_3 = \ol{\nu_3-k}} B_{3,\zeta}^{(l)} G^\imath(\Ytil^k b) \in \UidotA v.
$$
By the definition of $\imath$divided powers, we can write as
\begin{align}
\begin{split}
B_{3,p(\nu_3-k)}^{(l)} G^\imath(\Ytil^{k} b) &= \sum_{j=0}^l c_j B_3^{(j)} G^\imath(\Ytil^{k} b) \\
&\in \sum_{j=0}^l c_j B_3^{(j)} \Ytil^{k} v + \sum_{\substack{l'+m' < k \OR \\ l'+m' > 2\nu_1-k}} \bfK \Ytil^{l',m'} v
\end{split} \nonumber
\end{align}
for some $c_j \in \delta_{j,l} + q\inv \bfK_\infty$. For each $j \in [0,l]$, set
$$
u'_j := c_j B_3^{(j)} \Ytil^k v = \sum_{l' + m' = k} c_j r^{l',m',j} \Ytil^{l',m'} v,
$$
and
$$
l' := \min \{ l'' \neq l \mid \deg(c_j r^{l'',m'',j}) \geq 0 \}.
$$
By \eqref{degree estimate}, we must have $l' < j$. Set $m' := k-l'$. Then, there exists $a_{l',j} \in \bfA_{\Inv}$ such that $u'_j-a_{l',j} G^\imath(\Ytil^{l',m'} b)$, expanded by the basis $B$, contains $\Ytil^{l',m'} v$ with coefficient in $q\inv \bfK_\infty$. Note that we have
$$
\deg(a_{l',j}) = \deg(c_j r^{l',m',j}).
$$
Then, for each $l'' > l'$, setting $m'' := k-l''$, we compute as
\begin{align}
\begin{split}
&\deg(c_j r^{l'',m'',j}) - \deg(a_{l',j} p^{l'',m''}_{l',m'}) \\
&\qu = \deg(c_j) + \deg(r^{l'',m'',j}) - (\deg(c_j) + \deg(r^{l',m',j}) + \deg(p^{l'',m''}_{l',m'})) \\
&\qu = \deg(r^{l'',m'',j}) - \deg(r^{l',m',j}) - \deg(p^{l'',m''}_{l',m'}) \\
&\qu \geq \deg(r^{l'',m'',j}) - \deg(r^{l',m',j}) - \deg(r^{l'',m'',l'}) \\
&\qu = d_\nu^{l'',m''} - d_\nu^{l_0,m_0} + j(\nu_3+l''-m''-1)-\hf j(j+1) \\
&\qu \qu -(d_\nu^{l',m'} - d_\nu^{l_0,m_0} + j(\nu_3+l'-m'-1)-\hf j(j+1)) \\
&\qu \qu -(d_\nu^{l'',m''} - d_\nu^{l_0,m_0} + l'(\nu_3+l''-m''-1)-\hf l'(l'+1)) \\
&\qu = (j-l')(\nu_3+l''-m''-1) - \deg(r^{l',m',l'}) -(j-l')(\nu_3+l'-m'-1) \\
&\qu =(j-l')((l''-m'')-(l'-m')) > 0.
\end{split} \nonumber
\end{align}
Therefore, both $u'$ and $u'-a_{l',j} G^\imath(\Ytil^{l',m'} b)$ contain $\Ytil^{l'',m''}v$, $l'' > l'$, $l''+m'' = k$ with coefficient in $q^{\deg(c_j r^{l'',m'',j})} \bfK_\infty$.

Replacing $u'$ with $u'-a_{l',j} G^\imath(\Ytil^{l',m'} b)$, and repeating the procedure above, we finally obtain that
$$
u_j := c_j B_3^{(j)} \Ytil^k v - \sum_{l' = 0}^{j-1} a_{l',j} G^\imath(\Ytil^{l',m'} b)
$$
contains $\Ytil^{l'',m''} v$, $l'' < j$, $l''+m'' = k$ with coefficient in $q\inv \bfK_\infty$, and $\Ytil^{l'',m''}v$, $l'' \geq j$, $l''+m'' = k$ with coefficient in $q^{\deg(c_j r^{l'',m'',j})} \bfK_\infty$, where $a_{l',j} \in \bfA_{\Inv}$.

By above, we see that
$$
u := B_{3,p(\nu_3+l-m)}^{(l)} G^\imath(\Ytil^{k} b) - \sum_{j=0}^l \sum_{l' = 0}^{j-1} a_{l',j} G^\imath(\Ytil^{l',k-l'} b) \in \sum_{j=0}^l u_j + \sum_{\substack{l'+m' < k \OR\\ l'+m' > 2\nu_1-k}} \bfK \Ytil^{l',m'} v
$$
contains $\Ytil^{l'',m''}v$, $l'' < l$, $l''+m'' = k$ with coefficient in $q\inv \bfK_\infty$, and $\Ytil^{l'',m''}v$, $l'' \geq l$, $l''+m'' = k$ with the degree of coefficient being
$$
\max_{0 \leq j \leq l}\{ \deg(c_j r^{l'',m'',j}) \}.
$$
Since $m'' \leq k \leq \nu_3$, we have $\nu_3 + l''-m''  \geq l'' \geq l$. Then, by equation \eqref{deg of rmnk}, we see that
$$
\max_{0 \leq j \leq l}\{ \deg(c_j r^{l'',m'',j}) \} = \deg(r^{l'',m'',l}) \leq 0;
$$
the equality holds if and only if $l'' = l$ by equation \eqref{degree estimate}.

By the construction above, $u - G^\imath(\Ytil^{l,m} b)$ is contained in $M_{\bfA}$, and is invariant under $\psii_M$. Furthermore, the vector $u$, expanded by $\Btil$, contains $w \in \Btil$ with coefficient in $q\inv \bfK_\infty$ unless $w = \Ytil^{l',m'} v$ with $l'+m' < k$ or $l'+m' > 2\nu_1-k$. From these observations and condition \eqref{Property3}, we can find $a'_{l',m'} \in \bfA_{\Inv}$ such that
$$
u - G^\imath(\Ytil^{l,m} b) - \sum_{\substack{l'+m' < k \OR \\ l'+m' > 2\nu_1-k}} a'_{l',m'} G^\imath(\Ytil^{l',m'} b) \in q\inv \clL_M \cap M_{\bfA}.
$$
Then, by Lemma \ref{characterization of icanonical basis element}, we see that
$$
u - G^\imath(\Ytil^{l,m} b) - \sum_{\substack{l'+m' < k \OR \\ l'+m' > 2\nu_1-k}} a'_{l',m'} G^\imath(\Ytil^{l',m'} b) = 0.
$$
This shows that
$$
G^\imath(\Ytil^{l,m} b) \in \Ui v,
$$
as desired. By the construction, the degree of $\Ytil^{l',m'} v$, $l' > l$, $l'+m'=k$ in $G^\imath(\Ytil^{l,m} b)$ is at most $\deg(r^{l'',m'',l})$. Hence, we can proceed our induction.

Next, suppose that $\nu_3 \leq k \leq \nu_1$ and $k-\nu_3 \in \Z_{\ev}$ (the case when $k-\nu_3 \in \Z_{\odd}$ is similar). By the same argument as above, we obtain
$$
G^\imath(\Ytil^{l_0,m_0} b) = G^\imath(\Ytil^k b) \in \Ui v.
$$
Computing $B_{3,\ev}^{(\tilde{k}_{m,n})} G^\imath(\Ytil^k b_\nu)$, one can construct $G^\imath((\Ytil^{l_+,m_+} \pm \Ytil^{l_-,m_-}) b)$ if $|\nu_3+l-m| \leq k-\nu_3$, and $G^\imath(\Ytil^{l,m} b)$ if $\nu_3+l-m > k-\nu_3$ in the same way as above, where $(l_\pm,m_\pm)$ is such that $l_\pm + m_\pm = k$ and
$$
\nu_3+l_+-m_+ = -(\nu_3+l_--m_-) = |\nu_3+l-m|
$$
This completes the proof.
\end{proof}

\begin{theo}\label{Vnu is a based Ui-module for n=4}
Let $\nu \in X_{\frk,\Int}^+$. Then, $V(\nu)$ is a based $\Ui$-module.
\end{theo}

\begin{proof}
Let $M$ be a finite-dimensional based classical weight standard $X^\imath$-weight module satisfying condition \eqref{Property3} such that $M[\nu] \neq 0$. Such $M$ surely exists by Propositions \ref{Tensor power contains all highest weight vectors} and \ref{canonical Xi-weight structure is standard}. Let us write an irreducible decomposition of $M$ as
$$
M = \bigoplus_{k=1}^r M_k, \qu M_k \simeq V(\nu_k), \qu \nu_k \in X_{\frk,\Int}^+.
$$
Without loss of generality, we may assume that $k \leq l$ implies $\nu_k \leq \nu_l$. By Lemma \ref{existence of based structure on highest component}, we see that $M_r$ is a based $\Ui$-submodule of $M$, and hence, $V(\nu_r)$ is a based $\Ui$-module. Then, $M/M_r$ is a finite-dimensional based classical weight standard $X^\imath$-weight module satisfying condition \eqref{Property3} by Propositions \ref{sub and quot of Xi-weight module} \eqref{sub and quot of Xi-weight module 2} and \ref{based submodule and based quotient module} \eqref{based submodule and based quotient module 2}. Replacing $M$ with $M/M_r$, we see that $V(\nu_{r-1})$ is a based $\Ui$-module. Proceeding in this way, we conclude that $V(\nu_k)$ is a based $\Ui$-module for all $k \in [1,r]$. This completes the proof.
\end{proof}

\begin{prop}\label{preparation for branching rule n=4}
Let $M$ be a finite-dimensional classical weight $\Ui$-module equipped with a contragredient Hermitian inner product. Set
\begin{align}
\begin{split}
&L_1 := \{ b \in \ol{\clL}_M \mid \Xtil_+ b = \Xtil_- b = 0 \}, \\
&L_2 := \C \{ b \in \ol{\clL}_M \mid \text{$b$ is $B_1$-homogeneous,}\ \Btil_2 b = 0, \AND \deg_3((\Btil_2 \Btil_1)^{\deg_3(b)} b) = 0 \}.
\end{split} \nonumber
\end{align}
Then, the linear map
$$
L_2 \rightarrow L_1;\ b \mapsto (1 + \Btil_1)b
$$
is an isomorphism of $\C$-vector spaces with inverse
$$
L_1 \rightarrow L_2;\ b \mapsto \sum_{\substack{\zeta \in X^\imath \\ \zeta_2 = \ol{0}}} \mathbf{1}_\zeta b.
$$
\end{prop}

\begin{proof}
As the proof of Proposition \ref{preparation for branching rule n=3}, it suffices to consider the case when $M = V(\nu)$ for some $\nu \in X_{\frk,\Int}^+$. Clearly, we have
$$
L_1 = \C b_\nu.
$$

Let $b \in \ol{\clL}(\nu)$ be $\Btil_1$-homogeneous such that $\Btil_2 b = 0$ and $\deg_3((\Btil_2 \Btil_1)^{\deg_3(b)} b) = 0$. Set $l_1 := \deg_1(b)$ and $l_3 := \deg_3(b)$. By Proposition \ref{preparation for branching rule n=3}, there exists $b' \in \ol{\clL}(\nu)$ such that $b'$ is a sum of $X_{\frk}$-weight vectors of weights $\xi \in X_{\frk,\Int}$ with $\xi_1 = l_1$, and such that $\Xtil b' = 0$ and $b = (1 + (-1)^{l_1} \Ytil^{2l_1}) b'$. Then, from Corollary \ref{Btil on Vnu for rank 2}, we see that
$$
(\Btil_2 \Btil_1)^{l_3} b \in \C (\Ytil^{l_3} + (-1)^{l_1+l_3} \Ytil^{2l_1-l_3})b'.
$$
On the other hand, by Proposition \ref{deg3 leq 1}, we see that for each $b'' \in \ol{\clL}(\nu)$, we have $\deg_3(b'') = 0$ only if $b'' \in \bigoplus_{n=0}^{2\nu_1}\C \Ytil^n b_\nu$. Therefore, we obtain $l = \nu_1$, $b' \in \C b_\nu$, and hence,
$$
L_2 \subseteq \C (1 + (-1)^{\nu_1} \Ytil^{2\nu_1}) b_\nu.
$$
Conversely, by Proposition \ref{deg3 leq 1} again, we have
$$
\deg_3((\Btil_2 \Btil_1)^{|\nu_3|}(1 + (-1)^{\nu_1} \Ytil^{2\nu_1}) b_\nu) = \deg_3((\Ytil^{|\nu_3|} + (-1)^{\nu_1+|\nu_3|} \Ytil^{2\nu_1-|\nu_3|}) b_\nu) = 0.
$$
Thus, we see that
$$
L_2 = \C (1 + (-1)^{\nu_1} \Ytil^{2\nu_1}) b_\nu.
$$

Now, by the same calculation as in the proof of Proposition \ref{preparation for branching rule n=3}, the assertion follows.
\end{proof}

\section{General case}\label{section: general}
In this section, we consider the general $n \geq 3$ case, and set
$$
m := \begin{cases}
  \frac{n}{2} \qu & \IF n \in \mathbb{Z}_{\ev}, \\
  \frac{n-1}{2} \qu & \IF n \in \mathbb{Z}_{\odd}.
\end{cases}
$$

Let $j \in \Itil_{\frk}$, and set
$$
J = \begin{cases}
\{ 2i-1,2i,2i+1 \} \qu & \IF j = (2i,\pm), \\
\{ 2m-1,2m \} \qu & \IF j = 2m.
\end{cases}
$$
Then, $\Ui_J$ is canonically isomorphic to the $\imath$quantum group of type AI with $n = 2$ or $3$. Hence, we can define linear operators $\Xtil_j,\Ytil_j$ to correspond $\Xtil,\Ytil$ if $j = 2m$ or $\Xtil_\pm,\Ytil_\pm$ if $j = (2i,\pm)$.

\subsection{Based module structure of irreducible modules}
In this subsection, we prove that $V(\nu)$, $\nu \in X_{\frk,\Int}^+$ is a based $\Ui$-module.

\begin{lem}\label{induction step for constructing icanonical basis}
Let $M$ be a finite-dimensional based classical weight standard $X^\imath$-weight module satisfying condition \eqref{Property3}. Let $N \subseteq M$ be a $\Uidot$-submodule, $\nu \in X_{\frk,\Int}$, and $j \in \Itil_{\frk}$. Suppose that for each $\xi \geq \nu$ and $b \in \ol{\clL}_{N,\xi}$ we have $G^\imath(b) \in N$. Then, for each $b \in \ol{\clL}_{N,\nu}$, we have $G^\imath(\Ytil_j b) \in N$.
\end{lem}

\begin{proof}
Define $J \subseteq I$ as above. Let us view $M$ and $N$ as $\Ui_J$-modules. Then, we have an orthogonal decomposition
$$
\ol{\clL}_N = \bigoplus_{\xi \in X_{\frk,\Int}^+(J)} \ol{\clL}_N[\xi].
$$
Let $\{ b_1,\ldots,b_r \}$ be an orthonormal basis of $\bigoplus_{\xi \geq \nu} (\ol{\clL}_{N}[\xi] \cap \ol{\clL}_{N,\xi})$ consisting of $X_{\frk}$-weight vectors. Let $\xi_k \in X_{\frk,\Int}^+(J)$ be such that $b_k \in \ol{\clL}_{N,\xi_k}$. We may assume that $k \leq l$ implies $\xi_k \leq \xi_l$. Note that we have
$$
\ol{\clL}_{N,\nu} \subseteq \Span_{\C} \{ \Ytil_{j_1} \cdots \Ytil_{j_s} b_k \mid s \geq 0,\ k \in [1,r],\ j_1,\ldots,j_s \in \Itil_{\frk}(J) \}.
$$
Hence, to prove the assertion, it suffices to show that $G^\imath(\Ytil_{j_1} \cdots \Ytil_{j_s} b_k) \in N$ for all $s,k,j_1,\ldots,j_s$.

By weight consideration, $G^\imath(b_r) \in N$ is a highest weight vector of weight $\xi_r$. Then, by Theorems \ref{Vnu is a based Ui-module for n=3} and \ref{Vnu is a based Ui-module for n=4}, $\Ui_J G^\imath(b_r)$ is a based $\Ui_J$-submodule of $M$, and $G^\imath(\Ytil_{j_1} \cdots \Ytil_{j_s} b_r) \in \Ui_J G^\imath(b_r) \subseteq N$. Replacing $M$ with $M/\Ui_J G^\imath(b_r)$ and repeating the same argument as above, we see that $G^\imath(\Ytil_{j_1} \cdots \Ytil_{j_s} b_{r-1}) \in N$ for all $j_1,\ldots,j_s \in J$. Proceeding in this way, we conclude that $G^\imath(\Ytil_{j_1} \cdots \Ytil_{j_s} b_k) \in N$ for all $s \geq 0$, $k \in [1,r]$, $j_1,\ldots,j_s \in \Itil_{\frk}(J)$, as desired.
\end{proof}

\begin{theo}\label{V(nu) is a based Ui-module; general}
Let $\nu \in X_{\frk,\Int}^+$. Then, $V(\nu)$ is a based $\Ui$-module. Moreover, we have the following:
\begin{enumerate}
\item\label{V(nu) is a based Ui-module; general 1} $\clL(\nu) = \Span_{\bfK_\infty}\{ \Ytil_{j_1} \cdots \Ytil_{j_r} v_\nu \mid r \geq 0,\ j_1,\ldots,j_r \in \Itil_{\frk} \}$.
\item\label{V(nu) is a based Ui-module; general 2} $\ol{\clL}(\nu) = \Span_{\C}\{ \Ytil_{j_1} \cdots \Ytil_{j_r} b_\nu \mid r \geq 0,\ j_1,\ldots,j_r \in \Itil_{\frk} \}$.
\item\label{V(nu) is a based Ui-module; general 3} For each $b \in \ol{\clL}(\nu)$, if $\Xtil_j b = 0$ for all $j \in \Itil_{\frk}$, then we have $b \in \C b_\nu$.
\end{enumerate}
\end{theo}

\begin{proof}
Let $M$ be a finite-dimensional based classical weight standard $X^\imath$-weight module satisfying condition \eqref{Property3} such that $M[\nu] \neq 0$. We proceed by descending induction on $\{ \nu' \in X_{\frk,\Int}^+ \mid M[\nu'] \neq 0 \}$ that $V(\nu')$ is a based $\Ui$-module. Assuming that our claim is true for all $\nu' > \nu$ and replacing $M$ with $M/M[> \nu]$, we may assume that $M[> \nu] = 0$. Let $b_0 \in \ol{\clL}_{M,\nu}$. Then, $G^\imath(b_0)$ is a highest weight vector of weight $\nu$. Hence, we can identify $V(\nu)$ with an $X^\imath$-weight submodule $\Ui G^\imath(b_0)$ of $M$. Under this identification, we have $b_0 = b_\nu$.

Let $\xi \in X_{\frk,\Int}$ be such that $\xi < \nu$. Assume that for all $\xi' > \xi$ and $b \in \ol{\clL}(\nu)_{\xi'}$, we have $G^\imath(b) \in V(\nu)$. Then, $V(\nu)_{\xi'}$ is spanned by $G^\imath_0(\ol{\clL}(\nu)_{\xi'})$. Actually, for a basis $\clB$ of $\ol{\clL}(\nu)_{\xi'}$, $G^\imath_0(\clB)$ forms a basis of $V(\nu)_{\xi'}$.

Let $v \in \clL(\nu)_\xi$, and set $b := \ev_\infty(v) \in \ol{\clL}(\nu)_\xi$. Since we have
$$
V(\nu)_\xi = \sum_{j \in \Itil_{\frk}} Y_j V(\nu)_{\xi + \gamma_j} = \sum_{j \in \Itil_{\frk}} \Ytil_j V(\nu)_{\xi + \gamma_j}
$$
by Theorem \ref{facts about classical weight Ui-modules} \eqref{facts about classical weight Ui-modules 7} and Propositions \ref{Im Yj = Im Ytilj for rank 1 AI} and \ref{Im Yj = Im Ytilj for rank 2 AI}, we can write as
$$
v = \sum_{j \in \Itil_{\frk}} v_j, \qu v_j \in \Ytil_j V(\nu)_{\xi + \gamma_j}.
$$
We show that we can take $v_j \in \clL(\nu)$ for all $j \in \Itil_{\frk}$.

For each $j \in \Itil_{\frk}$, choose an orhtonormal basis $\clB_j$ of $\Ytil_j \ol{\clL}(\nu)_{\xi + \gamma_j}$. By Lemma \ref{induction step for constructing icanonical basis}, $G^\imath(\clB_j) \subseteq V(\nu)$. Consequently, $G^\imath_0(\clB_j)$ forms an almost orthonormal basis of $\Ytil_j V(\nu)_{\xi + \gamma_j}$. Let us write
$$
v_j = \sum_{b \in \clB_j} c_b G^\imath_0(b)
$$
with $c_b \in \bfK$. Set $d := \max\{ \deg(c_b) \mid b \in \bigsqcup_j \clB_j \}$. Assume that $d > 0$. Then, we have
$$
v_j \in \sum_{\substack{b \in \clB_j \\ \deg(c_b) = d}} \lt(c_b) G^\imath_0(b) + q^{d-1} \clL_M.
$$
This shows that
$$
\sum_{\substack{b \in \bigsqcup_j \clB_j \\ \deg(c_b) = d}} \lc(c_b) G^\imath_0(b) \equiv q^{-d} v \equiv 0
$$
modulo $q\inv \clL_M$. Taking $\ev_\infty$, we obtain
$$
\sum_{\substack{b \in \bigsqcup_j \clB_j \\ \deg(c_b) = d}} \lc(c_b) b = 0,
$$
and hence,
$$
\sum_{\substack{b \in \bigsqcup_j \clB_j \\ \deg(c_b) = d}} \lt(c_b) G^\imath_0(b) = q^d G^\imath_0(\sum_{\substack{b \in \bigsqcup_j \clB_j \\ \deg(c_b) = d}} \lc(c_b) b) = 0.
$$

Replacing $v_j$ with
$$
v_j - \sum_{\substack{b \in \clB_j \\ \deg(c_b) = d}} \lt(c_b) G^\imath_0(b) \in \Ytil_j V(\nu)_{\xi + \gamma_j},
$$
we still have $v = \sum_j v_j$. Furthermore, if we write $v_j = \sum_{b \in \clB_j} c_b G^\imath_0(b)$ again, then, we have
$$
\max\{ \deg(c_b) \mid b \in \bigsqcup_j \clB_j \} < d.
$$
Repeating this procedure, we conclude that $v_j \in \clL(\nu)$ for all $j$, as desired.

Now, we have $v = \sum_j v_j$ with $v_j \in \Ytil_j \clL(\nu)_{\xi + \gamma_j}$. Taking $\ev_\infty$, we see that $b \in \sum_j \Ytil_j \ol{\clL}(\nu)_{\xi + \gamma_j}$. This shows that
$$
\clL(\nu)_\xi = \sum_{j \in \Itil_{\frk}} \Ytil_j \clL(\nu)_{\xi+\gamma_j}, \qu \ol{\clL}(\nu)_\xi = \sum_{j \in \Itil_{\frk}} \Ytil_j \ol{\clL}(\nu)_{\xi+\gamma_j}.
$$
Then, by Lemma \ref{induction step for constructing icanonical basis}, we obtain
$$
G^\imath(\ol{\clL}(\nu)_\xi) \subseteq V(\nu)
$$
Now, the assertions, except \eqref{V(nu) is a based Ui-module; general 3}, follow.

Let us prove \eqref{V(nu) is a based Ui-module; general 3}. Let $b \in \ol{\clL}(\nu)$ be such that $\Xtil_j b = 0$ for all $j \in \Itil_{\frk}$. Then, by Propositions \ref{Ftil and Etil are almost adjoint AI} and \ref{Ftil and Etil are almost adjoint AI 2}, for each $b' \in \ol{\clL}(\nu)$ we have
$$
(b, \Ytil_j b')_\nu = (\Xtil_j b, b')_\nu = 0.
$$
This shows that
$$
(b, \bigoplus_{\xi < \nu} \ol{\clL}(\nu)_\xi)_\nu = 0.
$$
By weight consideration, we obtain $b \in \ol{\clL}(\nu)_\nu = \C b_\nu$, as desired.
\end{proof}

\begin{cor}
Let $M$ be a finite-dimensional based classical weight standard $X^\imath$-weight module satisfying condition \eqref{Property3}. Then, for each $\nu \in X_{\frk,\Int}^+$, the following hold.
\begin{enumerate}
\item $M[\geq \nu]$ is a based $\Ui$-submodule of $M$.
\item $M[> \nu]$ is a based $\Ui$-submodule of $M$ and $M[\geq \nu]$.
\item $M[\geq \nu]/M[> \nu]$ is a based $\Ui$-module isomorphic to $V(\nu)^{m_\nu}$, where $m_\nu$ denotes the multiplicity $\dim_{\bfK} \Hom_{\Ui}(V(\nu), M)$ of $V(\nu)$ in $M$.
\end{enumerate}
\end{cor}

\begin{proof}
The assertion follows from Lemma \ref{cellularity of highest component} and Theorem \ref{V(nu) is a based Ui-module; general}.
\end{proof}

\subsection{Branching rule}
In this subsection, we prove a combinatorial formula describing the branching rule from $\U$ to $\Ui$, which coincides with that from $\g$ to $\frk$. Before doing so, let us explain that our formula also describes the branching rule from $\frsl_n$ to $\frso_n$. Recall that $\g(I_{\frk})$ denote the complex simple Lie algebra of type $D_m$ if $n \in \Z_{\ev}$ or $B_m$ if $n \in \Z_{\odd}$. Let $e'_i,f'_i,h'_i \in \g(I_{\frk})$, $i \in I_{\frk}$ denote the Chevalley generators. Then, we have seen that there exist isomorphisms $\g(I_{\frk}) \rightarrow \frk$ and $\g(I_{\frk}) \rightarrow \frso_n$ which send $h'_i$ to $w_i$ and $w'_i$, respectively (see Subsection \ref{subsection: symmetric pair of type AI} for the definitions of $w_i$ and $w'_i$). Namely, we have two realizations of $\g(I_{\frk})$ inside $\g = \frsl_n$. By \cite[2.4]{NS05} and the following lemma, we see that the branching rules from $\g$ to $\frk$ and from $\g$ to $\frso_n$ coincide.

\begin{lem}
Let $V := \C^n$ denote the natural representation of $\g = \frsl_n$. Then, as $\g(I_{\frk})$-modules, we have $V|_{\frk} \simeq V|_{\frso_n}$.
\end{lem}

\begin{proof}
It suffices to show that $V|_{\frk}$ and $V|_{\frso_n}$ have the same character. Recall that the subspaces $\C \{ b_{2i-1} \mid i \in [1,m] \}$ and $\C \{ b'_{2i-1} \mid i \in [1,m] \}$ form Cartan subalgebras of $\frk$ and $\frso_n$, respectively. Under the isomorphisms $\g(I_{\frk}) \rightarrow \frk$ and $\g(I_{\frk}) \rightarrow \frso_n$, for each $i \in [1,m]$, $b_{2i-1}$ and $\sqrt{-1}b'_{2i-1}$ correspond to the same vector in $\g(I_{\frk})$. Hence, in order to prove the assertion, it suffices to show that for each $i \in [1,m]$, $b_{2i-1}$ and $\sqrt{-1}b'_{2i-1}$ have the same spectra on $V$.

Let $\{ v_1,\ldots,v_n \}$ denote the standard basis of $V = \C^n$. Then, for each $i \in [1,m]$ such that $2i \leq n$, the vector $v_{2i-1} \pm v_{2i}$ (resp., $v_{2i-1} \pm \sqrt{-1} v_{2i}$) is an engenvector of $b_{2j-1}$ (resp., $\sqrt{-1}b'_{2j-1}$) of eigenvalue $\pm \delta_{i,j}$. Also, if $n \in \Z_{\odd}$, the vector $v_n$ is an eigenvector of both $b_{2j-1}$ and $\sqrt{-1}b'_{2j-1}$ of eigenvalue $0$. Therefore, our claim follows. Thus, the proof completes.
\end{proof}

\begin{rem}\normalfont
Another often used realization of $\g(I_{\frk})$ is
$$
\frk' := \{ X = (x_{i,j})_{1 \leq i,j \leq n} \in \g \mid x_{i,j} = - x_{n-j+1,n-i+1} \},
$$
with $h'_i \in \g(I_{\frk})$ corresponding to $h_i + h_{n-i}$ if $i \neq m$, $h_{m-1}+2h_m+h_{m+1}$ if $i = m$ and $n \in \Z_{\ev}$, or $2(h_m+h_{m+1})$ if $i = m$ and $n \in \Z_{\odd}$. By the same way as above, we see that the branching rule from $\g$ to $\frk'$ is the same as that from $\g$ to $\frk$.
\end{rem}

Let $\lm \in X^+$ and consider the irreducible $\g$-module $V_{\g}(\lm)$ of highest weight $\lm$. As a $\frk$-module, it decomposes as
$$
V_{\g}(\lm) \simeq \bigoplus_{\nu \in X_{\frk,\Int}^+} V_{\frk}(\nu)^{\oplus [\lm:\nu]}
$$
for some $[\lm:\nu] \geq 0$.

Also, consider the irreducible $\U$-module $V(\lm)$, and set $V(\lm)_{\bfK_1} := \U_{\bfK_1} v_\lm$. We have
$$
V(\lm)_1 := V(\lm)_{\bfK_1}/(q-1) V(\lm)_{\bfK_1} \simeq V_{\g}(\lm).
$$
Let $v \in V(\lm)_{\bfK_1}$ and write its $X_{\frk}$-weight vector decomposition as $v = \sum_{\xi \in X_{\frk,\Int}} v_\xi$ with $v_\xi \in V(\lm)_\xi$. Then, for each $\xi \in X_{\frk,\Int}$, we have for sufficiently large $N \geq 0$,
$$
\prod_{i \in I_{\otimes}} \prod_{a \in [-N,N] \setminus \{ \xi_i \}} (B_i - [a]) v = \prod_{i \in I_{\otimes}} \prod_{a \in [-N,N] \setminus \{ \xi_i \}} ([\xi_i]-[a]) v_\xi.
$$
Hence, we obtain
$$
v_\xi = \prod_{i \in I_{\otimes}} \prod_{a \in [-N,N] \setminus \{ \xi_i \}} \frac{1}{[\xi_i]-[a]} (B_i - [a]) v \in \U_{\bfK_1} v_\lm \subseteq V(\lm)_{\bfK_1}.
$$
This shows that
$$
V(\lm)_{\bfK_1} = \bigoplus_{\xi \in X_{\frk,\Int}} (V(\lm)_\xi \cap V(\lm)_{\bfK_1}),
$$
and hence,
$$
\ch_{\Ui} V(\lm) = \ch_{\frk} V_{\g}(\lm).
$$
By character consideration, we have
$$
V(\lm) \simeq \bigoplus_{\nu \in X_{\frk,\Int}^+} V(\nu)^{\oplus [\lm:\nu]},
$$
and hence,
$$
\ol{\clL}(\lm) \simeq \bigoplus_{\nu \in X_{\frk,\Int}^+} \ol{\clL}(\nu)^{\oplus [\lm:\nu]}.
$$
This, together with Theorem \ref{V(nu) is a based Ui-module; general} \eqref{V(nu) is a based Ui-module; general 3}, shows that
$$
[\lm:\nu] = \dim_{\C} \{ b \in \ol{\clL}(\lm)_\nu \mid \Xtil_j b = 0 \Forall j \in \Itil_{\frk} \}.
$$
However, we do not know how $\Xtil_j$ acts on $\ol{\clL}(\lm)$ explicitly. In the sequel, we aim to describe this multiplicity $[\lm;\nu]$ in terms of $\deg_i$ and $\Btil_i$, which in turn, can be described in terms of the crystal structure of $\clB(\lm)$ by Corollary \ref{deg and Btil on usual crystals}.

\begin{lem}\label{Xk-weight decomposition of crystal basis vector}
Let $M$ be a finite-dimensional $\U$-module with a crystal basis $\clB$. Let $b \in \clB$, and write its $X_{\frk}$-weight vector decomposition as $b = \sum_{\nu \in X_{\frk,\Int}} b_\nu$ with $b_\nu \in \ol{\clL}_{M,\nu}$. Then, we have $b_\nu \neq 0$ if and only if $|\nu_{2i-1}| = \deg_{2i-1}(b)$ for all $i \in [1,m]$. Furthermore, if $b_\nu \neq 0$, then we have
$$
b_\nu = \frac{1}{2^s} \prod_{\substack{i \in [1,m] \\ \nu_{2i-1} > 0}} (1+\Btil_{2i-1}) \prod_{\substack{i \in [1,m] \\ \nu_{2i-1} < 0}} (1-\Btil_{2i-1})b,
$$
where $s := \sharp \{ i \mid \nu_{2i-1} \neq 0 \}$.
\end{lem}

\begin{proof}
Let $\nu \in X_{\frk,\Int}$ be such that $|\nu_{2i-1}| = \deg_{2i-1}(b)$ for all $i \in [1,m]$, and set $s$ to be the number of $i$ such that $\nu_{2i-1} \neq 0$. By Lemma \ref{weight decomposition of homogeneous vectors} and Corollary \ref{deg and Btil on usual crystals}, we see that
$$
\frac{1}{2^s} \prod_{\substack{i \in [1,m] \\ \nu_{2i-1} > 0}} (1+\Btil_{2i-1}) \prod_{\substack{i \in [1,m] \\ \nu_{2i-1} < 0}} (1-\Btil_{2i-1})b
$$
is a nonzero $X_{\frk}$-weight vector of weight $\nu$. Also, by an elementary calculation, we have
$$
b = \sum_{\substack{\nu \in X_{\frk,\Int} \\ |\nu_{2i-1}| = \deg_{2i-1}(b) \Forall i \in [1,m]}} \frac{1}{2^s} \prod_{\substack{i \in [1,m] \\ \nu_{2i-1} > 0}} (1+\Btil_{2i-1}) \prod_{\substack{i \in [1,m] \\ \nu_{2i-1} < 0}} (1-\Btil_{2i-1})b.
$$
Therefore, this is the $X_{\frk}$-weight decomposition of $b$. Thus, the proof completes.
\end{proof}

For a finite-dimensional classical weight module $M$ equipped with a contragredient Hermitian inner product, and $b \in \ol{\clL}_M$, consider the following condition:
\begin{align}\label{combinatorial characterization}
\begin{split}
&\text{$b$ is homogeneous, i.e., $B_i$-homogeneous for all $i \in [1,n-1]$}, \\
&\deg_{2i}(b) = 0 \qu \Forall i \in [1,m] \text{ such that } 2i < n, \\
&\deg_{2i+1}((\Btil_{2i} \Btil_{2i-1})^{\deg_{2i+1}(b)} b) = 0 \qu \Forall i \in [1,m] \text{ such that } 2i+1 < n.
\end{split}
\end{align}

\begin{prop}\label{preparation for branching rule; general}
Let $M$ be a finite-dimensional classical weight module equipped with a contragredient Hermitian inner product. Set
\begin{align}
\begin{split}
&L_1 := \{ b \in \ol{\clL}_M \mid \Xtil_j b = 0 \Forall j \in \Itil_{\frk} \}, \\
&L_2 := \C\{ b \in \ol{\clL}_M \mid \text{ $b$ satisfies condition \eqref{combinatorial characterization}} \}.
\end{split} \nonumber
\end{align}
Then, the linear map
$$
L_2 \rightarrow L_1;\ b \mapsto (1 + \Btil_1)(1 + \Btil_3) \cdots (1 + \Btil_{2m'-1}) b,
$$
where
$$
m' := \begin{cases}
m-1 \qu & \IF n \in \Z_{\ev}, \\
m \qu & \IF n \in \Z_{\odd}
\end{cases}
$$
is an isomorphism of $\C$-vector spaces with inverse
$$
L_1 \rightarrow L_2;\ b \mapsto \sum_{\substack{\zeta \in X^\imath \\ \zeta_{2i} = 0 \Forall i \in [1,m']}} \mathbf{1}_\zeta b.
$$
\end{prop}

\begin{proof}
We prove by induction on $m'$. The $m' = 1$ case follows from Propositions \ref{preparation for branching rule n=3} and \ref{preparation for branching rule n=4}.

As the proof of Proposition \ref{preparation for branching rule n=3}, it suffices to consider the case when $M = V(\nu)$ for some $\nu \in X_{\frk,\Int}^+$. Clearly, we have
$$
L_1 = \C b_\nu.
$$

As the proof of Lemma \ref{V(nu) is standard if it is contained in a standard module}, we have
\begin{align}
\begin{split}
c^\nu &:= \sum_{\substack{\zeta \in X^\imath \\ \zeta_{2i} = 0 \Forall i \in [1,m']}} \mathbf{1}_\zeta b_\nu \\
&= \ev_\infty(\sum_{\substack{\zeta \in X^\imath \\ \zeta_{2i} = 0 \Forall i \in [1,m']}} \mathbf{1}_\zeta v_\nu) \\
&= \ev_\infty(\mathbf{1}_{1,\ol{\nu_1}} \mathbf{1}_{2,\ol{0}} \mathbf{1}_{3,\ol{\nu_3}} \mathbf{1}_{4,\ol{0}} \cdots \mathbf{1}_{2m'-1,\ol{\nu_{2m'-1}}} \mathbf{1}_{2m',\ol{0}} v_\nu) \\
&= \ev_\infty(\frac{1}{2^{m'}}(1 + (-1)^{\nu_1}Y_2^{(2\nu_1)})(1 +(-1)^{\nu_3} Y_4^{(2\nu_3)}) \cdots (1 +(-1)^{\nu_{2m'-1}} Y_{2m'}^{(2\nu_{2m'-1})}) v_\nu) \\
&= \frac{1}{2^{m'}}(1 + (-1)^{\nu_1}\Ytil_2^{2\nu_1})(1 +(-1)^{\nu_3} \Ytil_4^{2\nu_3}) \cdots (1 +(-1)^{\nu_{2m'-1}} \Ytil_{2m'}^{2\nu_{2m'-1}}) b_\nu \neq 0.
\end{split} \nonumber
\end{align}
We show that $c^\nu \in L_2$. As a $\Ui_{1,2,\ldots,2m'-1}$-module vector, both $v_\nu$ and $Y_{2m'}^{(2\nu_{2m'-1})} v_\nu$ are highest weight vectors of weights $(\nu_1,\nu_3,\ldots,\nu_{2m'-3},\nu_{2m'-1})$ and $(\nu_1,\nu_3,\ldots,\nu_{2m'-3},-\nu_{2m'-1})$, respectively. By the induction hypothesis, we see that $c^\nu$ satisfies condition \eqref{combinatorial characterization} as a $\Ui_{1,2,\ldots,2m'-1}$-module vector.

On the other hand, by the definition of $\mathbf{1}_{i,p}$'s ($i \in I$, $p \in \{ \ev,\odd \}$), they pairwise commute. Hence, we have
\begin{align}
\begin{split}
&\mathbf{1}_{1,\ol{\nu_1}} \mathbf{1}_{2,\ol{0}} \mathbf{1}_{3,\ol{\nu_3}} \mathbf{1}_{4,\ol{0}} \cdots \mathbf{1}_{2m'-1,\ol{\nu_{2m'-1}}} \mathbf{1}_{2m',\ol{0}} v_\nu \\
&= \mathbf{1}_{2m'-1,\ol{\nu_{2m'-1}}} \mathbf{1}_{2m',\ol{0}} \mathbf{1}_{1,\ol{\nu_1}} \mathbf{1}_{2,\ol{0}} \mathbf{1}_{3,\ol{\nu_3}} \mathbf{1}_{4,\ol{0}} \cdots \mathbf{1}_{2m'-3,\ol{\nu_{2m'-3}}} \mathbf{1}_{2m'-2,\ol{0}} v_\nu \\
&= \frac{1}{2^{m'-1}} \mathbf{1}_{2m'-1,\ol{\nu_{2m'-1}}} \mathbf{1}_{2m',\ol{0}} \\
&\quad \cdot (1 + (-1)^{\nu_1}Y_2^{(2\nu_1)})(1 +(-1)^{\nu_3} Y_4^{(2\nu_3)}) \cdots (1 +(-1)^{\nu_{2m'-3}} Y_{2m'-2}^{(2\nu_{2m'-3})}) v_\nu,
\end{split} \nonumber
\end{align}
and consequently,
\begin{align*}
\begin{split}
  c^\nu &= \frac{1}{2^{m'-1}} \mathbf{1}_{2m'-1,\ol{\nu_{2m'-1}}} \mathbf{1}_{2m',\ol{0}}\\
  &\quad \cdot (1 + (-1)^{\nu_1}\Ytil_2^{2\nu_1})(1 +(-1)^{\nu_3} \Ytil_4^{2\nu_3}) \cdots (1 +(-1)^{\nu_{2m'-3}} \Ytil_{2m'-2}^{2\nu_{2m'-3}}) b_\nu.
\end{split}
\end{align*}
Now, set
$$
J := \begin{cases}
\{ 2m'-1,2m',2m'+1 \} \qu & \IF n \in \Z_{\ev}, \\
\{ 2m'-1,2m' \} \qu & \IF n \in \Z_{\odd}.
\end{cases}
$$
By weight consideration, as a $\Ui_J$-module vector,
$$
(1 + (-1)^{\nu_1}Y_2^{(2\nu_1)})(1 +(-1)^{\nu_3} Y_4^{(2\nu_3)}) \cdots (1 +(-1)^{\nu_{2m'-3}} Y_{2m'-2}^{(2\nu_{2m'-3})}) v_\nu
$$
is a sum of highest weight vectors of weight $(\nu_{2m'-1},\nu_{2m'+1})$ and a lowest weight vector of weight $(-\nu_{2m'-1},\nu_{2m'+1})$ if $n \in \Z_{\ev}$, or a sum of highest weight vectors of weight $\nu_{2m'-1}$ and lowest weight vectors of weight $-\nu_{2m'-1}$ if $n \in \Z_{\odd}$. Here, lowest weight vector means a vector of the form $Y^{(2\nu_{2m'-1})} v$ for some highest weight vector $v$. Then, by Proposition \ref{Xi weight decomposition of highest weight vector}, we see that $c^\nu$ is a linear combination of vectors of the form
$$
(1 +(-1)^{\nu_{2m'-1}} \Ytil) b'
$$
with $b'$ being a highest weight vector of weight $(\nu_{2m'-1},\pm \nu_{2m'+1})$ if $n \in \Z_{\ev}$, or $\nu_{2m'-1}$ if $n \in \Z_{\odd}$. Therefore, by the $n = 3,4$ cases, $c^\nu$ satisfies condition \eqref{combinatorial characterization} as a $\Ui_J$-module vector. Thus, we conclude that
$$
c^\nu \in L_2,
$$
as desired.

Next, let $c \in \ol{\clL}(\nu)$ satisfy condition \eqref{combinatorial characterization}. As the proof of Lemma \ref{Xk-weight decomposition of crystal basis vector}, we can write the $X_{\frk}$-weight vector decomposition of $c$ as
$$
c = \sum_{\substack{\xi \in X_{\frk,\Int} \\ |\xi_{2k-1}| = \deg_{2k-1}(c) \Forall k \in [1,m]}} c_\xi, \qu c_\xi \in \ol{\clL}(\nu)_\xi.
$$
Let $i \in [1,m']$. Then, we have
$$
(1 + \Btil_{2i-1})c = \sum_{\substack{|\xi_{2k-1}| = \deg_{2k-1}(c) \Forall k \in [1,m] \\ \xi_{2i-1} \geq 0}} 2^{1-\delta_{\deg_{2i-1}(c),0}} c_\xi.
$$
Define $J_i \subseteq I$ by
$$
J_i := \begin{cases}
\{ 2i-1,2i,2i+1 \} \qu & \IF 2i+1 < n ,\\
\{ 2i-1,2i \} \qu & \IF 2i+1 = n.
\end{cases}
$$
Then, by the $n = 3,4$ cases, we have
$$
\Xtil_j (1 + \Btil_{2i-1}) c = 0 \qu \Forall j \in \Itil_{\frk}(J_i).
$$
By weight consideration, we see that
\begin{align}\label{each weight vector is highest}
\Xtil_j c_\xi = 0 \qu \Forall j \in \Itil_{\frk}(J_i) \AND \xi \in X_{\frk,\Int} \text{ such that } \xi_{2i-1} = \deg_{2i-1}(c).
\end{align}

Now, remark that
$$
c' := (1 + \Btil_1)(1 + \Btil_3) \cdots (1 + \Btil_{2m'-1}) c = \sum_{\substack{|\xi_{2k-1}| = \deg_{2k-1}(c) \Forall k \in [1,m] \\ \xi_{2i-1} \geq 0 \Forall i \in [1,m']}} 2^s c_\xi,
$$
where $s = \sharp\{ i \in [1,m'] \mid \deg_{2i-1}(c) \neq 0 \}$. By \eqref{each weight vector is highest}, we conclude that
$$
\Xtil_j c' = 0 \Forall j \in \Itil_{\frk},
$$
and hence,
$$
c' \in \C b_\nu = L_1
$$
by Theorem \ref{V(nu) is a based Ui-module; general} \eqref{V(nu) is a based Ui-module; general 3}.

This far, we have verified that the linear maps in consideration are well-defined, and that $\dim L_2 \geq \dim L_1 = 1$. Hence, in order to complete the proof, it suffices to show that
$$
\sum_{\substack{\zeta \in X^\imath \\ \zeta_{2i} = 0 \Forall i \in [1,m']}} \mathbf{1}_\zeta (1 + \Btil_1)(1 + \Btil_3) \cdots (1 + \Btil_{2m'-1}) c = c.
$$
Set
$$
J' := \{ 1,2,3 \}, \qu J'' := \{ 3,4,\ldots,n-1 \}.
$$
It is easily verified that $(1+\Btil_1) c$ satisfies condition \eqref{combinatorial characterization} as a $\Ui_{J''}$-module vector. Then, by the induction hypothesis, we have
$$
\sum_{\substack{\zeta'' \in X^\imath(J'') \\ \zeta_{2i} = 0 \Forall i \in [2,m']}} \mathbf{1}_{J'',\zeta''} (1+\Btil_3)(1+\Btil_5) \cdots (1+\Btil_{2m'-1})(1+\Btil_1)c = (1+\Btil_1)c.
$$
Hence, we compute as
\begin{align}
\begin{split}
&\sum_{\substack{\zeta \in X^\imath \\ \zeta_{2i} = 0 \Forall i \in [1,m']}} \mathbf{1}_\zeta (1 + \Btil_1)(1 + \Btil_3) \cdots (1 + \Btil_{2m'-1}) c \\
&\qu = \sum_{\substack{\zeta' \in X^\imath(J') \\ \zeta_{2} = 0}} \sum_{\substack{\zeta'' \in X^\imath(J'') \\ \zeta_{2i} = 0 \Forall i \in [2,m']}} \mathbf{1}_{J',\zeta'} \mathbf{1}_{J'',\zeta''} (1 + \Btil_3)(1+\Btil_5) \cdots (1 + \Btil_{2m'-1})(1+\Btil_1) c \\
&\qu = \sum_{\substack{\zeta' \in X^\imath(J') \\ \zeta_{2} = 0}} \mathbf{1}_{J',\zeta'} (1+\Btil_1)c.
\end{split} \nonumber
\end{align}
By the $n = 4$ case, the last term equals $c$, as desired. Thus, the proof completes.
\end{proof}

\begin{lem}
Let $\lm \in X^+$. Let $b \in \ol{\clL}(\lm)$ and write $b = \sum_{k=1}^r a_k b_k$ with $a_k \in \C^\times$, $b_k \in \clB(\lm)$. If $b$ satisfies condition \eqref{combinatorial characterization}, then so does $b_k$ for all $k \in [1,r]$.
\end{lem}

\begin{proof}
Since each element in $\clB(\lm)$ is homogeneous, we must have $\deg_{i}(b_k) = \deg_i(b)$ for all $i \in [1,n-1]$. In particular, we obtain
$$
\deg_{2i}(b_k) = 0 \Forall i \in [1,m'].
$$

Next, let $i \in I_{\frk}$ be such that $2i+1 < n$, and set $\{ k_1,\ldots,k_l \} \subseteq [1,r]$ to be the subset such that $(\Btil_{2i} \Btil_{2i-1})^{\deg_{2i+1}(b)} b_{k_j} \neq 0$ and $\deg_{2i+1}((\Btil_{2i} \Btil_{2i-1})^{\deg_{2i+1}(b)} b_{k_j}) \neq 0$. Since $\deg_{2i+1}((\Btil_{2i} \Btil_{2i-1})^{\deg_{2i+1}(b)} b) = 0$, we must have
$$
\sum_{j=1}^l a_{k_j} (\Btil_{2i} \Btil_{2i-1})^{\deg_{2i+1}(b)} b_{k_j}  = 0.
$$
Since $(\Btil_{2i} \Btil_{2i-1})^{\deg_{2i+1}(b)} b_{k_j}$ are distinct vectors in $\clB(\lm)$, we see that $(\Btil_{2i} \Btil_{2i-1})^{\deg_{2i+1}(b)} b_{k_j} = 0$ for all $j \in [1,l]$. Therefore, we conclude that $\deg_{2i+1}((\Btil_{2i} \Btil_{2i-1})^{\deg_{2i+1}(b)} b_k) = 0$ for all $k \in [1,r]$. Thus, the proof completes.
\end{proof}

\begin{theo}\label{Branching rule from U to Ui}
Let $\lm \in X^+$ and $\nu \in X_{\frk,\Int}^+$.
\begin{enumerate}
\item Suppose $n \in \Z_{\ev}$ and $\nu_{2m-1} \neq 0$. Then, we have
\begin{align}
\begin{split}
[\lm:\nu] = \hf \sharp \{ b \in \clB(\lm) \mid &\deg_{2i-1}(b) = |\nu_{2i-1}| \Forall i \in [1,m], \\
&\deg_{2i}(b) = 0 \Forall i \in [1,m-1], \\
&\deg_{2i+1}((\Btil_{2i}\Btil_{2i-1})^{|\nu_{2i+1}|} b) = 0 \Forall i \in [1,m-1] \}.
\end{split} \nonumber
\end{align}
\item Suppose $n \in \Z_{\ev}$ and $\nu_{2m-1} = 0$. Then, we have
\begin{align}
\begin{split}
[\lm:\nu] = \sharp \{ b \in \clB(\lm) \mid &\deg_{2i-1}(b) = |\nu_{2i-1}| \Forall i \in [1,m], \\
&\deg_{2i}(b) = 0 \Forall i \in [1,m-1], \\
&\deg_{2i+1}((\Btil_{2i}\Btil_{2i-1})^{\nu_{2i+1}} b) = 0 \Forall i \in [1,m-1] \}.
\end{split} \nonumber
\end{align}
\item Suppose $n \in \Z_{\odd}$. Then, we have
\begin{align}
\begin{split}
[\lm:\nu] = \sharp \{ b \in \clB(\lm) \mid &\deg_{2i-1}(b) = |\nu_{2i-1}| \Forall i \in [1,m], \\
&\deg_{2i}(b) = 0 \Forall i \in [1,m], \\
&\deg_{2i+1}((\Btil_{2i}\Btil_{2i-1})^{\nu_{2i+1}} b) = 0 \Forall i \in [1,m-1] \}.
\end{split} \nonumber
\end{align}
\end{enumerate}
\end{theo}

\begin{proof}
Let $b \in \clB(\lm)$ be such that
$$
|\deg_1(b)| \geq |\deg_3(b)| \geq \cdots \geq |\deg_{2m-1}(b)|.
$$
By Lemma \ref{Xk-weight decomposition of crystal basis vector}, we can write
$$
b = \sum_{\substack{\xi \in X_{\frk,\Int} \\ |\xi_{2i-1}| = \deg_{2i-1}(b) \Forall i \in [1,m]}} b_\xi, \qu b_\xi \in \ol{\clL}(\lm)_\xi \setminus \{0\}.
$$
Note that when $n \in \Z_{\ev}$ and $\deg_{2m-1}(b) \neq 0$, there are exactly two dominant weights $\xi^+,\xi^- \in X_{\frk,\Int}^+$ such that $\xi^+_{2m-1} > 0$, $\xi^-_{2m-1} < 0$, and $b_{\xi^+},b_{\xi^-} \neq 0$. Otherwise, there is exactly one dominant weight $\xi \in X_{\frk,\Int}^+$ such that $b_\xi \neq 0$.

From now on, we concentrate on the case when $n \in \Z_{\ev}$ and $\nu_{2m-1} \neq 0$. The other cases can be proved in a similar way more easily.

Set
$$
L' := \{ b \in \ol{\clL}(\lm)_{\nu^+} \oplus \ol{\clL}(\lm)_{\nu^-} \mid \Xtil_j b = 0 \Forall j \in \Itil_{\frk} \},
$$
and
\begin{align}
\begin{split}
B := \{ b \in \clB(\lm) \mid &\deg_{2i-1}(b) = |\nu_{2i-1}| \Forall i \in [1,m], \\
&\deg_{2i}(b) = 0 \Forall i \in [1,m-1], \\
&\deg_{2i+1}((\Btil_{2i}\Btil_{2i-1})^{|\nu_{2i+1}|} b) = 0 \Forall i \in [1,m-1] \}.
\end{split} \nonumber
\end{align}
By Theorem \ref{preparation for branching rule; general}, we have an isomorphism
$$
\C B \rightarrow L'
$$
which sends
$$
b \mapsto (1 + \Btil_1)(1 + \Btil_3) \cdots (1 + \Btil_{2m-3}) b,
$$
with inverse given by
$$
b \mapsto \sum_{\substack{\zeta \in X^\imath \\ \zeta_{2i} = \ol{0} \Forall i \in [1,m-1]}} \mathbf{1}_\zeta b.
$$

Let $b \in B$ and write its $X_{\frk}$-weight vector decomposition as
$$
b = \sum_{\substack{\xi \in X_{\frk,\Int} \\ |\xi_{2i-1}| = \deg_{2i-1}(b) \Forall i \in [1,m]}} b_\xi, \qu b_\xi \in \ol{\clL}(\lm)_\xi \setminus \{0\}.
$$
Set
$$
b^+ := \sum_{\xi_1 \cdot \xi_3 \cdots \xi_{2m-1} > 0} b_\xi, \qu b^- := \sum_{\xi_1 \cdot \xi_3 \cdots \xi_{2m-1} < 0} b_\xi.
$$
Let us show that $b^\pm \in \C B$.

By the definitions, we have
$$
b = b^+ + b^-,
$$
and
$$
(1+\Btil_1)(1+\Btil_3) \cdots (1+\Btil_{2m-3}) b^\pm = 2^{m-1} b_{\nu^\pm}.
$$
Since we have $\Xtil_j (b_{\nu^+} + b_{\nu^-}) = 0$ for all $j \in \Itil_{\frk}$, by weight consideration, we obtain
$$
\Xtil_j b_{\nu^\pm} = 0 \Forall j \in \Itil_{\frk}.
$$
In particular, we have
$$
b_{\nu^\pm} \in L'.
$$

On the other hand, we have
$$
\sum_{\substack{\zeta \in X^\imath \\ \zeta_{2i} = \ol{0} \Forall i \in [1,m-1]}} \mathbf{1}_\zeta 2^{m-1}(b_{\nu^+} + b_{\nu^-}) = b = b^+ + b^-.
$$
Also, since $b_{\nu^\pm}$ is a highest weight vector of weight $\nu^\pm$, we see that
\begin{align*}
  \begin{split}
    \sum_{\substack{\zeta \in X^\imath \\ \zeta_{2i} = \ol{0} \Forall i \in [1,m-1]}} &\mathbf{1}_\zeta 2^{m-1} b_{\nu^\pm} \\
    &= (1+(-1)^{\nu_1}\Ytil_2^{2\nu_1})(1+(-1)^{\nu_3}\Ytil_4^{2\nu_3}) \cdots (1+(-1)^{\nu_{2m-3}}\Ytil_{2m-2}^{2\nu_{2m-3}}) b_{\nu^\pm}
  \end{split}
\end{align*}
is a sum of $X_{\frk}$-weight vectors of weights $\xi' \in X_{\frk,\Int}$ with $\pm \xi'_1 \cdot \xi'_3 \cdots \xi'_{2m-1} > 0$. Therefore, we see that
$$
\sum_{\substack{\zeta \in X^\imath \\ \zeta_{2i} = \ol{0} \Forall i \in [1,m-1]}} \mathbf{1}_\zeta 2^{m-1} b_{\nu^\pm} = b^\pm.
$$
This shows that
$$
b^\pm \in \C B,
$$
as desired.

Now, we show that $\Btil_1 \Btil_3 \cdots \Btil_{2m-1}$ restricts to an involution on $B$. In fact, for each $b = b^+ + b^- \in B$, we have
$$
\Btil_1 \Btil_3 \cdots \Btil_{2m-1}b = b^+ - b^- \in \C B.
$$
Since $\Btil_1 \Btil_3 \cdots \Btil_{2m-1} b \in \clB(\lm)$, we conclude that $\Btil_1 \Btil_3 \cdots \Btil_{2m-1} b \in B$.

Also, we have $\Btil_1 \Btil_3 \cdots \Btil_{2m-1} b = b$ if and only if $b^- = 0$, which cannot occur by Lemma \ref{Xk-weight decomposition of crystal basis vector}. Thus, we obtain an isomorphism
$$
\C \{ (1\pm \Btil_1 \Btil_3 \cdots \Btil_{2m-1})b \mid b \in B \} \rightarrow L'_\pm := \{ b \in \ol{\clL}(\lm)_{\nu^\pm} \mid \Xtil_j b = 0 \Forall j \in \Itil_{\frk} \}.
$$
This implies
$$
[\lm:\nu^\pm] = \dim L'_\pm = \hf |B|.
$$
Since we have either $\nu = \nu^+$ or $\nu = \nu^-$, the assertion follows.
\end{proof}


\begin{thebibliography}{99}
\bibitem{AV20} A. Appel and B. Vlaar, Universal $K$-matrices for quantum Kac-Moody algebras, Represent. Theory 26 (2022), 764--824.

\bibitem{BK15} M. Balagovi\'{c} and S. Kolb, The bar involution for quantum symmetric pairs, Represent. Theory 19 (2015), 186--210.

\bibitem{BK19} M. Balagovi\'{c} and S. Kolb, Universal $K$-matrix for quantum symmetric pairs, J. Reine Angew. Math. 747 (2019), 299--353.

\bibitem{BW18a} H. Bao and W. Wang, A New Approach to Kazhdan-Lusztig Theory of Type B via Quantum Symmetric Pairs, Ast\'{e}risque 2018, no. 402, vii+134 pp.

\bibitem{BW18b} H. Bao and W. Wang, Canonical bases arising from quantum symmetric pairs, Invent. Math. 213 (2018), no. 3, 1099--1177.

\bibitem{BW18c} H. Bao and W. Wang, Canonical bases arising from quantum symmetric pairs of Kac-Moody type, Compos. Math. 157 (2021), no. 7 1507--1537.

\bibitem{BeW18} C. Berman and W. Wang, Formulae of $\imath$-divided powers in $U_q(\frsl_2)$, J. Pure Appl. Algebra 222 (2018), no. 9, 2667--2702.

\bibitem{DCM20} K. De Commer and M. Matassa, Quantum flag manifolds, quantum symmetric spaces and their associated universal K-matrices, Adv. Math. 366 (2020), 107029, 100 pp.

\bibitem{D96} J. Dixmier, Enveloping Algebras, Revised reprint of the 1977 translation. Graduate Studies in Mathematics, 11. American Mathematical Society, Providence, RI, 1996. xx+379 pp.

\bibitem{ES18} M. Ehrig and C. Stroppel, Nazarov-Wenzl algebras, coideal subalgebras and categorified skew Howe duality, Adv. Math. 331 (2018), 58--142.


\bibitem{GK91} A. M. Gavrilik and A. U. Klimyk, $q$-deformed orthogonal and pseudo-orthogonal algebras and their representations, Lett. Math. Phys. 21 (1991), no. 3, 215--220. 

\bibitem{HK02} J. Hong and S.-J. Kang, Introduction to Quantum Groups and Crystal Bases, Graduate Studies in Mathematics, 42. American Mathematical Society, Providence, RI, 2002. xviii+307 pp.

\bibitem{H72} J. E. Humphreys, Introduction to Lie Algebras and Representation Theory, Graduate Texts in Mathematics, Vol. 9. Springer-Verlag, New York-Berlin, 1972. xii+169 pp.

\bibitem{JK21} I.-S. Jang and J.-H. Kwon, Flagged Littlewood-Richardson tableaux and branching rule for classical groups, J. Combin. Theory Ser. A 181 (2021), 105419.

\bibitem{J96} J. C. Jantzen, Lectures on Quantum Groups, Graduate Studies in Mathematics, 6. American Mathematical Society, Providence, RI, 1996. viii+266 pp.

\bibitem{Ka90} M. Kashiwara, Crystalizing the $q$-analogue of universal enveloping algebras, Comm. Math. Phys. 133 (1990), no. 2, 249--260.

\bibitem{Ka91} M. Kashiwara, On crystal bases of the $q$-analogue of universal enveloping algebras, Duke Math. J. 63 (1991), no. 2, 465--516.

\bibitem{KL79} D. Kazhdan and G. Lusztig, Representations of Coxeter groups and Hecke algebras, Invent. Math. 53 (1979), no. 2, 165--184.

\bibitem{Ko14} S. Kolb, Quantum symmetric Kac-Moody pairs, Adv. Math. 267 (2014), 395--469. 


\bibitem{Ko90} T. H. Koornwinder, Orthogonal polynomials in connection with quantum groups, Orthogonal polynomials (Columbus, OH, 1989), 257--292, NATO Adv. Sci. Inst. Ser. C Math. Phys. Sci., 294, Kluwer Acad. Publ., Dordrecht, 1990. 

\bibitem{Le99} G. Letzter, Symmetric pairs for quantized enveloping algebras, J. Algebra 220 (1999), no. 2, 729--767. 

\bibitem{L90} G. Lusztig, Canonical bases arising from quantized enveloping algebras, J. Amer. Math. Soc. 3 (1990), no. 2, 447--498.

\bibitem{L10} G. Lusztig, Introduction to Quantum Groups, Reprint of the 1994 edition. Modern Birkh\"{a}user Classics. Birkh\"{a}user/Springer, New York, 2010. xiv+346 pp.

\bibitem{NS05} S. Naito and D. Sagaki, An approach to the branching rule from sl2n(C) to sp2n(C) via Littelmann's path model, J. Algebra 286 (2005), no. 1, 187--212.

\bibitem{N96} M. Noumi, Macdonald's symmetric polynomials as zonal spherical functions on some quantum homogeneous spaces, Adv. Math. 123 (1996), no. 1, 16--77.

\bibitem{RV20} V. Regelskis and B. Vlaar, Quasitriangular coideal subalgebras of $U_q(\mathfrak{g})$ in terms of generalized Satake diagrams, Bull. Lond. Math. Soc. 52 (2020), no. 4, 693--715.

\bibitem{ST19} A. Sartori and D. Tubbenhauer, Webs and $q$-Howe dualities in types $BCD$, Trans. Amer. Math. Soc. 371 (2019), no. 10, 7387--7431.

\bibitem{St20} J. V. Stokman, Generalized Onsager algebras, Algebr. Represent. Theory 23 (2020), no. 4, 1523--1541.

\bibitem{Wa17} H. Watanabe, Crystal basis theory for a quantum symmetric pair $(\U,\U^\jmath)$, Int. Math. Res. Not. IMRN 2020, no. 22, 8292--8352.

\bibitem{Wa21} H. Watanabe, Global crystal bases for integrable modules over a quantum symmetric pair of type AIII, Represent. Theory 25 (2021), 27--66.

\bibitem{W19} H. Watanabe, Classical weight modules over $\imath$quantum groups, J. Algebra 578 (2021), 241--302.
\end{thebibliography}
\end{document}